\documentclass[a4paper,11pt]{article}

\usepackage[utf8]{inputenc}
\usepackage{amsfonts}
\usepackage{amssymb}
\usepackage{amsthm}
\usepackage{amsmath}
\usepackage{a4wide}
\usepackage{mathrsfs}
\usepackage{epsfig}
\usepackage{palatino}
\usepackage{tikz}
\usepackage{fp,ifthen}
\usepackage{nicefrac}
\usetikzlibrary{decorations.pathreplacing}
\usepackage{float}
\usepackage{enumerate} %Enumerate package must be always written before enumitem
\usepackage{enumitem} %Customize enumerate and itemize
\usepackage[a4paper,twoside,top=2.5cm, bottom=2cm, left=1.7cm, right=1.7cm]{geometry}

\usepackage[colorinlistoftodos,textwidth=2.3cm]{todonotes} %Package for comments, WARNING: it does not compile if you write Italian stresses like è,ì,à,ò,.. inside the command \todo{...}

\newcommand{\pdfgraphics}{\ifpdf\DeclareGraphicsExtensions{.pdf,.jpg}\else\fi}
\usepackage{graphicx}

\usepackage{color}
\definecolor{hanblue}{rgb}{0.27, 0.42, 0.81}
\definecolor{red}{rgb}{1.0, 0.0, 0.0}
\usepackage[colorlinks, citecolor=blue,linkcolor=blue, urlcolor = blue]{hyperref}
\usepackage[english,capitalize]{cleveref}

\usepackage{mathtools}
\theoremstyle{plain}

\newtheorem{teo}{Theorem}[section]
\newtheorem{lemma}[teo]{Lemma}
\newtheorem{prop}[teo]{Proposition}
\newtheorem{cor}[teo]{Corollary}

\theoremstyle{definition}
\newtheorem{defn}[teo]{Definition}

\theoremstyle{remark}
\newtheorem{rem}[teo]{Remark}

\numberwithin{equation}{section}

\newcommand{\de}{\ensuremath{\,\mathrm d}} % Il de degli integrali (contiene lo spazio da mettere tra integranda e misura)
\renewcommand{\d}{\ensuremath{\mathrm d}} % Il de delle forme
\newcommand{\st}{\ensuremath{\ :\ }} % Such that in formulas
\newcommand{\eqdef}{\ensuremath{\vcentcolon=}}
\newcommand \eps{\ensuremath{\varepsilon}} 
\renewcommand{\epsilon}{\varepsilon}
\newcommand{\N}{\ensuremath{\mathbb N}}%Natural numbers
\newcommand{\R}{\ensuremath{\mathbb R}}%%Real numbers
\newcommand{\TT}{\mathsf{T}}
\newcommand{\NN}{\mathsf{N}}
\DeclarePairedDelimiter\scal{\langle}{\rangle} %Scalar product
\newcommand{\Ell}{\mathscr{L}}

\makeatletter
\renewcommand*\env@matrix[1][*\c@MaxMatrixCols c]{%
  \hskip -\arraycolsep
  \let\@ifnextchar\new@ifnextchar
  \array{#1}}
\makeatother

\begin{document}

\pdfgraphics % Use this command right after \begin{document}

\title{\L ojasiewicz--Simon inequalities for minimal networks: stability and convergence}

\author{Alessandra Pluda \footnote{\href{mailto:alessandra.pluda@unipi.it}{alessandra.pluda@unipi.it}, Dipartimento di Matematica, Universit\`a di Pisa, Largo Bruno Pontecorvo 5, 56127 Pisa, Italy.} \and Marco Pozzetta \footnote{\href{mailtto:marco.pozzetta@unina.it}{marco.pozzetta@unina.it}, Dipartimento di Matematica e Applicazioni, Universit\`a di Napoli Federico II, Via Cintia, Monte S. Angelo 80126 Napoli, Italy}}

%\author{Alessandra Pluda}
%\address{Dipartimento di Matematica, Universit\`a di Pisa, Largo Bruno Pontecorvo 5, 56127 Pisa, Italy.}
%\email{alessandra.pluda@unipi.it}

%\author{Marco Pozzetta}
%\address{Dipartimento di Matematica e Applicazioni, Universit\`a di Napoli Federico II, Via Cintia, Monte S. Angelo 80126 Napoli, Italy.}
%\email{marco.pozzetta@unina.it}

%\subjclass{...}
%\keywords{...}

\date{\today}

\maketitle

\vspace{-0.5cm}
\begin{abstract}
\noindent We investigate stability properties of the motion by curvature of planar networks.
We prove \L ojasiewicz--Simon gradient inequalities for the length functional of planar networks with triple junctions. In particular, such an inequality holds for networks with junctions forming angles equal to $\tfrac23\pi$ that are close in $H^2$-norm to minimal networks, i.e., networks whose edges also have vanishing curvature. The latter inequality bounds a concave power of the difference between length of a minimal network $\Gamma_*$ and length of a triple junctions network $\Gamma$ from above by the $L^2$-norm of the curvature of the edges of $\Gamma$.\\
We apply this result to prove the stability of minimal networks in the sense that a motion by curvature starting from a network sufficiently close in $H^2$-norm to a minimal one exists for all times and smoothly converges.\\
We further rigorously construct an example of a motion by curvature having uniformly bounded curvature that smoothly converges to a degenerate network in infinite time.
\end{abstract}


% {\color{red}
% \vspace{1cm}
% Correzioni a typo fatte in questo file rispetto alla PRIMA versione arxiv:

% \begin{itemize}
%     \item Corretto l'altezza del rettangolo dell'esempio della Sez 6 da $1$ a $2/\sqrt{3}$. Cambiato quindi in: \cref{Fig:Network}, \cref{rem:FamigliaEsempi}, inizio testo sez 6, \cref{Fig:Symbols}, \cref{Fig:Graph}.
    
%     \item Corretto la derivata $\partial_x u(t,1)$ nello \ref{Step3} e coerentemente quello che c'è dopo.
    
%     \item Nella \cref{lem:Fredholmness}, tolto l'operatore ${\rm I}:Z\to Z^\star$ nella definizione di alcuni operatori ausiliari. Era mal posto nel senso che lì ${\rm I}$ non prendeva in ingresso elementi di $Z$.
    
% \end{itemize}
% }

\tableofcontents

\section{Introduction}

A \emph{planar network} is a pair $(\Gamma,G)$, where $G$ is an abstract connected graph with edges homeomorphic to the interval $[0,1]$ and $\Gamma:G\to\R^2$ is a continuous map, see \cref{network}. We shall mostly consider \emph{triple junctions networks}, that are networks such that the edges of $G$ either meet in junctions of order three, or end on terminal points of the graph, and the restriction $\gamma^i$ of $\Gamma$ to each edge is a $C^1$-embedding, see \cref{def:TripleJunctionsRegularMinimal}. 
If we further require that embedded edges meet forming angles equal to $\tfrac23\pi$, the network is said to be \emph{regular}, and if also such embeddings are straight segments, it is said to be \emph{minimal}, see \cref{def:TripleJunctionsRegularMinimal}.

Minimal networks are easily seen to be critical points of the length functional ${\rm L}$, which is defined by the sum of the lengths of the embedded edges via the parametrization $\Gamma$ of a network. This class includes Steiner trees of finitely many points in the plane, that is, networks minimizing the length among those connecting such given points \cite{PaoliniSteiner}.

In this paper we investigate functional analytic and stability properties of the length functional and of the $L^2$-gradient flow of ${\rm L}$.

The first of our main results consists in proving \emph{\L ojasiewicz--Simon gradient inequalities} for the length functional in $H^2$-neighborhoods of minimal networks.

\begin{teo}[cf. \cref{cor:LojaNetworkMinimali}]\label{thm:LojaIntro}
Let $\Gamma_*:G\to\R^2$ be a minimal network. Then there exist $C_{\rm LS},\sigma >0$ and $\theta \in (0,\tfrac12]$ such that the following holds.

If $\Gamma:G\to\R^2$ is a regular network of class $H^2$ such that $\Gamma$ and $\Gamma_*$ have the same endpoints and
\begin{align*}
     \sum_i \| \gamma^i_* - \gamma^i \|_{H^2(\d x)} \le \sigma,
 \end{align*}
where $\gamma^i_*, \gamma^i$ are the restrictions of $\Gamma_*, \Gamma$ to the $i$-th edge of $G$, respectively, then
 \begin{equation}\label{eq:LojaIntro}
     \left|{\rm L}(\Gamma)-{\rm L}(\Gamma_*) \right|^{1-\theta}
     \le C_{\rm LS} \left( \sum_i \int_0^1 |\boldsymbol{k}^i|^2 \de s \right)^{\frac12},
\end{equation}
where $\boldsymbol{k}^i$ is the curvature of $\gamma^i$.
\end{teo}

Estimates like \eqref{eq:LojaIntro} are named after \L ojasiewicz and Simon due to their seminal works \cite{Loja63, Loja84, Si83}, where they firstly proved and employed analogous inequalities for analytic functionals over finite or infinite dimensional linear spaces. As we shall see, the validity of a \L ojasiewicz--Simon inequality is sufficient to imply strong stability properties of critical points of the energy under consideration.

\cref{thm:LojaIntro} is actually a particular case of a more general result yielding a \L ojasiewicz--Simon inequality for the length functional among triple junctions networks, that is, one does not need to ask that edges form angles equal to $\tfrac23\pi$ at junctions to get a version of \eqref{eq:LojaIntro}, see \cref{thm:LojaGenerale}. As discussed in \cref{rem:Loja2}, it is even possible to generalize the inequality to triple junctions networks letting endpoints free to vary.

Proving \cref{thm:LojaIntro} as a consequence of a \L ojasiewicz--Simon inequality holding among triple junctions networks not only gives an inequality for a much larger class of networks, but also simplifies its proof, since triple junctions networks do not need to satisfy an additional \emph{nonlinear} requirement on the angles at junctions.
Indeed the proof of \cref{thm:LojaGenerale}, which implies \cref{thm:LojaIntro}, eventually follows employing a by now established method for proving these kinds of inequalities for extrinsic geometric functionals \cite{ChFaSc09, DaPoSp16, PozzettaLoja, MantegazzaPozzettaLoja, MantegazzaPozzettaHypersurfaces}, the method relying on \emph{linear} Functional Analysis.\\
Once one is able to parametrize the considered competitors as normal graphs over the critical point, the inequality eventually follows from a general functional analytic result, see \cref{prop:LojaAbstract}, based on \cite{Ch03}. However, differently from the previous cases, the nonsmooth structure of networks necessarily implies technical complications, as networks close to a fixed minimal one $\Gamma_*$ cannot be written as normal graphs over the critical point. Hence, in order to perform a \emph{graph parametrization} of networks close to $\Gamma_*$ we need to allow for graphs having both a normal and a tangential component with respect to $\Gamma_*$. This would generally violate the assumptions needed to produce a \L ojasiewicz--Simon inequality, cf. \cref{prop:LojaAbstract}, since variations of $\Gamma_*$ via tangential directions are equivalent to reparametrizations of the curves of a network and thus generate an infinite dimensional kernel for a geometric functional like the length. We fix this issue by prescribing that tangential components of these graph parametrizations linearly depend on normal ones. Since a relation between such normal and tangential components is naturally satisfied at the junctions, the chosen dependence of tangential components on normal ones is given by a suitable prolongation of the relations at junctions on the interior of the edges, see \cref{prop:ParametrizzazioneTN}. We mention that an analogous construction has been also employed in \cite{garckegoesswein}.

\medskip

It is possible to exploit the \L ojasiewicz--Simon inequality in \cref{thm:LojaIntro} for proving the \emph{stability} of minimal networks with respect to the $L^2$-gradient flow of ${\rm L}$, the so-called \emph{motion by curvature of networks}. Along such flow, a regular network evolves keeping its endpoints fixed and moving with normal velocity equal to the curvature vector along each edge, see \cref{sec:Motion}. The motion by curvature generalizes the one-dimensional mean curvature flow, called curve shortening flow, to the realm of singular one-dimensional objects given by planar networks.

Bronsard and Reitich
\cite{Bronsardreitich} firstly attempted to find strong solutions to the motion by curvature, providing local existence and uniqueness of solutions for admissible initial regular networks of class $C^{2+\alpha}$ with the sum of the curvature at the junctions equal to zero. 
Then the basic theory concerning short time existence and uniqueness of the motion by curvature was carried out in~\cite{mannovplusch}, and further improved in~\cite{GMP} in order prove existence of the flow starting from any regular network without extra assumptions on the initial datum. In~\cite{GMP} also the parabolic regularization of the flow has been addressed.  
It is known that the flow develops singularities, see \cref{thm:LongTimeGeneral},
and a great deal of work has been done to understand the nature of these singularities
and to define the flow past singularities~\cite{MantegazzaNovagaTortorelli,MaMaNo16,mannovplusch,BaldiHausMan,IlmNevSch,liramazplusae,MaNoPl21,Chang,Chang2}.

However, exploiting the \L ojasiewicz--Simon inequality we can prove that a flow starting sufficiently close to a minimal network in $H^2$ exists for every time and smoothly converges to a (possibly different) minimal network. We mention that global existence of the flow starting close to critical points and convergence along diverging sequence of times has been firstly studied in \cite{KinLiu}. Hence the next theorem recovers and improves the main results of \cite{KinLiu}, see also \cref{thm:Convergence} below.

\begin{teo}[cf. \cref{thm:Stability}]\label{thm:StabilityIntro}
Let $\Gamma_*:G\to \R^2$ be a minimal network.
Then there exists $\delta_{\Gamma_*}>0$ such that the following holds.

Let $\Gamma_0:G\to \R^2$ be a smooth regular network having the same endpoints of   $\Gamma_*$ and such that $\|\gamma^i_0-\gamma^i_*\|_{H^2(\d x)}\le \delta_{\Gamma_*}$. Then the motion by curvature $\Gamma_t:G\to \R^2$ starting from $\Gamma_0$ exists for all times and it smoothly converges to a minimal network $\Gamma_\infty$ such that ${\rm L}(\Gamma_\infty)={\rm L}(\Gamma_*)$, up to reparametrization.
\end{teo}

The fact that $\Gamma_\infty$ may be different from $\Gamma_*$ in the above theorem cannot be avoided, as there exist examples of one-parameter families of minimal networks having fixed endpoints. Consider for instance networks given by concentric regular hexagons with segments connecting the six vertices
to fixed vertices of a bigger regular hexagon, as
depicted in~\cref{Fig:Ragnatela} below.

% (see also the discussion in~\cite[Fig. 1]{KinLiu})

%and even worse  is the situation of the flat torus $T^2$ where any translation of a minimal network is still a minimal network.

Observe that in \cref{thm:StabilityIntro} the initial datum $\Gamma_0$ does not need to satisfy further geometric properties at junctions or endpoints except being regular; this is possible as we shall employ the short time existence theory recently developed in \cite{GMP}, which removes the additional geometric assumptions required by previous existence theorems in \cite{MantegazzaNovagaTortorelli, mannovplusch}.  On the other hand, it is an open problem to understand whether a stability result like \cref{thm:StabilityIntro} holds for an initial datum $\Gamma_0$ that is possibly non-regular, such as a $\Gamma_0$ with only triple junctions but with angles possibly different from $\tfrac23\pi$ at junctions, and such that it is sufficiently close in $H^2$ to a reference minimal network. For such a network $\Gamma_0$, it is possible to define a motion by curvature $\Gamma_t$ starting from $\Gamma_0$ \cite{IlmNevSch, liramazplusae} and $\Gamma_t$ is instantaneously regular for $t>0$, however the crucial short time properties we need (see \cref{wellposedness} and \cref{lem:Regularization2}) are delicate in this setting.

Finally, we observe that no minimizing properties on $\Gamma_*$ in \cref{thm:StabilityIntro} are required. Instead, by means of a simple comparison argument it is possible to show that minimal networks automatically minimize the length among suitably small $C^0$ perturbations, see \cref{lemma:MinimalMinimizing}. Once such minimality property is established, the proof of \cref{thm:StabilityIntro} follows adapting a general argument outlined in \cite{MantegazzaPozzettaLoja, MantegazzaPozzettaHypersurfaces}.

As a consequence of \cref{thm:StabilityIntro}, it immediately follows that if a motion by curvature smoothly converges to a minimal network along a sequence of times, then it smoothly converges as time increases, see \cref{thm:Convergence}.

\medskip

The final main contribution of this work is given by the rigorous construction of an example of motion by curvature presenting a topological singularity in infinite time.

It is known, see \cref{thm:LongTimeGeneral}, that the motion by curvature may develop singularities consisting in the blow up of the $L^2$-norm of the curvature or in the disappearance of a curve whose length tends to zero. There exist well known examples where both or one of the previous alternatives occur in finite time, see for instance \cite{mannovplusch} and \cite[Section 6]{MaNoPl21}. Here we firstly construct an example of a motion by curvature existing for every time and such that a curve of the evolving network vanishes in infinite time while the curvature remains uniformly bounded.

\begin{teo}[cf. \cref{thm:EsempioConvergenceInfiniteTime}]\label{thm:ExampleIntro}
There exists a smooth regular network $\Gamma_0:G\to \R^2$ such that the motion by curvature $\Gamma_t$ starting from $\Gamma_0$ exists for every time, the length of each curve $\gamma^i_t$ is strictly positive for any time, the curvature of each curve $\gamma^i_t$ is uniformly bounded from above, and $\Gamma_t$ smoothly converges to a degenerate network $\Gamma_\infty$ as $t\to+\infty$, up to reparametrization. Specifically, the length of a distinguished curve $\gamma^0_t$ tends to zero as $t\to+\infty$.
\end{teo}

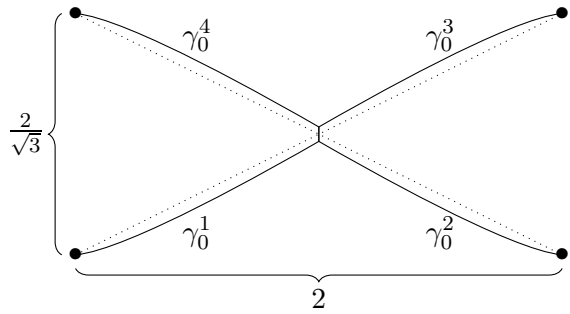
\begin{figure}[H]
	\begin{center}
		\begin{tikzpicture}[scale=3.2]
		\draw
		(-1,-0.5)to[out=5, in=210, looseness=0.5]
		(0,-0.03)to[out=90, in=270, looseness=0]
		(0,0);
		\draw[cm={cos(180) ,-sin(180) ,sin(180) ,cos(180) ,(0 cm,0 cm)}]
		(-1,-0.5)to[out=5, in=210, looseness=0.5]
		(0,-0.03)to[out=90, in=270, looseness=0]
		(0,0);
		\draw[xscale=-1]
		(-1,-0.5)to[out=5, in=210, looseness=0.5]
		(0,-0.03)to[out=90, in=270, looseness=0]
		(0,0);
		\draw[yscale=-1]
		(-1,-0.5)to[out=5, in=210, looseness=0.5]
		(0,-0.03)to[out=90, in=270, looseness=0]
		(0,0);
		\path[font=\normalsize]
		(-0.5,-0.3)node[below]{$\gamma^1_0$};
		\path[font=\normalsize]
		(0.5,-0.3)node[below]{$\gamma^2_0$};
		\path[font=\normalsize]
		(-0.5,0.3)node[above]{$\gamma^4_0$};
		\path[font=\normalsize]
		(0.5,0.3)node[above]{$\gamma^3_0$};
		\path[font=\normalsize]
		(-1,-0.5)node[]{$\bullet$};
		\path[font=\normalsize]
		(-1,0.5)node[]{$\bullet$};
		\path[font=\normalsize]
		(1,-0.5)node[]{$\bullet$};
		\path[font=\normalsize]
		(1,0.5)node[]{$\bullet$};
		\draw[dotted]
		(-1,-0.5)--(1,0.5);
		\draw[dotted]
		(-1,0.5)--(1,-0.5);
		\draw [decorate,decoration={brace,amplitude=5pt,mirror,raise=1ex}]
        (-1,-0.5) -- (1,-0.5)
        node[midway,yshift=-1.5em]{$2$};
        \draw [decorate,decoration={brace,amplitude=5pt,raise=1ex}]
        (-1,-0.5) -- (-1,0.5)
        node[midway,xshift=-1.7em]{$\frac{2}{\sqrt{3}}$};
		\end{tikzpicture}
	\end{center}
	\caption{Continuous lines: initial datum $\Gamma_0$. Dotted lines: limit degenerate network $\Gamma_\infty$. By the choices of the mutual distances of endpoints, the dotted lines intersect forming angles equal to $\tfrac{\pi}{3}$ and $\tfrac{2}{3}\pi$.}\label{Fig:Network}
\end{figure}

For proving the above theorem we will specifically consider the evolution of the network sketched in \cref{Fig:Network}.

For $\Gamma_0$ as in \cref{Fig:Network}, the endpoints determine a rectangle whose diagonals intersect forming angles exactly equal to $\tfrac{\pi}{3}$ and $\tfrac{2}{3}\pi$. We will prove that the central curve in \cref{Fig:Network} shrinks along the flow in infinite time and the four remaining curves converge to the diagonals connecting the four endpoints.

We will study the motion by curvature starting from such $\Gamma_0$ explicitly. The symmetry chosen for producing this example allows to generalize some ideas from \cite{EckerHuisken89} to the context of networks in order to uniformly bound the curvature along the motion exploiting a monotonicity-type formula, see \cref{thm:Monotonicity}. Eventually a comparison with solutions of heat-type equations show that convergence occurs in infinite time.

We stress that the example of \cref{thm:ExampleIntro} yields a simple and explicit flow converging to a degenerate critical point of the length, implying that the topology of the evolving networks changes in the limit. It is the change of topology the fundamental reason why such an example has to be studied individually and whose convergence cannot follow from the \L ojasiewicz--Simon inequalities proved in \cref{thm:LojaIntro} or \cref{thm:LojaGenerale}. This simple example motivates to look for improvements on the general method for proving convergence, as well as the search of possibly weaker variants of \L ojasiewicz--Simon inequalities, able to take into account these changes of topology in the limit. This project goes beyond the scope of the present paper and it is left for future investigation.

\medskip

We finally observe that the analysis carried out for the study of the example proving \cref{thm:ExampleIntro} actually provides a \emph{family} of examples of evolving networks whose curvature is uniformly bounded presenting \emph{every possible} long time behavior:  collapsing in infinite time of the length of a curve, collapsing in finite time of the length of a curve, convergence in infinite time to a regular network. The first case is the one claimed in \cref{thm:ExampleIntro}, the other cases are simple modifications of it and they are described more precisely in the next remark.

\begin{rem}\label{rem:FamigliaEsempi}
Let $L>0$. Consider a smooth regular initial datum $\Gamma_{0,L}$ completely analogous to the one in \cref{Fig:Network}, whose four curves connected to the endpoints are always the same, but the initial length of the central vertical edge equals $L$ at time zero. In particular there is $\overline{L}$ such that $\Gamma_{0,\overline{L}}$ is exactly the one in \cref{Fig:Network}, that is, the resulting rectangle determined by the endpoints has sides of length $2/\sqrt{3}$ and $2$.

It can be proved (see \ref{Step1} below) that the convexity of the four external curves is preserved along the flow, and then the central curve remains vertical and its length is strictly decreasing. Moreover the curvature of each curve is uniformly bounded (see first part of \ref{Step2} below) and, by comparison, each of the four external curves always remains on the same side of the suitable line; for instance, the evolved curve $\gamma^1_t$ stays below the line passing through the same endpoint of $\gamma^1_t$ forming an angle equal to $\tfrac\pi6$ with the horizontal axis (see the argument in \ref{Step3} below).

Let $T_L>0$ be the maximal time of existence of the motion by curvature starting from $\Gamma_{0,L}$. By uniqueness and locality of the flow, it is clear that the evolution of the four external curves of $\Gamma_{t,L}$ coincide with the one of the four external curves of $\Gamma_{t,L'}$ for any $t \in [0,\min\{T_L,T_{L'}\})$ (in particular the evolution of such curves is independent of the length of the central curve).

It follows that for $L>\overline{L}$ the length of the central curve is always bounded away from zero and then the above observations imply that $T_L=+\infty$ and the flow smoothly converges in infinite time to a minimal network.

On the contrary, if $L\in(0,\overline{L})$, the length of the central curve vanishes in finite time, leading to a topological singularity of the flow in finite time.
\end{rem}

We conclude this introduction by mentioning some further contributions related to the topic.

The use of \L ojasiewicz--Simon inequalities has become a prominent tool for understanding stability properties and convergence of geometric flows.
Apart from the above mentioned references, let us also recall the recent results on the uniqueness of blow-ups for the mean curvature flow \cite{Schulze, ColdingMinicozziUniqueness, ChodoshSchulze}, and the application of the method to constrained high order extrinsic flows \cite{Rupp1,Rupp2,garckegoesswein}. The use of these inequalities has seen successful applications also in the context of intrinsic geometric flows, namely in \cite{CarlottoChodoshRub, FeehanGlobalYM}.

Apart from stability results, which is the main focus of the current paper,
there are several questions
concerning the 
motion by curvature of networks:
study of singularities, 
global existence, extension to class of weaker objects.
As we said above there is an
extensive amount of literature concerning the analysis of the flow
in the framework of classical PDEs (see~\cite{mannovplusch} and references
within).
There are also several generalized weak
notions of the flow, for instance~\cite{brakke,es,KiTo,LaOt,StTo}.
Recently interesting progresses has been
made both in the direction of proving regularity of weak solution~\cite{KiTo,KiTo2} and in establishing the so-called weak-strong uniqueness theorem~\cite{FiHeLaSi,HeLa2}.

Worth to mention the fact that the motion by curvature of network was first proposed for modelling reasons~\cite{mullins} 
and recently has again 
attracted the applied mathematical community 
~\cite{applicato1,applicato2,applicato3,applicato4}.

\bigskip
\textbf{Organization.} In \cref{sec:Preliminaries} we collect basic definitions and results on networks and on the motion by curvature. In \cref{sec:Loja} we establish the graph parametrization of networks over minimal ones and we prove the \L ojasiewicz--Simon inequality implying \cref{thm:LojaIntro}. In \cref{sec:MinimalMinimizing} we prove that minimal networks locally minimize the length in $C^0$. \cref{sec:StabilityConvergence} is devoted to the proof of the stability of minimal networks, implying \cref{thm:StabilityIntro}. In \cref{sec:Example} we prove \cref{thm:ExampleIntro} by analyzing the motion of networks like the one in \cref{Fig:Network}. In \cref{sec:Appendix} we collect some tools needed in the proofs, namely a well known quantitative implicit function theorem and a monotonicity-type formula. In \cref{sec:AppendixSurfaces} we discuss extensions of our results to the case of networks on Riemannian surfaces.

\bigskip
\textbf{Acknowledgements.} The authors are partially supported by the INdAM - GNAMPA Project 2022 CUP \_ E55F22000270001 ``Isoperimetric problems: variational and geometric aspects''. The authors are grateful to Matteo Novaga for fruitful discussions on the topic of this work.

\bigskip

\section{Preliminaries}\label{sec:Preliminaries}

\subsection{Networks}\label{sec:Networks}

For a regular curve $\gamma:[0,1]\to \R^2$ of class $H^2$, define
\[
\tau_\gamma\eqdef\frac{\gamma'}{|\gamma'|},
\quad
\nu_\gamma \eqdef {\rm R}(\tau_\gamma),
\]
the tangent and the normal vector, respectively, where ${\rm R}$ denotes counterclockwise rotation of $\tfrac\pi2$. We define $\de s_\gamma\eqdef |\gamma'| \de x$ the arclength element and $\partial_s\eqdef |\gamma'|^{-1}\partial_x$ the arclength derivative. The curvature of $\gamma$ is the vector
\[
\boldsymbol{k}_\gamma\eqdef \partial_s^2 \gamma.
\]
We shall usually drop the subscript $\gamma$ when there is no risk of confusion.

\medskip

Fix $N\in\mathbb{N}$ and let $i\in\{1,\ldots, N\}$, $E^i:=[0,1]\times\{i\}$,
$E:=\bigcup_{i=1}^N E^i$ and 
$V:=\bigcup_{i=1}^N \{0,1\}\times\{i\}$.

\begin{defn}\label{defgraph}
Let $\sim$ be an equivalence relation that identifies points of $V$.
A \emph{graph} $G$ is the topological quotient space of $E$ induced by $\sim$, that is
\begin{equation*}
G:=E\Big{/}\sim\,,
\end{equation*}
and we assume that $G$ is connected.
\end{defn}

\begin{defn}\label{network}
A \emph{(planar) network} is a pair $\mathcal{N}=(G,\Gamma)$
where
\begin{equation*}
\Gamma: G\to \mathbb{R}^2
\end{equation*}
is a continuous map and $G$ is a graph.
We say that $\mathcal{N}$ is of class $W^{k,p}$ (resp. $C^{k,\alpha}$) if each map $\gamma^i:=\Gamma_{\vert E^i}$ is
either a constant map (\emph{singular} curve) or a regular curve of class $W^{k,p}$ (resp. $C^{k,\alpha}$ up to the boundary). A network is \emph{smooth} if it is of class $C^\infty$.  A network is \emph{degenerate} if there is at least one singular curve.

Denoting by $\pi:E\to G$ the projection onto the quotient, an \emph{endpoint} is a point $p \in G$ such that $\pi^{-1}(p) \subset V$ and it is a singleton, a \emph{junction} is a point $m \in G$ such that $\pi^{-1}(m) \subset V$ and it is not a singleton. The \emph{order} of a junction if the cardinality $\sharp \pi^{-1}(m)$.

We denote by $J_G$, and by $P_G$ respectively, the set of junctions, and endpoints respectively, of a graph $G$.
A graph $G$ is said to be \emph{regular} if each junction has order $3$.
\end{defn}

Without loss of generality, if $\mathcal{N}=(G,\Gamma)$ is a network and $p \in G$ is an endpoint with $\pi^{-1}(p)=(e,i)$,
we will implicitly assume that $e=1$.

\begin{defn}
Let $\mathcal{N}=(G,\Gamma)$ be a network of class $C^1$ and let $e \in \{0,1\}$. The \emph{inner tangent vector} of a regular curve $\gamma^i$ of $\mathcal{N}$ at $e$ is the vector
\[
(-1)^{e} \frac{(\gamma^i)'(e)}{|(\gamma^i)'(e)|}.
\]
\end{defn}

In this paper we will be interested in the classes of \emph{triple junctions networks} and \textit{regular networks}. 

\begin{defn}\label{def:TripleJunctionsRegularMinimal}
A network $\mathcal{N}=(G,\Gamma)$ is a \emph{triple junctions network} if $G$ is regular and each map $\gamma^i:=\Gamma_{\vert E^i}$ is a regular embedding of class $C^1$,
for every $i\ne j$ the curves $\gamma^i$
and $\gamma^j$ do not intersect in their interior, and $\pi(0,i)\neq \pi(1,i)$ for any $i$.

A network $\mathcal{N}=(G,\Gamma)$ is said to be
\emph{regular} if it is 
a triple junctions network such that whenever $\pi(e^i,i)=\pi(e^j,j)=\pi(e^k,k)$ is a junction then any two inner tangent vectors of $\gamma^i,\gamma^j,\gamma^k$ at $e^i,e^j,e^k$, respectively, form an angle equal to $\tfrac23 \pi$.

A network $\mathcal{N}=(G,\Gamma)$ is said to be \emph{minimal} if it is regular and the curvature of the parametrization of each edge is identically zero. Moreover, we assume that the parametrizations of minimal networks have constant speed.
\end{defn}

We shall usually denote a network by directly writing the map $\Gamma:G\to \R^2$.

\subsection{Function spaces}

We introduce the space $W^{1,2}_p$,
which is a natural choice to 
define the motion by curvature of networks.

For $T>0$, $N \in \N$ with $N\ge1$, $p \in (3,\infty)$, we define
\begin{equation*}
    W^{1,2}_p\left((0,T)\times(0,1); \R^{2N} \right) \eqdef W^{1,p}\left((0,T); L^p\left((0,1);\R^{2N}\right) \right) \cap L^p\left((0,T); W^{2,p}\left((0,1);\R^{2N}\right) \right).
\end{equation*}

\begin{rem}
Elements in the space $W^{1,2}_p$
are functions $f\in L^p\left((0,T); L^p(0,1)\right)$ possessing one distributional derivative with respect to time $\partial_t f\in L^p\left((0,T); L^p(0,1)\right)$. Furthermore, for almost every $t\in(0,T)$, the function $f(t)$ lies in $W^{2,p}(0,1)$ and thus has two spacial derivatives $\partial_x (f(t))$, $\partial_x ^2\left(f(t)\right)\in L_p(0,1)$. One easily sees that the functions $t\mapsto \partial_x^k(f(t))$ for $k\in\{1,2\}$ lie in $L_p\left((0,T);L_p(0,1)\right)$. 

The space $W^{1,2}_p$ is defined as the intersection of two Bochner spaces, which are, in this case, Sobolev spaces of functions defined on a measure space with values in a Banach space.
We remind that 
$$
\left\lVert f\right\rVert_{ L^p\left((0,T);W^{2,p}((0,1);\mathbb{R}^{2N})\right)}:=\left\lVert \left\lVert f(\cdot)\right\rVert_{W^{2,p}((0,1);\mathbb{R}^{2N})}\right\rVert_{L_p\left((0,T);\mathbb{R}\right)}
$$
and
$$
\lVert f\rVert_{W^{1,p}\left((0,T); L^p\left((0,1);\R^{2N}\right) \right) }:=
\left(
\lVert f\rVert_{L^p((0,T);L^p((0,1);\mathbb{R}^{2N}))}^p+\lVert\partial_t f\rVert_{ L^p((0,T);L^p((0,1);\mathbb{R}^{2N}))}^p\right)^{\nicefrac{1}{p}}
$$
\end{rem}

\medskip

Let $\Gamma_t:[0,T)\times G\to \R^2$ be a  time-dependent network
parametrized by $(\gamma^1_t,\ldots,\gamma^N_t)$
with $\gamma^i_t\in W^{1,2}_p$.
We shall denote  $(G,\Gamma_t)$
by the symbol $(\mathcal N_t)_t$
and
\begin{equation*}
    \|(\mathcal N_t)_t\|_{W^{1,2}_p} \eqdef
    \bigg( \int_0^T  \|\partial_t \Gamma_t\|_{L^p}^p \de t
    \bigg)^{\frac1p}
    + \bigg(\int_0^T \|\Gamma_t\|_{W^{2,p}}^p \de t
    \bigg)^{\frac1p},
\end{equation*}
where
\begin{equation*}
    \|\partial_t \Gamma_t\|_{L^p}^p = \sum_i  \|\partial_t \gamma^i_t\|_{L^p(\d x)}^p,
    \qquad
     \|\Gamma_t\|_{W^{2,p}}^p =  \sum_i  \|\gamma^i_t\|_{W^{2,p}(\d x)}^p\,.
\end{equation*}

\medskip

We also need to introduce a suitable space for initial data.
Fixed $p\in (3,\infty)$ the \textit{Sobolev--Slobodeckij space} $W^{2-\nicefrac{2}{p},p}\left((0,1);\mathbb{R}^{2N}\right)$ is defined by
\begin{equation*}
W^{2-\nicefrac{2}{p},p}\left((0,1);\mathbb{R}^{2N}\right):=\left\{f\in W^{1,p}\left((0,1);\mathbb{R}^{2N}\right):\left[\partial_x f\right]_{1-\nicefrac{2}{p},p}<\infty\right\}\,
\end{equation*}
with 
\begin{equation*}
\left[\partial_x f\right]_{1-\nicefrac{2}{p},p}
:=\left(\int_{0}^{1}\int_{0}^{1}\frac{\left\lvert \partial_x f(x)-\partial_x f(y)\right\rvert^p}{\lvert x-y\rvert^{p-1}}\,\mathrm{d}x\,\mathrm{d}y\right)^{\nicefrac{1}{p}}\,.
\end{equation*}
We define the $W^{2-\nicefrac{2}{p}, p}$-norm of a network $\Gamma_0$ by
\begin{equation*}
    \|\Gamma_0\|_{W^{2-\nicefrac{2}{p}, p}} \eqdef
     \| \Gamma_0\|_{W^{1,p}}
    +\left[\partial_x\Gamma_0\right]_{1-\nicefrac{2}{p},p},
\end{equation*}
where
\begin{equation*}
    \|\Gamma_0\|_{W^{1,p}} = \sum_i  \| \gamma^i_0\|_{W^{1,p}},
    \qquad
     \left[\partial_x\Gamma_0\right]_{1-\nicefrac{2}{p},p} =  \sum_i  \left[\partial_x\gamma^i_0\right]_{1-\nicefrac{2}{p},p}\,.
\end{equation*}	
	
\begin{rem}
The temporal trace of the space
$W^{1,2}_p$ is the space $W^{2-\nicefrac{2}{p},p}$.
Since we would like to set the problem in 
$W^{1,2}_p$,  
it is then natural to chose $W^{2-\nicefrac{2}{p},p}$
as the space for the initial data.
Moreover we ask $p\in(3,\infty)$ to write the Herring condition at the junction.
Indeed, for any
$T>0$, $p\in(3,\infty)$ and $\alpha\in\left(0,1-\nicefrac{3}{p}\right]$ we have the continuous embeddings
\begin{equation*}
W_p^{1,2}\left((0,T)\times(0,1);\mathbb{R}^{2N}\right)\hookrightarrow C\left([0,T];W^{2-\nicefrac{2}{p},p}\left((0,1);\mathbb{R}^{2N}\right)\right)\hookrightarrow C\left([0,T];C^{1+\alpha}\left([0,1];\mathbb{R}^{2N}\right)\right)\,.
	\end{equation*}
The first embedding follows from~\cite[Lemma 4.4]{DenkSaalSeiler}, the second is an immediate consequence of the Sobolev Embedding Theorem.
Again by the Sobolev Embedding Theorem for  $p\in(3,6)$ we have the compact embedding
\begin{equation*}
H^2\left((0,1);\mathbb{R}^{2N}\right)\hookrightarrow 
W^{2-\nicefrac{2}{p},p}\left((0,1);\mathbb{R}^{2N}\right)\,.
	\end{equation*}
\end{rem}

\subsection{Motion by curvature of planar networks}\label{sec:Motion}

In this section we introduce the basic definitions and known results on the motion by curvature. 
Let $p\in(3,\infty)$ be fixed.

\begin{defn}[Admissible initial datum]\label{geomadm}
A network 
$\mathcal{N}_0=(G,\Gamma_0)$ is an 
\textit{admissible initial datum}
for the motion  by curvature
if it is regular network of class $W^{2-\nicefrac{2}{p},p}$. 
\end{defn}

\begin{defn}[Solutions to the motion by curvature]\label{def:SolutionMotion}
Let $\Gamma_0:G\to \R^2$ be an admissible initial datum with $P^i=\Gamma_0(p^i)\in\mathbb{R}^2$. A one-parameter family of regular networks $\Gamma_t:G\to \R^2$, for $t \in [0,T)$, is a \textit{solution to the motion by curvature
with initial datum $\Gamma_0$} if the parametrizations $\gamma^i_t$ of $\Gamma_t$ satisfy
\begin{equation}\label{eq:motion}
    \begin{split}
        \gamma^i_t(p)=P^i & \qquad \forall\, t \in [0,T), \,\,\forall\,p \text{ endpoint of $G$},\\
        \left\langle\partial_t\gamma^i(t,x),\nu^i(t,x)\right\rangle\nu^i(t,x)
        =\boldsymbol{k}^i(t,x)
        & \qquad\text{for a.e. $t\in(0,T)$, $x \in (0,1)$},
    \end{split}
\end{equation}
and the collection of parametrizations $(\gamma^1_t,\ldots,\gamma^N_t)$ belongs to $W^{1,2}_p\left((0,T)\times(0,1); \R^{2N} \right)$, with $\gamma^i_t|_{t=0}=\gamma^i_0$ for any $i$.

The solution is assumed to be maximal, i.e., we ask that it does not exist
another solution defined on $[0,\widetilde{T})$ with $\widetilde{T}>T$.
\end{defn}

\begin{rem}
We stress the fact that the evolving network must be regular for every time $t\in [0,T)$. From the PDE point of view this means that the system~\eqref{eq:motion} is a boundary value problem
with coupled boundary conditions:
whenever $\pi(e^i,i)=\pi(e^j,j)=\pi(e^\ell,\ell)$ is a junction then for all $t\in [0,T)$
\begin{align*}
&\gamma^i(t,e^i)=\gamma^j(t,e^j)= \gamma^\ell(t,e^\ell)\,,\\
&(-1)^{e^i}\tau^i(t,e^i)+(-1)^{e^j}\tau^j(t,e^j)+(-1)^{e^k}\tau^\ell(t,e^\ell)=0\,.
\end{align*}
\end{rem}

We quote here a parade of results on the motion by curvature of network that are relevant for the sequel of this paper.

\begin{teo}[Short time existence and parabolic smoothing~\cite{GMP}]\label{wellposedness}
If $\Gamma_0:G\to \R^2$ is an 
admissible initial datum, then there exists a solution $\Gamma_t:G\to\R^2$ to the motion by curvature starting from $\Gamma_0$, for $t \in [0,T)$. The solution is unique up to reparametrizations, it is smooth on $[\eps,T-\eps]\times G$ for any $\eps>0$, and $\gamma^i_t\to\gamma^i_0$ in $C^1([0,1])$ for $t\to0^+$ for any $i$.

Moreover, for any $c_1,c_2>0$, there are $\tau=\tau(c_1,c_2),M=M(p,c_1,c_2)>0$ such that if
\begin{equation*}
    \min_{i\in\{1,\ldots,N\},x\in[0,1]} |\partial_x \gamma^i_0| \ge c_1\quad\text{and}\quad\|\Gamma_0\|_{W^{2-\nicefrac{2}{p}, p}} \leq c_2\,,
\end{equation*} 
then $T\ge \tau$ and the solution $\mathcal N$ satisfies $\|(\mathcal N_t)_t\|_{W^{1,2}_p} \le M$ for any $t \in [0,\tau]$.
\end{teo}

To show existence of solutions one finds a unique solution to the \emph{special flow}, i.e., the evolution determined by the non-degenerate parabolic second order equation
\begin{equation}\label{eq:SpecialFlow}
\partial_t \gamma^i_t = 
\frac{\partial^2_x\gamma^i_t}{|\partial_x\gamma^i_t|^2}\,.
\end{equation}
Clearly $\left\langle
\partial_t \gamma^i_t,\nu^i \right\rangle \nu^i= \boldsymbol{k}^i$,
so a solution to the special flow is in particular a solution to
the motion by curvature.
Uniqueness up to reparametrizations 
is obtained by showing that any solution of the network
flow can be obtained as a reparametrization of the solution of the special flow.

\begin{rem}
\cref{wellposedness} yields existence
and uniqueness of a solution in the sense of \cref{def:SolutionMotion} starting from \emph{any} regular network of class $W^{2-\nicefrac{2}{p},p}$. This is a great advantage in comparison with the theory firstly developed in~\cite{Bronsardreitich,mannovplusch}, where an initial datum $\Gamma_0$ parametrized by $\gamma^i_0$ is required not only to be regular and $C^2$, but also suitably \emph{geometrically compatible}, that is such that
\begin{equation}\label{eq:GeomCompatibilityConditions}
    \boldsymbol{k}^i(p)=0,
    \qquad
(-1)^{e_i}k^i(e^i)+(-1)^{e_j}k^j(e^j)+(-1)^{e_\ell}k^\ell(e^\ell)=0,
\end{equation}
at any endpoint $p \in G$ and at any junction $m = \pi(e^i,i)=\pi(e^j,j)=\pi(e^\ell,\ell)$, where $k^n$ is the oriented curvature of $\gamma^n$ for any $n$.
\end{rem}

\begin{rem}\label{rem:SoluzioniConfrontoCompatibilita}

We stress that for any admissible initial datum in Theorem~\ref{wellposedness}, the results in \cite{GMP} imply that we can take as solution to the motion by curvature exactly the solution to the special flow. More precisely, for any positive time the parametrizations $\gamma^i_t$ of the solution verify the equation~\eqref{eq:SpecialFlow}
with no need of reparametrize the solution or the initial datum. We will always assume that the solution to the motion by curvature satisfies \eqref{eq:SpecialFlow} whenever nothing different is specified.

As a consequence, the solution satisfies the analytic compatibility conditions of
every order (see \cite[Definition 4.7, Definition 4.16]{mannovplusch}). In particular, the following compatibility conditions of order two hold:
\begin{equation}\label{eq:CompatibilityConditions}
    \partial^2_x\gamma^k(p)=0\,,
    \qquad
    \frac{\partial^2_x\gamma^i(e^i)}{|\partial_x\gamma^j(e^i)|^2} = \frac{\partial^2_x\gamma^j(e^j)}{|\partial_x\gamma^j(e^j)|^2}\,,
\end{equation}
at any endpoint $p \in G$ and at any junction $m = \pi(e^i,i)=\pi(e^j,j)=\pi(e^\ell,\ell)$.

%Once we get the smooth solution of the special flow, we can re-parametrized it to a smooth solution of the motion by curvature which fulfil the \emph{geometric} compatibility conditions of every order.

%Also the re-parametrization of the curves of the networks with constant speed equal to their length is smooth for positive times. 

The fact that for any positive times we can take a smooth solution to the special flow is a key point in our analysis, as it allows to apply the classical results and to use all the estimates derived in~\cite{MantegazzaNovagaTortorelli,MaMaNo16,mannovplusch}.
\end{rem}

In the next statement we recall the possible singularities happening at a singular time.

\begin{teo}[Long time behavior~\cite{mannovplusch}]\label{thm:LongTimeGeneral}
Let $(\mathcal{N}_t)_t$ be a solution to the motion by curvature with initial datum $\Gamma_0$ in the time interval $[0,T)$.
Then either 
$$
T=+\infty,
$$
or as $t\to T$ at least one of the following happens:
\begin{itemize}
    \item[i)]  the limit inferior of the length of at least one curve of the network is zero;
    \item[ii)] the limit superior of the $L^2$-norm of the curvature is $+\infty$.
\end{itemize}
\end{teo}

As mentioned in the introduction, the possibilities listed in the above theorem are not mutually exclusive.

\section{\L ojasiewicz--Simon inequalities}\label{sec:Loja}

\subsection{Graph parametrization of regular networks}\label{sec:GraphParametrization}

We will employ the following notation. Let $\Gamma_*:G\to \R^2$ be a regular network. Denote by $\gamma^i_*$ the parametrization of the $i$-th edge of $\Gamma_*$, and by $\tau^i_*,\nu^i_*$ the relative unit tangent and normal vectors. Whenever $m\eqdef\pi(e^i,i)=\pi(e^j,j)$ is a junction, denote
\begin{equation}\label{eq:DefAlphaBeta}
    \alpha^{ij}_m \eqdef \scal{\tau_*^i(e^i),\tau^j_*(e^j)}
 %   = - \frac12 (-1)^{e^i+e^j},
    \qquad
    \beta^{ij}_m \eqdef \scal{\tau^i_*(e^i),\nu^j_*(e^j)} 
  %  = -\frac{\sqrt{3}}{2}  (-1)^{e^i+e^j}.
\end{equation}
Observe that $\alpha^{ij}_m=\scal{\nu^i_*(e^i),\nu^j_*(e^j)}$, $\alpha^{ij}_m=\alpha^{ji}_m$ and $\beta^{ij}_m=-\beta^{ji}_m$.

% {\color{blue}
% \begin{rem}
% Suppose $m\eqdef\pi(e^i,i)=\pi(e^j,j)=\pi(e^k,k)$ is a junction. Then
% \begin{align*}
%  & (-1)^{e^i+e^j}\alpha^{ij}_m=(-1)^{e^j+e^k}\alpha^{jk}_m
% =(-1)^{e^k+e^i}\alpha^{ki}_m \,, \\
%  & (-1)^{e^i+e^j}\beta^{ij}_m=(-1)^{e^j+e^k}\beta^{jk}_m
% =(-1)^{e^k+e^i}\beta^{ki}_m  \,.
% \end{align*}
% Observe that $\alpha^{ij}_m=\scal{\nu^i_*(e^i),\nu^j_*(e^j)}$, $\alpha^{ij}_m=\alpha^{ji}_m$ and $\beta^{ij}_m=-\beta^{ji}_m$. 
% Moreover, if $m$ is a junction of a regular network
% $\alpha^{ij}\in\{\pm \frac{1}{2}\}$ and 
% $\beta^{ij}\in\{\pm \frac{\sqrt{3}}{2}\}$.
% \end{rem}
% }
We derive the necessary conditions holding at junctions for the parametrizations of a triple junctions network written as a graphs over a regular one. 

\begin{lemma}\label{lem:BoundaryRelationsTN}
Let $\Gamma_*:G\to \R^2$ be a regular network
and $\Gamma:G\to \R^2$ be a network such that
at a junction $m=\pi(e^i,i)=\pi(e^j,j)=\pi(e^k,k)$
we have
\begin{equation*}
    \gamma^\ell(e^\ell) = \gamma^\ell_*(e^\ell) + \NN^\ell(e^\ell)\nu^\ell_*(e^\ell) + \TT^\ell(e^\ell)\tau^\ell_*(e^\ell),\quad\text{with}\;
    \ell\in\{i,j,k\}
\end{equation*}
for some constants $\NN^\ell(e^\ell),\TT^\ell(e^\ell) \in \R$.

Then the following relations hold
\begin{equation}\label{eq:BoundaryRelationsTN}
\begin{split}
     \TT^i(e^i) 
      &=
      \frac{1}{1-\alpha^{ij}_m\alpha^{jk}_m\alpha^{ki}_m}
      \left(\NN^j(e^j)\beta_m^{ij} +  \NN^k(e^k)\alpha^{ij}_m\beta_m^{jk} +  \NN^i(e^i)\alpha^{ij}_m\alpha^{jk}_m\beta_m^{ki}  \right)  ,\\
           \TT^j(e^j) 
      &=
      \frac{1}{1-\alpha^{ij}_m\alpha^{jk}_m\alpha^{ki}_m}
      \left(\NN^k(e^k)\beta_m^{jk} +  \NN^i(e^i)\alpha^{jk}_m\beta_m^{ki} +  \NN^j(e^j)\alpha^{jk}_m\alpha^{ki}_m\beta_m^{ij}  \right)  ,\\
                 \TT^k(e^k) 
      &=
      \frac{1}{1-\alpha^{ij}_m\alpha^{jk}_m\alpha^{ki}_m}
      \left(\NN^i(e^i)\beta_m^{ki} +  \NN^j(e^j)\alpha^{ki}_m\beta_m^{ij} +  \NN^k(e^k)\alpha^{ki}_m\alpha^{ij}_m\beta_m^{jk}  \right)  ,\\
\end{split}
\end{equation}
\begin{equation}\label{eq:SommaNN}
    (-1)^{e^i}\NN^i(e^i) + (-1)^{e^j}\NN^j(e^j) + (-1)^{e^k}\NN^k(e^k) = 0.
\end{equation}
In particular 
\begin{equation}\label{eq:LineariNiNj}
\begin{split}
\TT^i(e^i) 
     &=
     \frac{1}{1-\alpha^{ij}_m\alpha^{jk}_m\alpha^{ki}_m}
     \left(
     \NN^i(e^i)\left(\alpha^{ij}_m\alpha^{jk}_m\beta_m^{ki} +(-1)^{1+e_i+e_k}\alpha^{ij}_m\beta_m^{jk} \right)+
    \NN^j(e^j)\left(\beta_m^{ij} 
      +(-1)^{1+e_j+e_k}\alpha^{ij}_m\beta_m^{jk}\right)\right)  ,\\
           \TT^j(e^j) 
      &=
      \frac{1}{1-\alpha^{ij}_m\alpha^{jk}_m\alpha^{ki}_m}
      \left(
    \NN^i(e^i)\left(\alpha^{jk}_m\beta_m^{ki} 
       +(-1)^{1+e_i+e_k}\beta_m^{jk}
      \right)+  \NN^j(e^j)\left(\alpha^{jk}_m\alpha^{ki}_m\beta_m^{ij}  +(-1)^{1+e_j+e_k}\beta_m^{jk}\right)\right)  ,\\
                 \TT^k(e^k) 
      &=
      \frac{1}{1-\alpha^{ij}_m\alpha^{jk}_m\alpha^{ki}_m}
      \left(\NN^i(e^i)
      \left(\beta_m^{ki} +(-1)^{1+e_i+e_k} \alpha^{ki}_m\alpha^{ij}_m\beta_m^{jk} \right)+ \NN^j(e^j)\left(\alpha^{ki}_m\beta_m^{ij}+(-1)^{1+e_j+e_k} \alpha^{ki}_m\alpha^{ij}_m\beta_m^{jk}\right)    \right)  .
\end{split}
\end{equation}

% \begin{equation}\label{eq:LineariNiNj}
% \begin{split}
% &\TT^i(e^i) = 
% -\frac{1}{\sqrt{3}}\NN^i(e^i)
% -\frac{2}{\sqrt{3}}(-1)^{e_i+e_j}\NN^j(e^j) \,,\\
% &\TT^j(e^j) = 
% \frac{2}{\sqrt{3}}(-1)^{e_i+e_j}\NN^i(e^i)   +\frac{1}{\sqrt{3}}\NN^j(e^j)\,,\\
% &\TT^k(e^k) = 
% -\frac{1}{\sqrt{3}}(-1)^{e_i+e_k}\NN^i(e^i)
% +\frac{1}{\sqrt{3}}(-1)^{e_j+e_k}\NN^j(e^j)
% \,,\\
% \end{split}
% \end{equation}
\end{lemma}

\begin{proof}
Let $m=\pi(e^i,i)=\pi(e^j,j)=\pi(e^k,k)$ be a junction of $\Gamma$. Then $\gamma^i(e^i)=\gamma^j(e^j)=\gamma^k(e^k)$, that is
\begin{align}
\NN^i(e^i)\nu^i_*(e^i) + \TT^i(e^i)\tau^i_*(e^i) 
&= 
\NN^j(e^j)\nu^j_*(e^j) + \TT^j(e^j)\tau^j_*(e^j)\,,\label{eq:zw1}\\
\NN^i(e^i)\nu^i_*(e^i) + \TT^i(e^i)\tau^i_*(e^i) 
&=
\NN^k(e^k)\nu^k_*(e^k) + \TT^k(e^k)\tau^k_*(e^k)\,.\label{eq:zw2}
\end{align}
Multiplying \eqref{eq:zw1} by $\tau^i_*(e^i)$ we get
\begin{equation}\label{eq:zz}
    \TT^i(e^i) = \NN^j(e^j)\beta_m^{ij} + \TT^j(e^j)\alpha^{ij}_m\,.
\end{equation}
Hence analogously
\begin{equation}\label{eq:zz1}
    \begin{split}
        \TT^j(e^j)& = \NN^k(e^k)\beta_m^{jk} + \TT^k(e^k)\alpha^{jk}_m\,, \\
        \TT^k(e^k) &= \NN^i(e^i)\beta_m^{ki} + \TT^i(e^i)\alpha^{ki}_m\,.
    \end{split}
\end{equation}
Combining \eqref{eq:zz} and \eqref{eq:zz1} we obtain
\[
\begin{split}
     \TT^i(e^i) &= \NN^j(e^j)\beta_m^{ij} + \alpha^{ij}_m \left[ \NN^k(e^k)\beta_m^{jk} + \TT^k(e^k)\alpha^{jk}_m\right] \\
     &= \NN^j(e^j)\beta_m^{ij} + \alpha^{ij}_m \left[ \NN^k(e^k)\beta_m^{jk} + \alpha^{jk}_m\left( \NN^i(e^i)\beta_m^{ki} + \TT^i(e^i)\alpha^{ki}_m \right) \right],\\
      &= \NN^j(e^j)\beta_m^{ij} +  \NN^k(e^k)\alpha^{ij}_m\beta_m^{jk} +  \NN^i(e^i)\alpha^{ij}_m\alpha^{jk}_m\beta_m^{ki} + \TT^i(e^i)\alpha^{ij}_m\alpha^{jk}_m\alpha^{ki}_m  ,
\end{split}
\]
that gives the first line in \eqref{eq:BoundaryRelationsTN}. The remaining identities in \eqref{eq:BoundaryRelationsTN} follow analogously.

Denoting by $p=\Gamma(m)= \gamma^i(e^i)=\gamma^j(e^j)=\gamma^k(e^k)$ the image of the junction, since
\[
(-1)^{e_i}\nu^i_*(e^i) + (-1)^{e_j}\nu^j_*(e^j) + (-1)^{e_k}\nu^k_*(e^k) =0,
\]
we get
\[
\begin{split}
    0 &= \scal{p ,(-1)^{e_i}\nu^i_*(e^i) + (-1)^{e_j}\nu^j_*(e^j) + (-1)^{e_k}\nu^k_*(e^k) } \\
    &= \scal{\gamma^i(e^i) , (-1)^{e_i}\nu^i_*(e^i)} + \scal{\gamma^j(e^j) , (-1)^{e_j}\nu^j_*(e^j)} + \scal{\gamma^k(e^k) , (-1)^{e_k}\nu^k_*(e^k)} \\
    &= \scal{\Gamma_*(m) ,(-1)^{e_i}\nu^i_*(e^i) + (-1)^{e_j}\nu^j_*(e^j) + (-1)^{e_k}\nu^k_*(e^k) } + (-1)^{e^i}\NN^i(e^i) + (-1)^{e^j}\NN^j(e^j) + (-1)^{e^k}\NN^k(e^k) \\
    &= (-1)^{e^i}\NN^i(e^i) + (-1)^{e^j}\NN^j(e^j) + (-1)^{e^k}\NN^k(e^k),
\end{split}
\]
that is \eqref{eq:SommaNN}.
Now plugging~\eqref{eq:SommaNN} into~\eqref{eq:BoundaryRelationsTN} readily implies \eqref{eq:LineariNiNj}.

%\begin{equation}\label{eq:LineariNiNj}
%\begin{split}
 %   \TT^i(e^i) &= 
  %  -\frac{1}{\sqrt{3}}\NN^i(e^i)
   % -\frac{2}{\sqrt{3}}(-1)^{e_i+e_j}\NN^j(e^j) \,,\\
%     \TT^j(e^j) &= 
 %   \frac{2}{\sqrt{3}}(-1)^{e_i+e_j}\NN^i(e^i)   +\frac{1}{\sqrt{3}}\NN^j(e^j)\,,\\
  %   \TT^k(e^k) &= 
   %  -\frac{1}{\sqrt{3}}(-1)^{e_i+e_k}\NN^i(e^i)
    % +\frac{1}{\sqrt{3}}(-1)^{e_j+e_k}\NN^j(e^j)\,.
%\end{split}
%\end{equation}
\end{proof}

If now $\Gamma_*:G\to \R^2$ is regular, in the next lemma we state the sufficient conditions that functions $\NN^\ell, \TT^\ell$ may satisfy to define a triple junctions network as a graph over $\Gamma_*$.

\begin{lemma}\label{lem:SufficientConditionsNetwork}
Let $\Gamma_*:G\to \R^2$ be a regular network. 
Then there exists $\eps_{\Gamma_*}>0$ such that 
for every $\NN^\ell,\TT^\ell \in C^1([0,1])$, 
with $\|\NN^\ell\|_{C^1}, \|\TT^\ell\|_{C^1} \le \eps_{\Gamma_*}$ fulfilling at any junction $m=\pi(e^i,i)=\pi(e^j,j)=\pi(e^k,k)$
the identities  
\[
\begin{split}
    \TT^i(e^i) 
     &=
     \frac{1}{1-\alpha^{ij}_m\alpha^{jk}_m\alpha^{ki}_m}    \left(    \NN^i(e^i)\left(\alpha^{ij}_m\alpha^{jk}_m\beta_m^{ki} +(-1)^{1+e_i+e_k}\alpha^{ij}_m\beta_m^{jk} \right)+  \NN^j(e^j)\left(\beta_m^{ij}       +(-1)^{1+e_j+e_k}\alpha^{ij}_m\beta_m^{jk}\right)\right)  ,\\
          \TT^j(e^j)     &=
    \frac{1}{1-\alpha^{ij}_m\alpha^{jk}_m\alpha^{ki}_m}    \left(   \NN^i(e^i)\left(\alpha^{jk}_m\beta_m^{ki} 
      +(-1)^{1+e_i+e_k}\beta_m^{jk}
     \right)+  \NN^j(e^j)\left(\alpha^{jk}_m\alpha^{ki}_m\beta_m^{ij}  +(-1)^{1+e_j+e_k}\beta_m^{jk}\right)\right)  ,\\
                \TT^k(e^k) 
     &=
     \frac{1}{1-\alpha^{ij}_m\alpha^{jk}_m\alpha^{ki}_m}
     \left(\NN^i(e^i)
     \left(\beta_m^{ki} +(-1)^{1+e_i+e_k} \alpha^{ki}_m\alpha^{ij}_m\beta_m^{jk} \right)+ \NN^j(e^j)\left(\alpha^{ki}_m\beta_m^{ij}+(-1)^{1+e_j+e_k} \alpha^{ki}_m\alpha^{ij}_m\beta_m^{jk}\right)    \right) ,
\end{split}
\]
% \begin{equation}\label{eq:ww}
% \begin{split}
% &\TT^i(e^i) = 
% -\frac{1}{\sqrt{3}}\NN^i(e^i)
% -\frac{2}{\sqrt{3}}(-1)^{e_i+e_j}\NN^j(e^j) \,,\\
% &\TT^j(e^j) = 
% \frac{2}{\sqrt{3}}(-1)^{e_i+e_j}\NN^i(e^i)   +\frac{1}{\sqrt{3}}\NN^j(e^j)\,,\\
% &\TT^k(e^k) = 
% -\frac{1}{\sqrt{3}}(-1)^{e_i+e_k}\NN^i(e^i)
% +\frac{1}{\sqrt{3}}(-1)^{e_j+e_k}\NN^j(e^j)\,,
% \end{split}
%\end{equation}
\[
(-1)^{e^i}\NN^i(e^i) + (-1)^{e^j}\NN^j(e^j) + (-1)^{e^k}\NN^k(e^k) = 0\,,
\]
the maps
\begin{equation*}
    \gamma^\ell(x) \eqdef \gamma^\ell_*(x) + \NN^\ell(x)\nu^\ell_*(x) + \TT^\ell(x)\tau^\ell_*(x),
\end{equation*}
define a triple junctions network.
\end{lemma}

\begin{proof}
It is sufficient to check that whenever $m=\pi(e^i,i)=\pi(e^j,j)=\pi(e^k,k)$ is a junction, then $\gamma^i(e^i)=\gamma^j(e^j)=\gamma^k(e^k)$, namely, we have to check that the three vectors 
$\NN^i(e^i)\nu^i_*(e^i) + \TT^i(e^i)\tau^i_*(e^i)$, 
$\NN^j(e^j)\nu^j_*(e^j) + \TT^j(e^j)\tau^j_*(e^j)$ and 
$\NN^k(e^k)\nu^k_*(e^k) + \TT^k(e^k)\tau^k_*(e^k)$ coincide.
To this aim, observe that the identities in the assumptions imply that $\TT^i(e^i), \TT^j(e^j), \TT^k(e^k)$ satisfy \eqref{eq:BoundaryRelationsTN}. Taking scalar products of $\NN^j(e^j)\nu^j_*(e^j) + \TT^j(e^j)\tau^j_*(e^j)$ and $\NN^k(e^k)\nu^k_*(e^k) + \TT^k(e^k)\tau^k_*(e^k)$ with $\tau^i_*(e^i)$ and $\nu^i_*(e^i)$, exploiting \eqref{eq:BoundaryRelationsTN} one easily checks that
\[
\begin{split}
\scal{\NN^j(e^j)\nu^j_*(e^j) &+ \TT^j(e^j)\tau^j_*(e^j), \tau^i_*(e^i)} =
    \scal{\NN^k(e^k)\nu^k_*(e^k) + \TT^k(e^k)\tau^k_*(e^k), \tau^i_*(e^i)}  = \TT^i(e^i),\\
    \scal{\NN^j(e^j)\nu^j_*(e^j) &+ \TT^j(e^j)\tau^j_*(e^j), \nu^i_*(e^i)} =
    \scal{\NN^k(e^k)\nu^k_*(e^k) + \TT^k(e^k)\tau^k_*(e^k), \nu^i_*(e^i)}  = \NN^i(e^i),
\end{split}
\]
\end{proof}

The previous \cref{lem:BoundaryRelationsTN} and \cref{lem:SufficientConditionsNetwork} motivate the following definition.

\begin{defn}\label{def:LinearOperators}
Let $G$ be a regular graph. At any junction $m=\pi(e^i,i)=\pi(e^j,j)=\pi(e^k,k)$, if $i<j<k$, we denote by $L^i,L^j,L^k$ the linear maps
\begin{equation*}
\begin{split}
L^i(a,b)
     &=
     \frac{1}{1-\alpha^{ij}_m\alpha^{jk}_m\alpha^{ki}_m}    \left(    \left(\alpha^{ij}_m\alpha^{jk}_m\beta_m^{ki} +(-1)^{1+e_i+e_k}\alpha^{ij}_m\beta_m^{jk} \right) a +  \left(\beta_m^{ij}       +(-1)^{1+e_j+e_k}\alpha^{ij}_m\beta_m^{jk}\right)b\right)  ,\\
L^j(a,b)    &=
    \frac{1}{1-\alpha^{ij}_m\alpha^{jk}_m\alpha^{ki}_m}    \left( \left(\alpha^{jk}_m\beta_m^{ki} 
      +(-1)^{1+e_i+e_k}\beta_m^{jk}
     \right)a+  \left(\alpha^{jk}_m\alpha^{ki}_m\beta_m^{ij}  +(-1)^{1+e_j+e_k}\beta_m^{jk}\right) b\right)  ,
    \\
L^k(a,b)  
     &=
     \frac{1}{1-\alpha^{ij}_m\alpha^{jk}_m\alpha^{ki}_m}
     \left(
     \left(\beta_m^{ki} +(-1)^{1+e_i+e_k} \alpha^{ki}_m\alpha^{ij}_m\beta_m^{jk} \right)a+ \left(\alpha^{ki}_m\beta_m^{ij}+(-1)^{1+e_j+e_k} \alpha^{ki}_m\alpha^{ij}_m\beta_m^{jk}\right)b    \right),
%\\
%     &L^i(a,b) = -\frac{1}{\sqrt{3}}a
%     -\frac{2}{\sqrt{3}}(-1)^{e_i+e_j}b\,,\\
%      &L^j(a,b) = \frac{2}{\sqrt{3}}(-1)^{e_i+e_j}a + 
%      \frac{1}{\sqrt{3}}b\,,\\
%      &L^k(a,b) =-\frac{1}{\sqrt{3}}(-1)^{e_i+e_k}a + \frac{1}{\sqrt{3}}(-1)^{e_j+e_k}b\,.
\end{split}
\end{equation*}
for any $(a,b)\in\R^2$.
Moreover, we denote by $I_m$ the set of indices $\ell$ such that $E^\ell$ has an endpoint at $m$, and we denote by $e^\ell_m\in\{0,1\}$ the endpoint of $E^\ell$ at $m$, for $\ell \in I_m$.\\
Furthermore, for $\ell \in I_m$, we denote by $\Ell^\ell_m$ the linear operator $L^i,L^j,L^k$, depending on whether $\ell$ is the minimal, intermediate, or maximal index $\ell \in I_m$.\\
Finally, for any endpoint $p$, we denote by $i_p$ the corresponding index such that $p$ is an endpoint of $E^{i_p}$.
\end{defn}

We are ready for proving the existence of a canonical graph parametrization of triple junctions networks over regular ones that are close in $H^2$-norm. As discussed in the introduction, we shall perform the construction fixing a dependence of the tangential component of the graph with respect to the normal one. Such dependence is naturally defined by suitably extending to the interior of the edges the relations we found holding at junctions in \cref{lem:BoundaryRelationsTN} and \cref{lem:SufficientConditionsNetwork}.

For this purpose, from now on and for the rest of the paper, we fix a nonincreasing smooth cut-off function
\begin{equation*}
    \chi:[0,1/2]\to [0,1], 
    \quad 
    \chi|_{[0,\frac18]}=1,
    \quad
    \chi|_{[\frac38,\frac12]}\equiv 0.
\end{equation*}

\begin{prop}\label{prop:ParametrizzazioneTN}
Let $\Gamma_*:G\to \R^2$ be a smooth regular network. Then there exists $\eps_{\Gamma_*}>0$ such that whenever $\Gamma:G\to \R^2$ is a triple junctions network of class $H^2$ such that
\begin{align}
    \label{eq:aa1} \sum_i \| \gamma^i_* - \gamma^i \|_{H^2(\d x)} \le \eps_{\Gamma_*},& \\
    \gamma^i_*(p)=\gamma^i(p) &\qquad \text{$\forall\,p\in G\st p$ is an endpoint,} \label{eq:EndpointFissi}
\end{align}
for any $i$, then there exist functions $\NN^i, \TT^i \in H^2(\d x)$
and reparametrizations $\varphi^i:[0,1]\to[0,1]$ of class $H^2(\d x)$ such that
\begin{equation*}
     \gamma^i\circ \varphi^i(x) = \gamma^i_*(x) + \NN^i(x)\nu^i_*(x) + \TT^i(x)\tau^i_*(x).
\end{equation*}

At any junction $\pi(e^i,i)=\pi(e^j,j)=\pi(e^k,k)$, where $i<j<k$, there holds
\begin{equation}\label{eq:TinFunzDiN}
\begin{split}
\TT^i(|e^i-x|) &= \chi(x) L^i(\NN^i (|e^i-x|), \NN^j(|e^j-x|) ) , \\
\TT^j(|e^j-x|) &= \chi(x) L^j(\NN^i(|e^i-x|) , \NN^j(|e^j-x|)), \\
\TT^k(|e^k-x|) &= \chi(x) L^k(\NN^i(|e^i-x|), \NN^j(|e^j-x|) ) ,
\end{split}
\end{equation}
for $x \in [0,\tfrac12]$.

If $\pi(1,i)$ is an endpoint, then
\begin{equation}\label{eq:TinFunzDiN2}
    \TT^i(x)=0,
\end{equation}
for $x \in [\tfrac12,1]$.

Moreover
\begin{itemize}
    \item for any $\delta>0$ there is $\eps \in (0,\eps_{\Gamma_*})$ such that
    \begin{equation}\label{eq:StimaH2Ni}
        \sum_i \| \gamma^i_* - \gamma^i \|_{H^2(\d x)} \le \eps
    \quad\implies\quad  
    \sum_i \|\NN^i\|_{H^2(\d x)}+ \|\varphi^i(x)-x\|_{H^2(\d x)} \le \delta;
    \end{equation}
    
    \item for any $\eta>0$ and $m\in \N$ there is $\eps_{\eta,m}\in(0,\eps_{\Gamma_*})$ such that if $\sum_i \| \gamma^i_* - \gamma^i \|_{C^{m+1}([0,1])} \le \eps_{\eta,m}$, then
    \begin{equation}\label{eq:StimaSobolevNi}
    \sum_i \| \NN^i\|_{H^m(\d x)} \le \eta.
    \end{equation}
\end{itemize}
\end{prop}

\begin{proof}
Without loss of generality, we can perform the construction at a junction $\pi(e^i,i)=\pi(e^j,j)=\pi(e^k,k)$, where $i<j<k$, assuming $e^i=e^j=e^k=0$. For the sake of clarity, we show how to construct $\varphi^i,\NN^i$ and  $\varphi^j,\NN^j$ on $[0,\frac12]$ only, the complete proof being a straightforward adaptation.

Consider the function $F:[0,\frac12]\times\R^2\times\R^2\to \R^4$ given by
\[
F(x,n^i,y^i,n^j,y^j)\eqdef 
\begin{pmatrix}
\gamma^i_*(x) + n^i\nu^i_*(x) + \chi(x) L^i(n^i,n^j) \tau^i_*(x) - \gamma^i(y^i) \\
\gamma^j_*(x) + n^j\nu^j_*(x) + \chi(x) L^j(n^i,n^j) \tau^j_*(x) - \gamma^j(y^j) 
\end{pmatrix}
\]
Reflecting by symmetry, we can assume that $F$ is also defined in an open neighborhood of $x=0$. Since $\{\tau^i_*(0),\nu^i_*(0)\}$ (and $\{\tau^j_*(0),\nu^j_*(0)\}$) is a basis of $\R^2$, there exist unique numbers $\NN^i(0),\TT^i(0)$ (and $\NN^j(0),\TT^j(0)$) such that $\gamma^i(0) = \gamma^i_*(0) + \NN^i(0)\nu^i_*(0) + \TT^i(0)\tau^i_*(0)$ (and $\gamma^j(0) = \gamma^j_*(0) + \NN^j(0)\nu^j_*(0) + \TT^j(0)\tau^j_*(0)$). Since $\Gamma$ is a triple junctions network, by \cref{lem:BoundaryRelationsTN}, we have that $F(0,\NN^i(0),0,\NN^j(0),0)=(0)$. Moreover, the matrix
\[
\begin{split}
M(x,n^i,&y^i,n^j,y^j) \eqdef 
\begin{pmatrix}[c|c|c|c]
\partial_{n^i} F & \partial_{y^i} F & \partial_{n^j} F & \partial_{y^j} F
\end{pmatrix} \big|_{(x,n^i,y^i,n^j,y^j)} \\
&=
\begin{pmatrix}[c|c|c|c]
\nu^i_*(x) + \chi(x) L^i(1,0)\tau^i_*(x) & -(\gamma^i)'(y^i)&
\chi(x)L^i(0,1)\tau^i_*(x)
& 0
\\
\chi(x)L^j(1,0)\tau^j_*(x) & 0 & \nu^j_*(x) + \chi(x) L^j(0,1)\tau^j_*(x) & -(\gamma^j)'(y^j)
\end{pmatrix}
\end{split}
\]
satisfies
\[
M(0,\NN^i(0),0,\NN^j(0),0)=
\begin{pmatrix}[c|c|c|c]
\nu^i_*(0) +  L^i(1,0)\tau^i_*(0) & -(\gamma^i)'(0)&
L^i(0,1)\tau^i_*(0)
& 0
\\
L^j(1,0)\tau^j_*(0) & 0 & \nu^j_*(0) +  L^j(0,1)\tau^j_*(0) & -(\gamma^j)'(0)
\end{pmatrix}.
\]
It is readily checked that any two columns are linearly independent. Hence we can apply the implicit function theorem to get the existence of $\varphi^i,\NN^i$ and  $\varphi^j,\NN^j$ defined on some interval $[0,\xi]\subset[0,1/2]$ such that
\[
\begin{split}
    \gamma^i(\varphi^i(x)) &= \gamma^i_*(x) + \NN^i(x)\nu^i_*(x) + \chi(x) L^i(\NN^i(x),\NN^j(x)) \tau^i_*(x)  ,\\
 \gamma^j(\varphi^j(x)) &= \gamma^j_*(x) + \NN^j(x)\nu^j_*(x) + \chi(x) L^j(\NN^i(x),\NN^j(x)) \tau^j_*(x) .
\end{split}
\]
From the identity
\begin{equation}\label{eq:DerivataFunzImplicita}
\begin{split}
   & \begin{pmatrix}
\partial_x \NN^i(x) \\
\partial_x \varphi^i(x) \\
\partial_x \NN^j(x) \\
\partial_x \varphi^j(x) \\
\end{pmatrix} 
=
-\left[ M\left(x, \NN^i(x),\varphi^i(x), \NN^j(x), \varphi^j(x)\right)\right]^{-1} \cdot\partial_x F\big|_{(x, \NN^i(x),\varphi^i(x), \NN^j(x), \varphi^j(x))} \\
&\quad = -\left[ M\left(x, \NN^i(x),\varphi^i(x), \NN^j(x), \varphi^j(x)\right)\right]^{-1} 
 \cdot
 \begin{pmatrix}
 (\gamma^i_*)' + \NN^i \partial_x \nu^i_* + [\chi' \tau^i_* + \chi \partial_x \tau^i_*] L^i(\NN^i, \NN^j) \\
 (\gamma^j_*)' + \NN^j \partial_x \nu^j_* + [\chi' \tau^j_* + \chi \partial_x \tau^j_*] L^j(\NN^i, \NN^j)
 \end{pmatrix}
\end{split}
\end{equation}
we estimate
\begin{equation*}
\begin{split}
       \| \partial_x \NN^i(x) \|_{L^\infty(0,\xi)} &+
\|\partial_x \varphi^i(x) \|_{L^\infty(0,\xi)} +
\|\partial_x \NN^j(x) \|_{L^\infty(0,\xi)} +
\|\partial_x \varphi^j(x)\|_{L^\infty(0,\xi)} \\
&\le C(\Gamma_*, \|\partial_x \gamma^i\|_\infty,  \|\partial_x \gamma^i\|_\infty) \left( 1 + \| \NN^i(x) \|_{L^\infty(0,\xi)} + \| \NN^j(x) \|_{L^\infty(0,\xi)}\right) \\
& \stackrel{\eqref{eq:aa1}}{\le} C(\Gamma_*, \eps_{\Gamma_*}) \left( 1 + \| \NN^i(x) \|_{L^\infty(0,\xi)} + \| \NN^j(x) \|_{L^\infty(0,\xi)}\right).
\end{split}
\end{equation*}
Since $\NN^i = \scal{\gamma^i\circ \varphi^i -\gamma^i_*, \nu^i_*}$, and analogously for $j$, recalling \eqref{eq:aa1} we get
\begin{equation}\label{eq:yyy}
    \begin{split}
        \| \partial_x \NN^i(x) \|_{L^\infty(0,\xi)} &+
\|\partial_x \varphi^i(x) \|_{L^\infty(0,\xi)} +
\|\partial_x \NN^j(x) \|_{L^\infty(0,\xi)} +
\|\partial_x \varphi^j(x)\|_{L^\infty(0,\xi)} \le C(\Gamma_*, \eps_{\Gamma_*}).
    \end{split}
\end{equation}

We claim that
\begin{equation}\label{eq:Claim}
\begin{split}
    \forall\,\delta>0 \,\exists\,\eps \in(0,\eps_{\Gamma_*}) \st \sum_i\| \gamma^i_* - \gamma^i \|_{H^2(\d x)} \le \eps \quad\implies\quad & \exists\, \overline{\xi}(\eps)>0\,\text{ st }\\
    &\xi\ge \overline{\xi}, \\ &\|\varphi^i(x)-x\|_{W^{1,\infty}(0,\overline\xi)} \le \delta.
\end{split}
\end{equation}
Indeed, suppose by contradiction that there is $\delta>0$ and a sequence of triple junctions networks $\Gamma_n$ such that $\sum_i\| \gamma^i_* - \gamma^i_n \|_{H^2(\d x)} \le 1/n $, but the implicit functions $\varphi^i_n:[0,\xi_n]\to[0,1/2]$ obtained as above do not verify \eqref{eq:Claim}.
Denote by $F_n,M_n, \NN^i_n,\NN^j_n$ the map, the matrices, and the functions defined by the above procedure applied on the network $\Gamma_n$ in place of $\Gamma$.
Since $\sum_i\| \gamma^i_* - \gamma^i_n \|_{H^2(\d x)} \le 1/n $, there are $S,\rho>0$ independent of $n$ such that
\[
\|[M_n(0,\NN_n^i(0),0,\NN_n^j(0),0)]^{-1}]\|\le S,
\]
\[
\left\|{\rm id} - [M_n(0,\NN_n^i(0),0,\NN_n^j(0),0)]^{-1} M_n(x,n^i,y^i,n^j,y^j) \right\| \le \frac12 ,
\]
whenever $|x|<\rho$ and $|(n^i,y^i,n^j,y^j) - (\NN_n^i(0),0,\NN_n^j(0),0)|<\rho$. Furthermore, since $\partial_x F_n (x,n^i,y^i,n^j,y^j)$ does not depend on $n$, there is $N>0$ such that
\[
\|\partial_x F_n (x,n^i,y^i,n^j,y^j)\| \le N,
\]
whenever $|x|<\rho$ and $|(n^i,y^i,n^j,y^j) - (\NN_n^i(0),0,\NN_n^j(0),0)|<\rho$, for any $n$. Hence the assumptions of \cref{thm:ImplicitFunction} are satisfied, and thus there is $\overline\xi>0$ such that $\xi_n\ge\overline\xi$ for any $n$. Then it must be that $\|\varphi^i(x)-x\|_{W^{1,\infty}(0,\overline\xi)} > \delta$.

Up to subsequence, recalling the uniform bounds \eqref{eq:yyy}, we can pass to the limit $n\to \infty$ in the identity
\[
\gamma^i_n(\varphi^i_n(x)) = \gamma^i_*(x) + \NN^i_n(x)\nu^i_*(x) + \chi(x) L^i(\NN^i_n(x),\NN^j_n(x)) \tau^i_*(x),
\]
to obtain
\begin{equation}\label{eq:yyy2}
\gamma^i_*(\varphi^i_\infty(x)) = \gamma^i_*(x) + \NN^i_\infty(x)\nu^i_*(x) + \chi(x) L^i(\NN^i_\infty(x),\NN^j_\infty(x)) \tau^i_*(x),
\end{equation}
where $\varphi^i_n\to \varphi^i_\infty$, $\NN^i_n\to \NN^i_\infty$, and $\NN^j_n\to \NN^j_\infty$ in $C^0([0,\overline\xi])$ and in $H^1(0,\overline\xi)$, and \eqref{eq:yyy2} holds pointwise on $[0,\overline\xi]$. By the uniqueness part of the implicit function theorem, we deduce that $\varphi^i_\infty(x)\equiv x$, $\NN^i_\infty(x)\equiv0$, and $\NN^j_\infty(x)\equiv0$.\\
Moreover, by \eqref{eq:DerivataFunzImplicita}, uniform convergence on the right hand side implies that $\varphi^i_n(x)\to x$, $\NN^i_n\to 0$, and $\NN^j_n\to 0$ in $C^1([0,\overline\xi])$.\\
Hence $0=\|\varphi^i_\infty(x)-x\|_{W^{1,\infty}(0,\overline\xi)} =\lim_n\|\varphi^i_n(x)-x\|_{W^{1,\infty}(0,\overline\xi)} >\delta$ gives a contradiction, and \eqref{eq:Claim} follows.

By \eqref{eq:Claim}, up to decreasing $\eps_{\Gamma_*}$, then $(\varphi^i)'\ge \tfrac12$ on $(0,\overline\xi)$ for any $i$. Hence further differentiating \eqref{eq:DerivataFunzImplicita} and arguing as before, we also derive that
\begin{equation}\label{eq:Claim2}
\begin{split}
    \forall\,\delta>0 \,\exists\,\eps \in(0,\eps_{\Gamma_*}) &\st \sum_i \| \gamma^i_* - \gamma^i \|_{H^2(\d x)} \le \eps \\
    &\quad\implies\quad  
    \exists\, \overline{\xi}(\eps)>0\,\text{ st }\,\,
    \xi\ge \overline{\xi}, \quad\|\NN^i\|_{H^2(0,\overline\xi)}+ \|\varphi^i(x)-x\|_{H^2(0,\overline\xi)} \le \delta.
\end{split}
\end{equation}
Since $\overline{\xi}$ only depends on $\eps_{\Gamma_*}$, we can iterate finitely many times the above argument to get the complete construction on the interval $[0,\frac12]$ as claimed.\\
Moreover, \eqref{eq:Claim2} eventually implies \eqref{eq:StimaH2Ni}.
Similarly, further differentiating \eqref{eq:DerivataFunzImplicita} leads to \eqref{eq:StimaSobolevNi}.
\end{proof}

Arguing as in \cref{prop:ParametrizzazioneTN}, one obtains the following analogous consequence for the parametrization of a time dependent family of networks in a neighborhood of a fixed one $\Gamma_*$.

\begin{cor}\label{cor:ParametrizzazioneTNtempo}
Let $\Gamma_*:G\to \R^2$ be a smooth regular network. Then there exist $\eps_{\Gamma_*}>0$ such that whenever $\Gamma_t:G\to \R^2$ is a one-parameter family of triple junctions networks of class $H^2$, differentiable with respect to $t$ for $t \in [t_0-h,t_0+h]$ and $h>0$, such that $(t,x)\mapsto \partial_t \gamma^i_t(x)$ is continuous for any $i$ and
\begin{align*}
     \sum_i \| \gamma^i_* - \gamma^i_t \|_{H^2(\d x)}
     \le \eps_{\Gamma_*},& \\
     \gamma^i_*(p)=\gamma^i_t(p) &\qquad \text{$\forall\,p\in G\st p$ is an endpoint,} 
\end{align*}
for any $i,t$, then there exist $h'\in(0,h)$ and functions $\NN^i_t, \TT^i_t \in H^2(\d x)$ and reparametrizations $\varphi^i_t:[0,1]\to[0,1]$ of class $H^2(\d x)$, continuously differentiable with respect to $t$ for $t \in [t_0-h',t_0+h']$ such that
\begin{equation*}
      \gamma^i_t\circ \varphi^i_t(x) = \gamma^i_*(x) + \NN^i_t(x)\nu^i_*(x) + \TT^i_t(x)\tau^i_*(x).
\end{equation*}

At any junction $\pi(e^i,i)=\pi(e^j,j)=\pi(e^k,k)$, where $i<j<k$, there holds
\begin{equation*}
\begin{split}
\TT^i_t(|e^i-x|) &= \chi(x) L^i(\NN^i_t (|e^i-x|), \NN^j_t(|e^j-x|) ) , \\
\TT^j_t(|e^j-x|) &= \chi(x) L^j(\NN^i_t(|e^i-x|) , \NN^j_t(|e^j-x|)), \\
\TT^k_t(|e^k-x|) &= \chi(x) L^k(\NN^i_t(|e^i-x|), \NN^j_t(|e^j-x|) ) ,
\end{split}
\end{equation*}
for $x \in [0,\tfrac12]$.

If $\pi(1,i)$ is an endpoint, then
\begin{equation*}
     \TT^i_t(x)=0,
\end{equation*}
for $x \in [\tfrac12,1]$.

Moreover
\begin{itemize}
    \item for any $\delta>0$ there is $\eps \in (0,\eps_{\Gamma_*})$ such that
    \begin{equation}\label{eq:StimaH2NiTempo}
        \sum_i \| \gamma^i_* - \gamma^i_t \|_{H^2(\d x)} \le \eps \quad\forall\,t
    \quad\implies\quad  
    \sum_i \|\NN^i_t\|_{H^2(\d x)}+ \|\varphi^i_t(x)-x\|_{H^2(\d x)} \le \delta,
    \end{equation}
    for any $t \in [t_0-h',t_0+h']$;
    
    \item for any $\eta>0$ and $m\in \N$ there is $\eps_{\eta,m}\in(0,\eps_{\Gamma_*})$ and $h_{\eta,m} \in (0,h)$ such that if $\sum_i \| \gamma^i_* - \gamma^i_t \|_{C^{m+1}([0,1])} \le \eps$ for any $t$, then
    \begin{equation}\label{eq:StimaSobolevNiTempo}
    \sum_i \| \NN^i_t\|_{H^m(\d x)} \le \eta,
    \end{equation}
    for any $t \in [t_0-h_{\eta,m}, t_0 + h_{\eta,m}]$.
    \end{itemize}
\end{cor}

The construction of the ``tangent functions'' $\TT^i$'s in \cref{prop:ParametrizzazioneTN} and \cref{cor:ParametrizzazioneTNtempo} depending on the ``normal functions'' $\NN^i$'s motivates the next definition.

\begin{defn}[Adapted tangent functions]\label{def:Adapted}
Let $G$ be a regular graph. Let $\NN^i,\TT^i:[0,1]\to \R$ be functions of class $C^1$, for $i=1,\ldots,N$. We say that the $\TT^i$'s are \emph{adapted} to the $\NN^i$'s whenever there hold the relations \eqref{eq:TinFunzDiN} and \eqref{eq:TinFunzDiN2}.\\
More explicitly, the $\TT^i$'s are adapted to the $\NN^i$'s whenever
\begin{equation}\label{eq:DefAdattati}
    \TT^\ell(|e^\ell_m - x|) = \chi(x) \Ell^\ell_m(\NN^i(|e^i_m - x|), \NN^j(|e^j_m-x|) ),
\end{equation}
for $x \in [0,\tfrac12]$ for any junction $m=\pi(e^i,i)=\pi(e^j,j)=\pi(e^k,k)$ with $i<j<k$, and
\begin{equation*}
    \TT^i(x)=0
\end{equation*}
for $x \in [\tfrac12,1]$ for any endpoint $\pi(1,i)$.
\end{defn}

\subsection{First and second variations}

In order to derive the desired \L ojasiewicz--Simon inequality, we need to compute first and second variations of the length functional taking variations determined by graph parametrizations over regular networks with tangent functions adapted to normal functions as in~\cref{def:Adapted}.

\begin{prop}\label{prop:FirstVariation}
Let $\Gamma_*:G\to \R^2$ be a smooth regular network. Then there is $\eps_{\Gamma_*}>0$ such that the following holds.

Let $\NN^i, X^i \in H^2$ with $\|\NN^i\|_{H^2}\le \eps_{\Gamma_*}$ such that
\begin{equation*}
\sum_{\ell \in I_m} (-1)^{e^\ell_m} \NN^\ell(e^\ell_m) = \sum_{\ell \in I_m} (-1)^{e^\ell_m} X^\ell(e^\ell_m)=0
\qquad
\forall m \in J_G.
\end{equation*}
Let $\Gamma^\eps:G\to \R^2$ be the triple junctions network defined by
\begin{equation}\label{eq:zzz}
\gamma^{i,\eps}(x) \eqdef \gamma^i_*(x) + (\NN^i(x)+ \eps X^i(x))\nu^i_*(x) + \TT^{i,\eps}(x)\tau^i_*(x),
\end{equation}
for any $i$, for any $|\eps|<\eps_0$ and some $\eps_0>0$, where the $\TT^{i,\eps}$'s are adapted to the $(\NN^i+ \eps X^i)$'s, for any $|\eps|<\eps_0$\footnote{Immersions $\gamma^{i,\eps}$ define a triple junctions network by \cref{lem:SufficientConditionsNetwork} for $\eps_{\Gamma_*}$ small enough.}.
%Assume further that $\gamma^{i,\eps}(p)=\gamma^i_*(p) $ for any $i$ at any endpoint $p$.

Call $\Gamma$ the network given by the immersions $\gamma^i\eqdef\gamma^{i,0}$.
%
% Let $\Gamma^\eps:G\to \R^2$ be a triple junctions network of class $H^2$ of the form
% \begin{equation}\label{eq:zzz}
% \gamma^{i,\eps}(x) \eqdef \gamma^i_*(x) + (\NN^i(x)+ \eps X^i(x))\nu^i_*(x) + \TT^{i,\eps}(x)\tau^i_*(x),
% \end{equation}
% for any $i$, for any $|\eps|<\eps_0$ and some $\eps_0>0$, and call $\Gamma$ the network given by the immersions $\gamma^i\eqdef\gamma^{i,0}$. Assume that the $\TT^{i,\eps}$'s satisfy the relations in \eqref{eq:TinFunzDiN} with respect to the $(\NN^i+ \eps X^i)$'s, for any $|\eps|<\eps_0$, and that $\gamma^{i,\eps}(p)=\gamma^i_*(p) $ for any $i$ at any endpoint $p$.
%
Then
\begin{equation}\label{eq:FirstVariation}
\begin{split}
    \frac{\d}{\d \eps}  {\rm L} (\Gamma^\eps) \bigg|_0 
    &= 
    \sum_{p \in P_G} \scal{\tau^{i_p}(1), \nu^{i_p}_*(1)} X^{i_p}(1) +
    \\& \phantom{=}+
    \sum_{m \in J_G}  \sum_{\ell \in I_m} 
    (-1)^{1+e^\ell_m} \left[ 
    \scal{\tau^\ell(e^\ell_m), \nu^\ell_*(e^\ell_m)} + \sum_{j \in I_m} h_{\ell j}\scal{\tau^j(e^j_m) , \tau^j_*(e^j_m)}
     \right] 
     X^\ell(e^\ell_m) + \\
    &\phantom{=} - \sum_i 
    \int_0^1 \bigg( \scal{\boldsymbol{k}^i,\nu^i_*}|\partial_x \gamma^i| 
    + \sum_j f_{ij}\chi \scal{\boldsymbol{k}^j,\tau^j_*}|\partial_x \gamma^j| +\\
    &\qquad\qquad\qquad
    + g_{ij}\chi(1-x) \scal{\boldsymbol{k}^j,\tau^j_*}(1-x)|\partial_x \gamma^j|(1-x) \bigg) X^i \de x,
\end{split}
\end{equation}
where $f_{ij}, g_{ij}, h_{\ell j} \in \R$ depend on the topology of $G$.

If also $\Gamma$ is regular and $\gamma^{i,\eps}(p)=\gamma^i_*(p) $ for any $i$ at any endpoint $p$, then
\begin{equation}\label{eq:FirstVariationRegular}
    \begin{split}
        \frac{\d}{\d \eps}  {\rm L} (\Gamma^\eps) \bigg|_0 
    &= 
    - \sum_i 
    \int_0^1 \bigg( \scal{\boldsymbol{k}^i,\nu^i_*}|\partial_x \gamma^i| 
    + \sum_j f_{ij}\chi \scal{\boldsymbol{k}^j,\tau^j_*}|\partial_x \gamma^j| +\\
    &\qquad\qquad\qquad
    + g_{ij}\chi(1-x) \scal{\boldsymbol{k}^j,\tau^j_*}(1-x)|\partial_x \gamma^j|(1-x) \bigg) X^i \de x.
    \end{split}
\end{equation}
\end{prop}

\begin{proof}
Let us assume first that there is a junction $m$ such that the functions $X^\ell$ appearing in~\eqref{eq:zzz} all vanish except for $\ell \in I_m$. Moreover, for $\ell \in I_m$, assume that $e^\ell_m=0$ and that $X^\ell$ has compact support in $[0,\tfrac58)$.

Let us denote $m=\pi(e^i,i)=\pi(e^j,j)=\pi(e^k,k)$, where $i<j<k$. By differentiating the length functional we get
\begin{equation}\label{eq:www}
    \begin{split}
        \frac{\d}{\d \eps} {\rm L} (\Gamma^\eps)
        =  \sum_{\ell \in I_m}\int_0^1 \frac{1}{|\partial_x \gamma^{\ell,\eps}|} \scal{\partial_x \gamma^{\ell,\eps}, \partial_x \partial_\eps\gamma^{\ell,\eps}} \de x,
    \end{split}
\end{equation}
indeed, since ${\rm spt}\,X^\ell \subset [0,\tfrac58)$, by \eqref{eq:DefAdattati} and definition of $\chi$, we have that $\TT^{n,\eps}$ does not depend on $\eps$ for all $n \not \in I_m$. Moreover
\begin{equation}\label{eq:zxzxzx}
    \begin{split}
    \TT^{\ell,\eps}(x) &= \chi(x) \Ell^\ell_m((\NN^i+ \eps X^i)(x), (\NN^j+\eps X^j)(x) ) \\
    &=  \chi(x) \Ell^\ell_m(\NN^i(x), \NN^j(x) ) + \eps \, \chi(x)\Ell^\ell_m(X^i(x), X^j(x) )
\end{split}
\end{equation}
for $\ell \in I_m$, hence, letting
\begin{equation}\label{eq:zww}
    Y^\ell \eqdef \partial_\eps\gamma^{\ell,\eps} = X^\ell\nu^\ell_* + \chi \Ell^\ell_m(X^i,X^j)\tau^\ell_*,
\end{equation}
we find
\begin{equation*}
    \begin{split}
        \frac{\d}{\d \eps} {\rm L} (\Gamma^\eps) \bigg|_0
        &= \sum_{\ell \in I_m} \int_0^1 \frac{1}{|\partial_x \gamma^\ell|} \scal{\partial_x \gamma^\ell, \partial_x Y^\ell} \de x 
        = \sum_{\ell \in I_m} \int_0^1  \scal{\tau^\ell, \partial_x Y^\ell} \de x  \\
        &=\sum_{\ell \in I_m}  - \scal{\tau^\ell(0),Y^\ell(0)} - \int_0^1  \scal{\boldsymbol{k}^\ell, Y^\ell} \de s^\ell.
    \end{split}
\end{equation*}
Since $\gamma^{\ell,\eps}(0)=\gamma^{l,\eps}(0)$ for any $\eps$ and $\ell,l \in I_m$, then $Y^\ell(0)=Y^l(0)$ for any $\ell,l \in I_m$. Hence, if $\Gamma$ is regular, then the boundary term $\sum_{\ell \in I_m} \scal{\tau^\ell(0),Y^\ell(0)}=0$.

Employing \eqref{eq:zww} we get
\begin{equation*}
    \begin{split}
        \frac{\d}{\d \eps} {\rm L} (\Gamma^\eps) \bigg|_0
        &= - \sum_{\ell \in I_m}   \scal{\tau^\ell(0),\nu^\ell_*(0)} X^\ell(0) + \Ell^\ell_m(X^i(0),X^j(0)) \scal{\tau^\ell(0), \tau^\ell_*(0)} +\\
        &\phantom{=}
        - \sum_{\ell \in I_m} \int_0^1  \scal{\boldsymbol{k}^\ell, \nu^\ell_*} X^\ell + \scal{\boldsymbol{k}^\ell, \tau^\ell_* }\chi \Ell^\ell_m(X^i,X^j) \de s^\ell.
    \end{split}
\end{equation*}
%Factoring out the functions $X^\ell$ in each term, the above formula takes the form given in \eqref{eq:FirstVariation} (and in \eqref{eq:FirstVariationRegular} in case $\Gamma$ is regular).

Suppose now that there is an endpoint $p \in P_G$ such that the functions $X^\ell$ appearing in \eqref{eq:zzz} all vanish except for $\ell = i_p$. Moreover, assume that $X^{i_p}$ has compact support in $(\tfrac38,1]$. Hence $Y^{i_p}\eqdef \partial_\eps \gamma^{i_p,\eps} = X^{i_p}\nu^{i_p}_*$ in this case, and the same computation performed above now yields
\begin{equation*}
\begin{split}
        \frac{\d}{\d \eps} {\rm L}(\Gamma^\eps)\bigg|_0 = \int_0^1  \scal{\tau^{i_p}, \partial_x Y^{i_p}} \de x  = \scal{\tau^{i_p}(1),Y^{i_p}(1)} - \int_0^1  \scal{\boldsymbol{k}^{i_p}, Y^{i_p}} \de s^{i_p},
\end{split}
\end{equation*}
which takes the form given in \eqref{eq:FirstVariation}. In case $\gamma^{i,\eps}(p)=\gamma^i_*(p) $ for any $i$ at any endpoint $p$, then $\NN^{i_p}(1)=X^{i_p}(1)=0$, and \eqref{eq:FirstVariationRegular} follows as well.

Considering now arbitrary variations as in \eqref{eq:zzz}, then \eqref{eq:FirstVariation} follows in the general case observing that the formula is linear with respect to the $X^i$'s and that each $X^i$ can be written as $X^i= \eta X^i + (1-\eta) X^i$ in a way that ${\rm spt} (\eta X^i) \subset [0,\tfrac58)$ and ${\rm spt} ((1-\eta) X^i )\subset (\tfrac38,1]$, recalling also that $\partial_\eps \gamma^{i,\eps}(p)=0$ at any endpoint $p$. Additive terms of the form $g_{ij}\chi(1-x) \scal{\boldsymbol{k}^j,\tau^j_*}(1-x)|\partial_x \gamma^j|(1-x)  X^i (x)$ appear by changing variables in order to factor out the function $X^i(x)$ in the $i$-th integral.
\end{proof}

\begin{prop}\label{prop:SecondVariationBilinear}
Let $\Gamma_*:G\to \R^2$ be a minimal network.
Let $ X^i, Z^i \in H^2$ such that
\begin{equation*}
\sum_{\ell \in I_m} (-1)^{e^\ell_m} X^\ell(e^\ell_m) = \sum_{\ell \in I_m} (-1)^{e^\ell_m} Z^\ell(e^\ell_m) =0
\qquad
\forall m \in J_G.
\end{equation*}
Let $\Gamma^{\eps,\eta}:G\to \R^2$ be the triple junctions network defined by
\begin{equation*}
\gamma^{i,\eps,\eta}(x) \eqdef \gamma^i_*(x) + (\eps X^i(x) + \eta Z^i(x))\nu^i_*(x) + \TT^{i,\eps, \eta}(x)\tau^i_*(x),
\end{equation*}
for any $i$, for any $|\eps|, |\eta|<\eps_0$ and some $\eps_0>0$, where the $\TT^{i,\eps,\eta}$'s are adapted to the $(\eps X^i+\eta Z^i)$'s, for any $|\eps|,|\eta|<\eps_0$\footnote{Immersions $\gamma^{i,\eps,\eta}$ define a triple junctions network by \cref{lem:SufficientConditionsNetwork} for $\eps_{\Gamma_*}$ small enough.}.

Then
\begin{equation}\label{eq:SecondVariationBilinear}
    \begin{split}
        \frac{\d}{\d\eps} \frac{\d}{\d\eta} {\rm L} (\Gamma^{\eps,\eta})\bigg|_{0,0}
        &= \sum_i \int_0^1 \partial_s X^i \partial_s Z^i |\partial_x \gamma^i_*| \de x
        \\
        &=  \sum_{p \in P_G} \partial_s X^{i_p}(1) Z^{i_p}(1) + 
        \sum_{m \in J_G} \sum_{\ell \in I_m} (-1)^{1+e^\ell_m} \partial_s X^\ell (e^\ell_m) Z^\ell (e^\ell_m) + \\
        &\phantom{=}- \sum_i \int_0^1 \partial^2_s X^i \, Z^i \, |\partial_x \gamma^i_*| \de x,
    \end{split}
\end{equation}
where $\partial_s X^i= |\partial_x\gamma^i_*|^{-1}\partial_x X^i$ and $\partial_s Z^i= |\partial_x\gamma^i_*|^{-1}\partial_x Z^i$ for any $i$.
\end{prop}

\begin{proof}
    By \cref{def:Adapted} and \cref{def:LinearOperators}, for any $i$ we have that
    \[
    \TT^{i,\eps,\eta} = \eps \TT^{i}_X + \eta \TT^{i}_Z,
    \]
    where the $\TT^{i}_X$'s are adapted to the $X^i$'s, and the $\TT^{i}_Z$'s are adapted to the $Z^i$'s. Denoting $\gamma^{i,\eps}\eqdef \gamma^{i,\eps,0}$, we compute
    \[
    \begin{split}
        \frac{\d}{\d\eps} \frac{\d}{\d\eta} &{\rm L} (\Gamma^{\eps,\eta})\bigg|_{0,0}
        = \sum_i \frac{\d}{\d\eps} \int_0^1 \frac{\scal{\partial_x \gamma^{i,\eps} , \partial_x (Z^i \nu^i_* + \TT^i_Z \tau^i_*) }}{|\partial_x \gamma^{i,\eps}|} \de x \bigg|_0 \\
        &= \sum_i \int_0^1 \left\langle 
        - \frac{\scal{\partial_x \gamma^i_*, \partial_x(X^i \nu^i_* + \TT^i_X \tau^i_*)}}{|\partial_x \gamma^i_*|^3} \partial_x \gamma^i_* + \frac{\partial_x(X^i \nu^i_* + \TT^i_X \tau^i_*)}{|\partial_x \gamma^i_*|} 
        , \partial_x (Z^i \nu^i_* + \TT^i_Z \tau^i_*) 
        \right\rangle \de x.
    \end{split}
    \]
    Since $\partial_x \tau^i_*=\partial_x \nu^i_*=0$ as $\Gamma_*$ is minimal, we get
    \[
    \begin{split}
         \frac{\d}{\d\eps} \frac{\d}{\d\eta} {\rm L} (\Gamma^{\eps,\eta})\bigg|_{0,0}
        &= \sum_i \int_0^1 \bigg\langle
        -\scal{\tau^i_*, \partial_s X^i \nu^i_* + \partial_s \TT^i_X \tau^i_*} \tau^i_* + \partial_s X^i \nu^i_* + \partial_s \TT^i_X \tau^i_*
        ,
        \partial_s Z^i \nu^i_* + \partial_s \TT^i_Z \tau^i_* \bigg\rangle |\partial_x \gamma^i_*|\de x \\
        &=\sum_i \int_0^1 \partial_s X^i \partial_s Z^i |\partial_x \gamma^i_*| \de x.
    \end{split}
    \]
    Integrating by parts, the claim follows.    
\end{proof}

\subsection{\L ojasiewicz--Simon inequalities for minimal networks}\label{sec:Loja}

We need to set up a functional analytic framework for proving the desired \L ojasiewicz--Simon inequalities.

For a fixed minimal network $\Gamma_*:G\to \R^2$, we denote by $M\eqdef \sharp J_G$ and $P\eqdef \sharp P_G$, and we define the Banach spaces
\begin{equation}\label{eq:DefV}
\begin{split}
    V \eqdef 
     \bigg\{ \overline{\NN}\eqdef(\NN^1,\ldots,\NN^N) \in [H^2(0,1)]^N  \st 
    &\sum_{\ell \in I_m} (-1)^{e^\ell_m} \NN^\ell(e^\ell_m) =0 \,\, \forall\, m \in J_G , \\
    & \NN^{i_p}(1)=0 \,\, \forall\, p \in P_G
    \bigg\},
\end{split}
\end{equation}
endowed with $\|\overline{\NN}\|_V^2 \eqdef \sum_i \|\NN^i\|_{H^2}^2$, and
\begin{equation}\label{eq:DefZ}
    Z\eqdef  W_1 \times \ldots \times W_M \times [L^2(0,1)]^N,
\end{equation}
endowed with the product norm, where
\begin{equation*}
    W_m \eqdef \left\{(v^\ell_m)_{\ell \in I_m} \in \R^3 \st \sum_{\ell \in I_m}(-1)^{e^\ell_m} v^\ell_m=0 \right\},
\end{equation*}
and $W_m$ is endowed with the Euclidean scalar product.

Observe that ${\rm j}:V\hookrightarrow Z$ compactly with the natural injection
\begin{equation}\label{eq:Injection}
\overline{\NN} \quad\overset{\rm j}{\mapsto}\quad \left( (\NN^\ell(e^\ell_m)) ,  \overline{\NN}\right)
\end{equation}

For $r_{\Gamma_*}>0$ small enough, we also define the energy $\boldsymbol{\rm L}:B_{r_{\Gamma_*}} (0) \subset V\to [0,+\infty)$ by
\begin{equation}\label{eq:DefEnergiaLineare}
    \boldsymbol{\rm L}(\overline{\NN}) \eqdef \sum_i {\rm L} \left(\gamma^i_* + \NN^i \nu^i_* + \TT^i \tau^i_* \right),
\end{equation}
where the $\TT^i$'s are adapted to the $\NN^i$'s (see \cref{def:Adapted}). We observe that, according to \cref{lem:SufficientConditionsNetwork}, the immersions $\gamma^i_* + \NN^i \nu^i_* + \TT^i \tau^i_* $ define a triple junctions network.

\begin{cor}
Let $\Gamma_*:G\to \R^2$ be a minimal network. Let $V,Z,\boldsymbol{\rm L}$ as above, and identify $Z^\star$ with ${\rm j}^\star(Z^\star)\subset V^\star$, for ${\rm j}$ as in \eqref{eq:Injection}.

Then the following hold.
\begin{enumerate}
    \item The first variation $\delta \boldsymbol{\rm L} : V\to Z^\star$ is $Z^\star$-valued by setting
    \begin{equation}\label{eq:FirstVariationBold}
    \begin{split}
        \delta \boldsymbol{\rm L} & (\overline{\NN})[((v^\ell_m), \overline{X})]
    \\&= 
   % \sum_{p \in P_G} \scal{\tau^{i_p}(1), \nu^{i_p}_*(1)} u^{i_p} +\\& \phantom{=}   +
    \sum_{m \in J_G}  \sum_{\ell \in I_m} 
    (-1)^{1+e^\ell_m} \left[ 
    \scal{\tau^\ell(e^\ell_m), \nu^\ell_*(e^\ell_m)} + \sum_{j \in I_m} h_{\ell j}\scal{\tau^j(e^j_m) , \tau^j_*(e^j_m)}
     \right] 
     v^\ell_m + \\
    &\phantom{=} - \sum_i 
    \int_0^1 \bigg( \scal{\boldsymbol{k}^i,\nu^i_*}|\partial_x \gamma^i| 
    + \sum_j f_{ij}\chi \scal{\boldsymbol{k}^j,\tau^j_*}|\partial_x \gamma^j| +\\
    &\qquad\qquad\qquad
    + g_{ij}\chi(1-x) \scal{\boldsymbol{k}^j,\tau^j_*}(1-x)|\partial_x \gamma^j|(1-x) \bigg) X^i \de x,
\end{split}
    \end{equation}
where $f_{ij}, g_{ij}, h_{\ell j} \in \R$ depend on the topology of $G$, and $\tau^i, \boldsymbol{k}^i$ are referred to the immersions $\gamma^i \eqdef \gamma^i_*  + \NN^i \nu^i_* + \TT^i\tau^i_*$, with $\TT^i$ adapted to $\NN^i$.

If also, the network defined by the immersions $\gamma^i$ is regular, then
\begin{equation}\label{eq:FirstVariationBoldRegular}
    \begin{split}
        \delta \boldsymbol{\rm L}  (\overline{\NN})[((v^\ell_m), \overline{X})]
    &= 
   - \sum_i 
    \int_0^1 \bigg( \scal{\boldsymbol{k}^i,\nu^i_*}|\partial_x \gamma^i|
    + \sum_j f_{ij}\chi \scal{\boldsymbol{k}^j,\tau^j_*}|\partial_x \gamma^j| +\\
    &\qquad\qquad\qquad
    + g_{ij}\chi(1-x) \scal{\boldsymbol{k}^j,\tau^j_*}(1-x)|\partial_x \gamma^j|(1-x) \bigg) X^i \de x,
    \end{split}
\end{equation}

    \item The second variation $\delta^2 \boldsymbol{\rm L}_0 : V\to Z^\star $ at $0$ is $Z^\star$-valued by setting
    \begin{equation}\label{eq:SecondVariationBold}
    \begin{split}
        \delta^2 \boldsymbol{\rm L}_0 ( \overline{X} ) [((v^\ell_m), \overline{Z}) ]
        &=
        %\sum_{p \in P_G} \partial_s X^{i_p} (1) u^{i_p} +
        \sum_{m \in J_G} \sum_{\ell \in I_m} (-1)^{1+e^\ell_m} \partial_s X^\ell (e^\ell_m) v^\ell_m + \\
        &\phantom{=}- \sum_i \int_0^1 \bigg(|\partial_x \gamma^i_*| \partial^2_s X^i 
        %+ \sum_j |\partial_x \gamma^i_*|\Omega_{ij}(\overline{X})
        \bigg) Z^i \de x,
    \end{split}
    \end{equation}
where $\partial_s X^n= |\partial_x\gamma^n_*|^{-1}\partial_x X^n$ for any $n$.
\end{enumerate}
\end{cor}

\begin{proof}
For the sake of precision, we maintain $Z^\star$ and ${\rm j}^\star(Z^\star)$ distinct in this proof.

The first item follows by \cref{prop:FirstVariation}. Let $\overline{\NN},\overline{X}\in V$. Equation \eqref{eq:FirstVariation} yields the expression for $\delta \boldsymbol{\rm L}(\overline{\NN}) \in V^\star $, and we notice that, since $\overline{X}\in V$, the sum over endpoints $p \in P_G$ in \eqref{eq:FirstVariation} vanishes. Hence \eqref{eq:FirstVariation} shows that there exists an element $\nabla \boldsymbol{\rm L}(\overline{\NN})$ of $Z$ such that $\delta \boldsymbol{\rm L}(\overline{\NN})[\overline{X}] = \scal{\nabla \boldsymbol{\rm L}(\overline{\NN}), {\rm j}(\overline{X})}_{Z}$. Letting ${\rm I}:Z\to Z^\star$ the natural isometry, this means
\[
\delta \boldsymbol{\rm L}(\overline{\NN})[\overline{X}] = \scal{\nabla \boldsymbol{\rm L}(\overline{\NN}), {\rm j}(\overline{X})}_{Z} 
= \scal{I\left(\nabla \boldsymbol{\rm L}(\overline{\NN}) \right),{\rm j}(\overline{X})}_{Z^\star,Z}
= \scal{\,{\rm j}^\star \left(I\left(\nabla \boldsymbol{\rm L}(\overline{\NN}) \right)\right),\overline{X}}_{V^\star,V},
\]
that is, $\delta \boldsymbol{\rm L}(\overline{\NN}) \in {\rm j}^\star(Z^\star)$, and \eqref{eq:FirstVariationBold} follows as well. By the same reasoning, \eqref{eq:FirstVariationBoldRegular} follows from \eqref{eq:FirstVariationRegular}.

The second item analogously follows from \cref{prop:SecondVariationBilinear}. In this case we notice that the sum over endpoint $p \in P_G$ in \eqref{eq:SecondVariationBilinear} vanishes whenever $\overline{Z} \in V$, leading to \eqref{eq:SecondVariationBold}.
\end{proof}

Now we start checking that the assumptions needed to imply a \L ojasiewicz--Simon inequality hold, see \cref{prop:LojaAbstract}. We start from the analyticity of the functional and of its first variation.

\begin{lemma}\label{lem:Analyticity}
Let $\Gamma_*:G\to \R^2$ be a minimal network. Let $V,Z,\boldsymbol{\rm L}, r_{\Gamma_*}$ as above, and identify $Z^\star$ with ${\rm j}^\star(Z^\star)\subset V^\star$, for ${\rm j}$ as in \eqref{eq:Injection}.

Then the maps $\boldsymbol{\rm L}:B_{r_{\Gamma_*}} (0) \subset V\to [0,+\infty)$ and $\delta \boldsymbol{\rm L} : V\to Z^\star$ are analytic.
\end{lemma}

\begin{proof}
The claim easily follows by recalling that multilinear continuous maps are analytic and that sum and compositions of analytic maps are analytic. Moreover if $T_j:U\subset B\to C_j$, for $j=1,2$, is analytic from an open set $U$ of a Banach space $B$ into a Banach space $C_j$, and $\cdot:C_1\times C_2\to D$ is a bilinear continuous map into a Banach space $D$, then the ``product operator'' $T(v,w)\eqdef T_1(v)\cdot T_2(w)$ is analytic from $U$ into $D$.

Concerning analyticity of $\boldsymbol{\rm L}$ we need to check that
\begin{equation*}
    B_{r_{\Gamma_*}} (0) \ni \,\NN \quad\mapsto\quad \int_0^1 \left|\partial_x \left(\gamma^i_* + \NN^i \nu^i_* + \TT^i \tau^i_* \right) \right| \de x,
\end{equation*}
is analytic for any $i$. Since the $\TT^i$'s are adapted, they depend linearly on the $\NN^i$'s, moreover differentiation with respect to $x$ is linear and continuous from $V$ to $[H^1(0,N)]^N$. Also, for $r_{\Gamma_*}$ sufficiently small, we have that $ \left|\partial_x \left(\gamma^i_* + \NN^i \nu^i_* + \TT^i \tau^i_* \right) \right|\ge c_*>0$, for $c_*$ depending on $\Gamma_*, r_{\Gamma_*}$ only. Finally integration is linear and continuous on $L^1(0,1)$. Putting together all these observations, we get that $\boldsymbol{L}$ is analytic.

The analyticity of $\delta \boldsymbol{\rm L} : V\to Z^\star$ follows by completely analogous observations, recalling the expression in \eqref{eq:FirstVariationBold}. Indeed, one can check that tangent and curvature vectors to an immersion $\gamma^i_* + \NN^i \nu^i_* + \TT^i \tau^i_*$ depend analytically on the parametrization, and then on $\NN$ (see for example the analogous treatment in \cite[Section 3.1, Appendix B]{DaPoSp16}). Moreover, the trace operator evaluating a tangent vector $\tau^\ell \in H^1(0,1)$ at junction points is linear and continuous. Recalling that product operators of analytic maps are analytic, the analyticity of $\delta \boldsymbol{\rm L}$ follows.
\end{proof}

Now we need to prove that the second variation is Fredholm of index zero. We recall that a continuous linear operator $T$ between Banach spaces is Fredholm of index zero if its kernel has finite dimension, its image has finite codimension, and such dimensions are equal.

\begin{lemma}\label{lem:Fredholmness}
Let $\Gamma_*:G\to \R^2$ be a minimal network. Let $V,Z,\boldsymbol{\rm L}$ as above, and identify $Z^\star$ with ${\rm j}^\star(Z^\star)\subset V^\star$, for ${\rm j}$ as in \eqref{eq:Injection}.

Then the second variation $\delta^2 \boldsymbol{\rm L}_0 : V\to Z^\star $ at $0$ is a Fredholm operator of index zero.
\end{lemma}

\begin{proof}
Denote by  ${\rm I}:Z\to Z^\star$ the natural isometry.
Recalling \eqref{eq:SecondVariationBold}, we see that the claim follows as long as we can prove that the following operator is Fredholm of index $0$:
\begin{equation*}
    V \ni \quad \overline{X} \quad\mapsto\quad {\rm I}\left(
    %(\partial_s X^{i_p}(1) ), 
    ((-1)^{1+e^{\ell}_m}\partial_s X^\ell(e^\ell_m) ),
    -|\partial_x\gamma^i_*| \partial^2_s X^i 
    \right) \quad\in Z^\star.
\end{equation*}
Let
\begin{equation*}
\begin{split}
    V_1 \eqdef \bigg\{ \overline{X}\eqdef(X^1,\ldots,X^N) \in [H^1(0,1)]^N \st 
    &\sum_{\ell \in I_m} (-1)^{e^\ell_m} X^\ell(e^\ell_m) =0 \,\, \forall\, m \in J_G  , \\
    & X^{i_p}(1)=0 \,\, \forall\, p \in P_G \bigg\},
\end{split}
\end{equation*}
and let $( (V^\ell_m), \overline{Z}) \in Z$ be fixed. We consider the operator $F: V_1 \to \R$ given by
\begin{equation*}
    F(\overline{Y}) \eqdef 
    %\sum_{p \in P_G} U^{i_p}Y^{i_p}(1) + 
    \sum_{m \in J_G} \sum_{\ell \in I_m} V^\ell_m Y^\ell(e^\ell_m) + \sum_i \int_0^1 |\partial_x \gamma^i_* |Z^i Y^i \de x.
\end{equation*}
We can endow $V_1$ with the scalar product $\scal{\overline{X},\overline{Y}} \eqdef \sum_i \int_0^1 \partial_s X^i \partial_s Y^i + X^iY^i \de s$, where $\d s=\d s_{\gamma^i_*}$ along the $i$th edge. Hence $F:(V_1,\scal{\cdot,\cdot})\to \R$ is linear and continuous, and then there exists a unique $\overline{X} \in V_1$ such that
\begin{equation}\label{eq:zz8}
    \sum_i \int_0^1 \partial_s X^i \partial_s Y^i + X^iY^i \de s =  
    %\sum_{p \in P_G} U^{i_p}Y^{i_p}(1) +
    \sum_{m \in J_G} \sum_{\ell \in I_m} V^\ell_m Y^\ell(e^\ell_m) + \sum_i \int_0^1 Z^i Y^i \de s,
\end{equation}
for any $\overline{Y} \in V_1$. Testing on $\overline{Y} \in V_1$ such that $Y^i\equiv 0$ for all $i$ except for a fixed index $j$, and $Y^j \in C^1_c(0,1)$, we see that
\[
\int_0^1 \partial_s X^j\partial_s Y^j + X^jY^j \de s =   \int_0^1 Z^j Y^j \de s,
\]
which implies that $X^j \in H^2(0,1)$ with $-\partial^2_s X^j + X^j = Z^j$, and thus $\overline{X}$ belongs to $V$.

For $m\in J_G$, we can now take $\overline{Y} \in V_1$ with $Y^\ell\equiv 0$ for all $\ell$ except for $\ell \in I_m$, with $Y^\ell \in C^1$ vanishing at the endpoint of $E^\ell$ different from the junction $m$. Integration by parts in \eqref{eq:zz8} then gives
\begin{equation*}
    \sum_{\ell \in I_m} (-1)^{1+ e^\ell_m} \partial_s X^\ell(e^\ell_m) Y^\ell(e^\ell_m) = \sum_{\ell \in I_m} V^\ell_m Y^\ell(e^\ell_m).
\end{equation*}
Arbitrariness of $\overline{Y}$ implies that $\sum_{\ell \in I_m} \left((-1)^{1+ e^\ell_m} \partial_s X^\ell(e^\ell_m) - V^\ell_m \right) v^\ell =0$ for any triple $\{ v^\ell \st \ell \in I_m, \, \sum_{\ell \in I_m} (-1)^{e^\ell_m} v^\ell =0 \}$. This means that there exists a constant $\alpha_m \in \R$ such that
\begin{equation*}
    (-1)^{1+ e^\ell_m} \partial_s X^\ell(e^\ell_m) - V^\ell_m  = \alpha_m (-1)^{e^\ell_m}
    \qquad
    \forall \ell \in I_m.
\end{equation*}
Multiplying by $(-1)^{e^\ell_m}$ and summing over $\ell$ implies that $3 \alpha_m = - \sum_{\ell \in I_m} \partial_s X^\ell(e^\ell_m)$, and then
\begin{equation*}
    (-1)^{1+e^{\ell}_m}\left(\partial_s X^\ell(e^\ell_m) - \frac13 \sum_{\ell \in I_m} \partial_s X^\ell(e^\ell_m) \right) = V^\ell_m
    \qquad\forall\, \ell \in I_m.
\end{equation*}
%Arguing analogously at each endpoint $p \in P_G$, eventually 
Therefore, we have proved that for arbitrary $( (V^\ell_m), \overline{Z}) \in Z$ there exists a unique $\overline{X} \in V$ satisfying
\begin{equation*}
    \begin{cases}
    -\partial^2_s X^i + X^i = Z^i & \forall \, i , \\
    %\partial_s X^{i_p}(1) = U^{i_p} & \forall\, p \in P_G, \\
    (-1)^{1+e^{\ell}_m}\left(\partial_s X^\ell(e^\ell_m) - \frac13 \sum_{\ell \in I_m} \partial_s X^\ell(e^\ell_m) \right)  = V^\ell_m & \forall\, m \in J_G, \, \ell \in I_m.
    \end{cases}
\end{equation*}
Therefore, if we further define the linear and continuous operator $\mathscr{F}: V \to Z^\star$ given by
\begin{equation*}
    \mathscr{F}(\overline{X}) \eqdef {\rm I} \left( 
    %(\partial_s X^{i_p}(1) ), 
    \left((-1)^{1+e^{\ell}_m}\left(\partial_s X^\ell(e^\ell_m) - \frac13 \sum_{\ell \in I_m} \partial_s X^\ell(e^\ell_m) \right)\right),
    -|\partial_x\gamma^i_*| \partial^2_s X^i + |\partial_x\gamma^i_*| X^i\right),
\end{equation*}
where ${\rm I}:Z\to Z^\star$ is the natural isometry, we see that $\mathscr{F}$ is invertible, and thus it is Fredholm of index $0$.\\
Recall that Fredholmness is stable under compact perturbations: a linear operator $T$ between Banach spaces is Fredholm of index $l$ if and only if $T+K$ is Fredholm of index $l$, for any compact operator $K$ (see \cite[Section 19.1]{HormanderIII}). Therefore, since
\begin{equation*}
    V \ni \quad \overline{X} \quad\mapsto\quad
    {\rm I}\left( 
    %(0), 
    \left( -(-1)^{1+e^{\ell}_m}\frac13 \sum_{\ell \in I_m} \partial_s X^\ell(e^\ell_m)  \right) ,|\partial_x\gamma^i_*| X^i \right) \quad \in Z^\star,
\end{equation*}
is compact, we conclude that
\begin{equation*}
    V \ni \quad \overline{X} \quad\mapsto\quad {\rm I}\left(
    %(\partial_s X^{i_p}(1) ), 
    ((-1)^{1+e^{\ell}_m}\partial_s X^\ell(e^\ell_m) ),
    -|\partial_x\gamma^i_*| \partial^2_s X^i 
    \right) \quad\in Z^\star,
\end{equation*}
is Fredholm of index $0$ as well, completing the proof.
\end{proof}

We can now apply the following abstract result stating sufficient conditions implying a \L ojasiewicz--Simon gradient inequality.

\begin{prop}[{\cite[Corollary~2.6]{PozzettaLoja}}]\label{prop:LojaAbstract}
Let $E:B_{\rho_0}(0)\subseteq V \to \R$ be an analytic map, where $V$ is a Banach space. Suppose that $0$ is a critical point for $E$, i.e., $\delta E_0 = 0$. Assume that there exists a Banach space $Z$ such that $V\hookrightarrow Z$, the first variation $\delta E : B_{\rho_0}(0)\to Z^\star$ is $Z^\star$-valued and analytic and the second variation $\delta^2 E_0 : V \to Z^\star$ evaluated at $0$ is $Z^\star$-valued and Fredholm of index zero.\\
Then there exist constants $C,\rho_1>0$ and $\theta \in (0,1/2]$ such that
\[
|E(v)- E(0)|^{1-\theta} \le C \| \delta E_v \|_{Z^\star},
\]
for every $v \in B_{\rho_1}(0) \subseteq V$.
\end{prop}

The above functional analytic result is a corollary of the useful theory developed in~\cite{Ch03} and it has been independently observed in~\cite{Rupp20}.

\begin{teo}[\L ojasiewicz--Simon inequality at minimal networks]\label{thm:LojaGenerale}
Let $\Gamma_*:G\to\R^2$ be a minimal network. Let $V,Z$ be as in \eqref{eq:DefV}, \eqref{eq:DefZ}, and define $\boldsymbol{\rm L}:B_{r_{\Gamma_*}} (0) \subset V\to [0,+\infty)$ as in \eqref{eq:DefEnergiaLineare}.

Then there exist $C_{\rm LS}>0$, $\theta \in (0,\tfrac12]$, and $r \in (0,r_{\Gamma_*}]$ such that
\begin{equation*}
    \left|\boldsymbol{\rm L}(\overline{\NN}) - {\rm L}(\Gamma_*)  \right|^{1-\theta} \le C_{\rm LS} \left\|\delta \boldsymbol{\rm L}(\overline{\NN}) \right\|_{Z^\star},
\end{equation*}
for any $\overline{\NN} \in B_r(0)\subset V$.
\end{teo}

\begin{proof}
The proof immediately follows by applying \cref{prop:LojaAbstract} recalling \cref{lem:Analyticity} and \cref{lem:Fredholmness}.
\end{proof}

We can finally derive the following more explicit \L ojasiewicz--Simon inequality for regular networks.

\begin{cor}\label{cor:LojaNetworkMinimali}
Let $\Gamma_*:G\to\R^2$ be a minimal network. Then there exist $C_{\rm LS},\sigma >0$ and $\theta \in (0,\tfrac12]$ such that the following holds.

If $\Gamma:G\to\R^2$ is a regular network of class $H^2$ such that
\begin{align}
     \sum_i \| \gamma^i_* - \gamma^i \|_{H^2(\d x)} \le \sigma,& \\
    \gamma^i_*(p)=\gamma^i(p) &\qquad \text{$\forall\,p\in G\st p$ is an endpoint,} \label{eq:wzw}
 \end{align}
 then
 \begin{equation}\label{eq:LojaNetworkMinimali}
     \left|{\rm L}(\Gamma)-{\rm L}(\Gamma_*) \right|^{1-\theta}
     \le C_{\rm LS} \left( \sum_i \int_0^1 |\boldsymbol{k}^i|^2 \de s \right)^{\frac12}.
\end{equation}
\end{cor}

\begin{proof}
For $\sigma$ small enough, applying \cref{prop:ParametrizzazioneTN} and recalling \eqref{eq:StimaH2Ni}, we know that there exist functions $\NN^i, \TT^i \in H^2(\d x)$, where $\TT^i$'s are adapted to the $\NN^i$'s, and reparametrizations $\varphi^i:[0,1]\to[0,1]$ such that
\begin{equation*}
     \gamma^i\circ \varphi^i(x) = \gamma^i_*(x) + \NN^i(x)\nu^i_*(x) + \TT^i(x)\tau^i_*(x) \eqqcolon \widetilde{\gamma}^i.
\end{equation*}
Moreover, by \eqref{eq:wzw}, \cref{lem:BoundaryRelationsTN}, and up to decreasing $\sigma$, we have that $\overline{\NN}\eqdef(\NN^1,\ldots,\NN^N)$ belongs to the ball $B_r(0)\subset V$, where $r,V$ are as in \cref{thm:LojaGenerale}.

For $\boldsymbol{\rm L}, Z$ as in \cref{thm:LojaGenerale}, since $\Gamma$ is regular, by \eqref{eq:FirstVariationBoldRegular} we get that
\begin{equation*}\label{eq:zzxx}
    \begin{split}
        \|\delta \boldsymbol{\rm L}(\overline{\NN})\|_{Z^\star}^2
        &=\sum_i \int_0^1 \bigg| \scal{\widetilde{\boldsymbol{k}}^i,\nu^i_*}|\partial_x \widetilde\gamma^i| 
    + \sum_j f_{ij}\chi \scal{\widetilde{\boldsymbol{k}}^j,\tau^j_*}|\partial_x \widetilde\gamma^j| +\\
    &\qquad+
     g_{ij}\chi(1-x) \scal{\widetilde{\boldsymbol{k}}^j,\tau^j_*}(1-x)|\partial_x \widetilde\gamma^j|(1-x) \bigg) \bigg|^2 \de x \\
     &\le C(\Gamma_*,\sigma)  \sum_i \int_0^1 |\widetilde{\boldsymbol{k}}^i|^2 \de s .
    \end{split}
\end{equation*}
Since ${\rm L}(\Gamma)=\boldsymbol{\rm L}(\overline{\NN})$ and the $L^2(\d s)$ norm of the curvature on the right hand side of \eqref{eq:zzxx} does not depend on the parametrization, the above estimate together with \cref{thm:LojaGenerale} imply \eqref{eq:LojaNetworkMinimali}.
\end{proof}

\begin{rem}[Further \L ojasiewicz--Simon inequalities at minimal networks]\label{rem:Loja2}
By an adaptation of the above arguments, we expect to be possible to prove a \L ojasiewicz--Simon inequality at minimal networks taking into account also variations at endpoints.

More precisely, removing the constraint $\NN^{i_p}(1)=0$ for $\overline{\NN} \in V$ in \eqref{eq:DefV}, and considering $\widetilde{Z}\eqdef \R^P \times Z$, for $Z$ as in \eqref{eq:DefZ}, employing the variation formulae in \cref{prop:FirstVariation} and \cref{prop:SecondVariationBilinear}, one can consider triple-junctions networks $\Gamma$ in a neighborhood of a minimal one $\Gamma_*$ having endpoints different from those of $\Gamma_*$.

Arguing as in the above propositions, one eventually deduces an analog of \cref{thm:LojaGenerale}. The resulting statement would formally read exactly as \cref{thm:LojaGenerale}, but in this case the norm $\left\|\delta \boldsymbol{\rm L}(\overline{\NN}) \right\|_{Z^\star}$ on the right hand side of the inequality also counts contributions from the varied endpoints. More precisely, all the terms in the first variation formula \eqref{eq:FirstVariationRegular} representing the operator $\delta \boldsymbol{\rm L}(\overline{\NN})$ do not vanish in general and thus contribute to its norm.
\end{rem}

\section{Minimal networks locally minimize length}\label{sec:MinimalMinimizing}

% \begin{defn}
% We say that a regular network $\Gamma_*:G\to \R^2$ \emph{locally minimizes the length in $H^k$} for $k\ge2$ (resp. \emph{in $C^{k,\alpha}$} for $k\ge1$, $\alpha \in [0,1]$) if there exists $\eta>0$ such that ${\rm L}(\Gamma) \ge {\rm L}(\Gamma_*)$ whenever $\Gamma:G\to \R^2$ is a regular network having the same endpoints of $\Gamma_*$ and such that $\|\gamma^i\circ \sigma^i-\gamma^i_*\|_{H^k} < \eta$ (resp. $\|\gamma^i\circ \sigma^i-\gamma^i_*\|_{C^{k,\alpha}} < \eta$), for some reparametrizations $\sigma^i$.
% \end{defn}

In this section we provide a simple proof of the fact that minimal networks are automatically local minimizers for the length with respect to perturbations sufficiently small in $C^0$.

More precisely, we say that a regular network $\Gamma_*:G\to \R^2$ \emph{locally minimizes the length in $C^0$} if there exists $\eta>0$ such that ${\rm L}(\Gamma) \ge {\rm L}(\Gamma_*)$ whenever $\Gamma:G\to \R^2$ is a regular network having the same endpoints of $\Gamma_*$ and such that $\|\gamma^i\circ \sigma^i-\gamma^i_*\|_{C^0} < \eta$, for some reparametrizations $\sigma^i$.

We mention that more general minimality properties of minimal networks can be proved, see \cite{Mor, MartelliNovagaPludaRiolo, Wh3, PludaPozzCalibrations, FischerHenLauxSimonCalibrazioni}.

\begin{lemma}\label{lemma:MinimalMinimizing}
Let $\Gamma_*:G\to \R^2$ be a minimal network. Then $\Gamma_*$ locally minimizes the length in $C^0$.
\end{lemma}

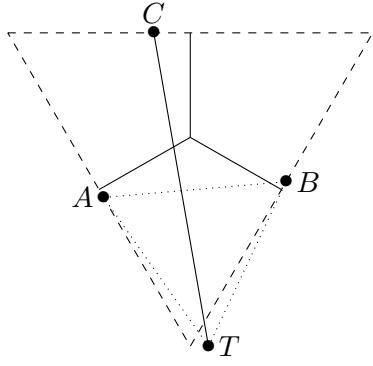
\begin{figure}[H]
	\begin{center}
		\begin{tikzpicture}[scale=1.2]
		\draw[dashed]
		(-2,1.732)to[out=0, in=180, looseness=1]
		(2,1.732)to[out=240, in=60, looseness=1]
		(0,-1.732)to[out=120, in=300, looseness=1](-2,1.732);
        \draw
        (-1,0)--(0,0.577)--(0,1.732);
        \draw
        (1,0)--(0,0.577);
        \path[font=\normalsize]
		(-0.95,-0.087)node[]{$\bullet$}; %s=-0.1
        \path[font=\normalsize]
		(1.05,0.087)node[]{$\bullet$}; %t=0.1
		\path[font=\normalsize]
		(-0.4,1.732)node[]{$\bullet$}; %x=-0.4
		\draw[dotted]
		(-0.95,-0.087)--(1.05,0.087)--(0.2,-1.732)--(-0.95,-0.087);
		\path[font=\normalsize]
		(0.2,-1.732)node[]{$\bullet$};
		\draw
		(0.2,-1.732)--(-0.4,1.732);
		\path[font=\normalsize]
		(-0.95,-0.087)node[left]{$A$};
		\path[font=\normalsize]
		(1.05,0.087)node[right]{$B$};
		\path[font=\normalsize]
		(-0.4,1.732)node[above]{$C$};
		\path[font=\normalsize]
		(0.2,-1.732)node[right]{$T$};
		\end{tikzpicture}
	\end{center}
	\caption{Standard triod joining mid points on the sides of an equilateral triangle (dashed lines). Dotted lines: equilateral triangle constructed over the side $AB$. For other endpoints $A, B, C$ close to such midpoints, the length of the Steiner tree joining $A, B, C$ is equal to the length of $CT$.}\label{Fig:MinimalMinimizing}
\end{figure}

\begin{proof}
For any $r>0$ and for any junction $m=\pi(e^i,i)=\pi(e^j,j)=\pi(e^k,k)$ of $G$, let $T_{r,m}$ be the closed equilateral triangle having $\Gamma_*(m)$ as barycenter and whose sides have length $r$ and are orthogonal to the inner tangent vectors at $m$, that is, the vectors $(-1)^{e^i}\tau^i(e^i)$, $(-1)^{e^j}\tau^j(e^j)$, $(-1)^{e^k}\tau^k(e^k)$.

Now fix $r>0$ small enough such that the set $T_{r,m} \cap \Gamma_*(G)$ is a standard triod for any junction $m$, i.e., such set is given by the union of three straight segments of the same length having one end in common forming angles equal to $\tfrac23\pi$ (see \cref{Fig:MinimalMinimizing}).

Let $\Gamma:G\to \R^2$ be a smooth regular network with the same endpoints of $\Gamma_*$. If, up to reparametrization, the immersions defining $\Gamma$ are close to the ones of $\Gamma_*$ in $C^0$, then for any edge $E_i$ if, say, $m=\pi(0,i)$ and $m'=\pi(1,i)$ are two junctions, we can fix times $0<t_{i,1}<t_{i,2}<1$ such that $t_{i,1}$ is the last time $\gamma^i$ intersects $\partial T_{r,m}$ and $t_{i,2}$ is the first time $\gamma^i$ intersects $\partial T_{r,m'}$. Such intersections define points close to $(\partial T_{r,m} ) \cap \Gamma_*(G)$ and $(\partial T_{r,m'} ) \cap \Gamma_*(G)$. In case $\pi(0,i)$ is an endpoint, we set $t_{i,1}=0$.

%$T_{r,m} \cap \Gamma(G)$ has the topology of a standard triod and intersects the boundary of $T_{r,m}$ in three points close to $(\partial T_{r,m} ) \cap \Gamma_*(G)$.

In order to complete the proof, if, say, $m=\pi(0,i)=\pi(0,j)=\pi(0,k)$ is a junction, it is sufficient to prove that the length of $\Gamma_*$ in $T_{r,m}$ is smaller than the sum $\sum_{\ell=i,j,k} {\rm L}(\gamma^\ell|_{(0,t_{\ell,1})})$. Indeed, $\Gamma_*(G)\setminus \cup_m T_{r,m}$ is given by straight segments orthogonal to the sides of the triangles $T_{r,m}$ and whose endpoints lay either on parallel sides of different triangles $T_{r,m'}, T_{r,m''}$, or on a side of a triangle $T_{r,m}$ and on an endpoint $\Gamma_*(p)$ of the network. Hence the length of $\Gamma_*$ outside $ \cup_m T_{r,m}$ is automatically smaller than the sum of the lengths of the curves of $\Gamma$ on intervals $(t_{1,i},t_{2,i})$.

Eventually, the argument reduces to prove that the length of a standard triod $\mathbb T$ whose endpoints are the mid points of the sides of an equilateral triangle is the least possible among the length of topological triods having endpoints on the sides of the same triangle close to the ones of $\mathbb T$ (see \cref{Fig:MinimalMinimizing}). Up to scaling and translation, let us assume that the endpoints of a standard triod are located at points $(-1,0), (1,0), (0,\sqrt{3})$ in the plane. Hence the endpoints of a competitor triod take the form $A=(-1,0)+s(-\tfrac12,\tfrac{\sqrt{3}}{2})$, $B=(1,0)+t(\tfrac12,\tfrac{\sqrt{3}}{2})$, $C=(x,\sqrt{3})$ for $s,t,x$ close to zero (see \cref{Fig:MinimalMinimizing}). The length of the competitor triod in greater or equal than the one of the Steiner tree joining $A, B, C$, which is another topological triod $\mathbb S$ whose total length can be shown to be equal to the length of the segment $CT$, where $T$ is the point (farthest from $C$) such that points $A,B,T$ are vertices of an equilateral triangle (see \cref{Fig:MinimalMinimizing} and \cite{PaoliniSteiner}). In the end, the proof follows if we prove that ${\rm L}(\mathbb T) \le {\rm L}(CT)$.

In our choice of coordinates we have that ${\rm L}(\mathbb T)=2\sqrt{3}$. On the other hand we have that $T=A + {\rm R}(B-A)$, where ${\rm R}$ is the clockwise rotation of an angle equal to $\tfrac{\pi}{3}$. Hence
\[
\begin{split}
    T=A+ \frac12\begin{pmatrix}
1 & \sqrt{3} \\
-\sqrt{3} & 1
\end{pmatrix}(B-A) = (t-s, -\sqrt{3}).
\end{split}
\]
Then ${\rm L}(CT)^2 = (t-s-x)^2+(-\sqrt{3}-\sqrt{3})^2 \ge (2\sqrt{3})^2 = {\rm L}(\mathbb T)^2$, which completes the proof.
\end{proof}

\section{Stability and convergence}\label{sec:StabilityConvergence}

In this section we prove our main stability theorem. First we need the next technical lemma, which is based on a simple contradiction argument implying that the motion by curvature starting sufficiently close to a minimal network $\Gamma_*$ in $H^2$ passes as close as prescribed to $\Gamma_*$ in $C^k$ at some positive time.

\begin{lemma}\label{lem:Regularization2}
Let $\Gamma_*:G\to \R^2$ be a minimal network. Hence, for any $\eta>0$ and $k \in \N$ there exists $\overline{\eps}=\overline{\eps}(\Gamma_*,\eta,k)>0$ such that the following holds.

For any smooth regular network $\Gamma_0:G\to \R^2$ such that $\|\gamma^i_0-\gamma^i_*\|_{H^2} < \overline{\eps}$, the motion by curvature $\Gamma_t:G\to \R^2$, for $t \in [0,T)$, starting from $\Gamma_0$ satisfies
\begin{equation}\label{eq:Regularization2Tesi}
    \|\gamma^i_\tau\circ\sigma^i-\gamma^i_*\|_{C^k} < \eta,
\end{equation}
for some $\tau \in (0,T)$ and smooth reparametrizations $\sigma^i$, for any $i$.
\end{lemma}

\begin{proof}
Suppose by contradiction that there are $\eta>0,k \in \N$ and a sequence of smooth regular networks $\Gamma_{n,0}:G\to \R^2$ such that $\|\gamma^i_{n,0}-\gamma^i_*\|_{H^2} < 1/n$, but the motions by curvature $\Gamma_{n,t}:G\to \R^2$, defined on maximal intervals $[0,T_n)$ and starting from $\Gamma_{n,0}$, satisfy
\begin{equation}\label{eq:zzContradiction2}
    \|\gamma^i_{n,t}\circ\sigma^i_t-\gamma^i_*\|_{C^k} \ge \eta_0,
\end{equation}
for any $t \in (0,T_n)$ and any reparametrizations $\sigma^i_t$, where $\sigma^i_t$ is smooth with respect to $x$.

By \cref{wellposedness}, since $\|\gamma^i_{n,0}-\gamma^i_*\|_{H^2} \to 0$ for any $i$ as $n\to\infty$, there exists $T>0$ such that $T_n> 2 T$ for any $n$.
Moreover the solutions $\mathcal N_n$ of the motion by curvature starting from $\Gamma_{n,0}$ satisfy a uniform bound $\|\mathcal N_n \|_{W^{1,2}_5}\le M=M(\Gamma_*)$. By the compact embedding $W^{1,2}_5\hookrightarrow W^{1,2}_4$, it follows that, up to subsequence, the solutions $\gamma^i_{n,t}$ converge in $W^{1,2}_4\left((0,T)\times(0,1);\R^2\right)$ to limit immersions $\gamma^i_{\infty,t}$. Moreover, $\gamma^i_{n,0}\to \gamma^i_{\infty,0}=\gamma^i_*$ in $H^2$ and passing to the limit at almost every $t,x$ in
\[
\scal{\partial_t \gamma^i_{n,t}, \nu^i_{n,t}} \nu^i_{n,t} = \boldsymbol{k}^i_{n,t},
\]
we deduce that the maps $\gamma^i_{\infty,t}$ give a solution to the motion by curvature starting from $\Gamma_*$.
Since $\Gamma_*$ is minimal, then $\gamma^i_{\infty,t}$ actually coincides with $\gamma^i_*$ up to reparametrization.

From the uniform bound in $W^{1,2}_5$, we can fix $s\in(0,T)$ such that $\gamma^i_{n,s}\to \gamma^i_{\infty,s}$ in $H^2$ for any $i$. Hence the $L^2(\d s)$-norm of the curvature of $\gamma^i_{n,s}$ is bounded from above and the length ${\rm L}(\gamma^i_{n,s})$ is bounded from below away from zero, independently of $n$. Recalling from \cref{wellposedness} and \cref{rem:SoluzioniConfrontoCompatibilita} that for positive times the flow is smooth and it evolves according to $\partial_t\gamma^i_{n,t} = \partial^2_x\gamma^i_{n,t}/|\partial_x\gamma^i_{n,t}|^2$, we can apply the regularity estimates from \cite[Proposition 5.10, Proposition 5.8]{mannovplusch} considering $\Gamma_{n,s}$ as a new initial datum.
%
%Se avessi capito male e $\Gamma_{n,s}$ va riparametrizatto per renderlo compatibile allora si fa così: prendi ste riparametrizzazioni $\widetilde\gamma^i_{n,s}$, i boud $L^2(\d s)$ e su ${\rm L}$ non dipendono dalla riparam, quindi riparti con la teoria di \cite{mannovplusch} da $\widetilde\Gamma_{n,s}$, e similmente concludi.
%
%
This implies that there are $s<T_1\le T$ and $C_m>0$, for any $m \in \N$, independent of $n$ such that $\|\boldsymbol{k}^i_{n,t}\|_{H^m}(\d s)\le C_m$ for any $t \in [s,T_1]$.

Therefore the sequence of flows $\Gamma_{n,t}$ converges smoothly on $[s,T_1]\times G$, up to reparametrizations, to the motion by curvature $\widehat\Gamma_{\infty,t}$ parametrized by $\widehat\gamma^i_{\infty,t}$, and $\widehat\gamma^i_{\infty,t}$ is a reparametrization of $\gamma^i_*$. As the convergence holds in $H^m$ for any $m \in \N$, we find a contradiction with \eqref{eq:zzContradiction2} at any $t \in [s,T_1]$ for large $n$.
\end{proof}

\begin{teo}\label{thm:Stability}
Let $\Gamma_*:G\to \R^2$ be a minimal network. Then there exists $\delta_{\Gamma_*}>0$ such that the following holds.

Let $\Gamma_0:G\to \R^2$ be a smooth regular network such that $\gamma^i_*(p)=\gamma^i(p)$ for any endpoint $,p\in G$ and such that $\|\gamma^i_0-\gamma^i_*\|_{H^2(\d x)}\le \delta_{\Gamma_*}$. Then the motion by curvature $\Gamma_t:G\to \R^2$ starting from $\Gamma_0$ exists for all times and it smoothly converges to a minimal network $\Gamma_\infty$ such that ${\rm L}(\Gamma_\infty)={\rm L}(\Gamma_*)$, up to reparametrization.
\end{teo}

\begin{proof}
We recall the following interpolation inequalities. For any $k \in \N$ with $k\ge 1$ there exist $\lambda_k>0,\zeta_k\in(0,1)$ such that
\begin{equation}\label{eq:InterpolationSobolevGenerale}
    \|u \|_{H^k(\d x)} \le \lambda_k \| u \|_{L^2}^{\zeta_k} \| u \|_{H^{k+1}}^{1-\zeta_k},
\end{equation}
for any $u \in H^{k+1}\left((0,1);\R^N\right)$. We shall drop the subscript $k$ when $k=2$.

\medskip

Let $\sigma,\theta,r,C_{\rm LS}$ be given by \cref{thm:LojaGenerale} and \cref{cor:LojaNetworkMinimali}, where $C_{\rm LS}$ is the maximum of the constants given by both the statements.

\medskip

Recalling \cref{lemma:MinimalMinimizing}, up to decrease $r>0$, we can assume that the following hold. Whenever $\widehat\gamma^i\eqdef \gamma^i_*+ \NN^i\nu^i_*+\TT^i\tau^i_*$ is a smooth regular network, for $\overline{\NN}\in B_r(0)\subset V$ in the notation of \cref{thm:LojaGenerale}, where the $\TT^i$'s are adapted, then
\begin{enumerate}[label={\normalfont{\color{blue}\arabic*)}}]
    \item \label{it:r1} there exists a constant $C_G>2$, depending only on the graph $G$ and $\Gamma_*$, such that
    \begin{equation*}
    \begin{split}
        &\scal{\widehat\nu^i , \nu^i_*}\ge \frac34 , \qquad\qquad |\scal{\widehat\nu^i , \tau^i_*}| < \frac{1}{C_G}, \\
        &\sum_{m \in J_G}\sum_{\ell \in I_m} \left| a^\ell (x) \scal{\widehat\nu^\ell , \nu^\ell_*}
        +
        \scal{\widehat\nu^\ell , \tau^\ell_*}
        \chi(x)\mathscr{L}^\ell_m (a^{i_m}(x),a^{j_m}(x) )
        \right|^2
        \ge \frac{2}{C_G}\sum_i |a^i(x)|^2
        ,
    \end{split}
    \end{equation*}
    where $\widehat\nu^i$ is the normal vector of $\widehat\gamma^i$, and $i_m$ (resp. $j_m$) denotes the minimal (resp. intermediate) element of $I_m$, for any continuous functions $a^1,\ldots,a^N$;
    
    \item \label{it:r2} there exist $c_1,c_2>0$ such that
    \begin{equation*}
        c_1\le |\partial_x \widehat\gamma^i|^{-1} \le c_2,
    \end{equation*}
    for any $i$;
    
    \item \label{it:r3} there is $C_G'>2$ such that
    \begin{itemize}
        \item if $\Xi$ is a smooth regular network having the same endpoints of $\Gamma_*$ defined by immersions $\xi^i$ such that $\|\xi^i-\gamma^i_*\|_{H^2}< C_G' r$, then ${\rm L}(\Xi)\ge {\rm L}(\Gamma_*)$;
        
        \item $\|\widehat\gamma^i - \gamma^i_*\|_{H^2}< \min\{(C_G'-1) r , \sigma/2 \}$.
    \end{itemize}
\end{enumerate}

\medskip

We claim that whenever $\widehat\gamma^i_t=\gamma^i_*+ \NN^i_t\nu^i_*+\TT^i_t\tau^i_*$ is a smooth solution to the motion by curvature, for $\overline{\NN}_t\in B_r(0)\subset V$ for any $t$, where we used the notation of \cref{thm:LojaGenerale} and the $\TT_t^i$'s are adapted, then for any $m \in \N$ with $m\ge3$ there exists $C_m=C_m(r, \Gamma_*)>0$ such that
\begin{equation}\label{eq:UniformBoundNTempo}
    \left\| \overline{\NN}_t \right\|_{H^m(\d x)} \le C_m,
\end{equation}
for any $t$. The claim easily follows by combining the fact that $\overline{\NN}_t\in B_r(0)$ ensures a uniform $C^1$-bound on the parametrizations with the fact that uniform upper bounds on the $L^2(\d s)$-norm of the curvature along a motion by curvature imply uniform $L^2(\d s)$-bounds on every derivative of the immersion. The proof of \eqref{eq:UniformBoundNTempo} is postponed to the end of the proof.

\medskip 

Taking into account \cref{lemma:MinimalMinimizing} and \cref{cor:ParametrizzazioneTNtempo}, we can fix $\eta>0$ such that:
\begin{enumerate}[label={\normalfont{\color{blue}\roman*)}}]
    %\item\label{it:etaA} $\eta<\tfrac\sigma2$;
    
    \item \label{it:etaB} if immersions $\widehat\gamma^i$ define a regular network $\widehat\Gamma$ with same endpoints of $\Gamma_*$ such that $\|\widehat\gamma^i - \gamma^i_*\|_{C^0} \le 2\eta$, then ${\rm L}(\Gamma_*)\le {\rm L}(\widehat\Gamma)$;
    
    \item  \label{it:etaC} if $\widehat\gamma^i_t$ define a one-parameter family of immersions satisfying the assumptions of \cref{cor:ParametrizzazioneTNtempo} and $\sum_i \| \widehat\gamma^i_t - \gamma^i_t \|_{C^5} \le \eta $ for any $t$ around some $t_0$, then the resulting $\NN^i_t$ verify $\sum_i \|\NN^i_t\|_{H^4(\d x)}< \tfrac{r}{2}$ for any $t$ around $t_0$;
    
    \item \label{it:etaD} if immersions $\widehat\gamma^i$ define a network $\widehat\Gamma$ such that $\|\widehat\gamma^i - \gamma^i_*\|_{C^1} \le \eta$, then 
    \begin{equation*}
        \left|{\rm L}(\widehat\Gamma) - {\rm L}(\Gamma_*) \right|^\theta \le \frac{\theta r^{\frac1\zeta} }{C_{\rm LS} \sqrt{c_2C_G} \left(100\, \lambda C_3\right)^{\frac1\zeta}}.
    \end{equation*}
\end{enumerate}

With the above choices, we want to show that the statement follows by choosing
\begin{equation*}
    \delta_{\Gamma_*}\eqdef \overline{\eps}\bigg(\Gamma_*,\frac{\eta}{N},5\bigg),
\end{equation*}
where $\overline{\eps}$ is given by \cref{lem:Regularization2}.

\medskip

So let $\Gamma_0$ be as in the statement.
By \cref{lem:Regularization2}, the flow $\Gamma_t$ starting from $\Gamma_0$ satisfies 
\begin{equation}\label{eq:VicinanzaTau}
     \sum_i\|\gamma^i_\tau\circ\sigma^i-\gamma^i_*\|_{C^5} <  \eta,
\end{equation}
for some $\tau \in [0,T)$ and smooth reparametrizations $\sigma^i$. Then by \ref{it:etaB} we have ${\rm L}(\Gamma_\tau)\ge {\rm L}(\Gamma_*)$. Moreover, if ${\rm L}(\Gamma_\tau)= {\rm L}(\Gamma_*)$, then \ref{it:etaB} implies that $\Gamma_\tau$ is a local minimizer for the length in $C^0$, and thus it is minimal up to reparametrization, and the resulting flow is stationary. Hence we can assume that ${\rm L}(\Gamma_\tau)> {\rm L}(\Gamma_*)$.

Moreover, by \cref{cor:ParametrizzazioneTNtempo} and \ref{it:etaC} we get the existence of $\NN^i_t,\TT^i_t,\varphi^i_t$ as in \cref{cor:ParametrizzazioneTNtempo} such that
\begin{equation}\label{eq:wzwz}
    \gamma^i_t\circ \sigma^i \circ \varphi^i_t = \gamma^i_* + \NN^i_t \nu^i_* + \TT^i_t \tau^i_* \eqqcolon \widetilde\gamma^i_t,
\end{equation}
with
\begin{equation*}
    \sum_i \|\NN^i_t\|_{H^4(\d x)} < \frac{r}{2},
\end{equation*}
for any $t \in [\tau,\tau_1)$ with $\tau_1>\tau$.

We define the nonincreasing function
\begin{equation}\label{eq:DefH}
    H(t)\eqdef ({\rm L}(\Gamma_t)- {\rm L}(\Gamma_*) )^\theta,
\end{equation}
for $t \in [0,T)$.

Let us further define $S$ the supremum of all $s \in [\tau,T)$ such that $\gamma^i_t$ can be written as in \eqref{eq:wzwz} for some reparametrizations $\varphi_t$ and functions $\NN^i_t$ continuously differentiable in time with $\sum_i \|\NN^i_t\|_{H^2(\d x)} < r$ for any $t \in[\tau,s]$.\\
We have that $S\ge\tau_1>\tau$. Moreover, we can assume that ${\rm L}(\Gamma_s)>{\rm L}(\Gamma_*)$ for any $s \in [\tau,S)$. Indeed, if instead ${\rm L}(\Gamma_s)={\rm L}(\Gamma_*)$ for some $s$, then $\Gamma_s$ locally minimizes the length in $H^2$: if immersions $\bar\gamma^i$ define a smooth regular network with $\|\bar\gamma^i-\widetilde\gamma^i_s\|_{H^2}< r$, then $\|\bar\gamma^i-\gamma^i_*\|_{H^2} \le \|\bar\gamma^i-\widetilde\gamma^i_s\|_{H^2} + \|\widetilde\gamma^i_s-\gamma^i_*\|_{H^2} < C_G'r$ by \ref{it:r3}, and then ${\rm L}(\Gamma_s)={\rm L}(\Gamma_*) \le {\rm L}(\bar\Gamma)$ by \ref{it:r3}. Hence in this case $\Gamma_s$ is minimal, up to reparametrization, and the resulting flow is stationary.

Therefore we can assume $H(t)>0$ for $t \in(\tau,S)$, and then $H$ is differentiable on $(\tau,S)$. We now want to show that $S=T=+\infty$.

We differentiate
\begin{equation*}
    \begin{split}
        -\frac{\d}{\d t} H 
        &= \theta H^{\frac{\theta-1}{\theta}} \sum_i \int_0^1 |\boldsymbol{k}^i_t|^2 \de s 
        = \theta H^{\frac{\theta-1}{\theta}} \left(\sum_i \int_0^1 |\boldsymbol{k}^i_t|^2 \de s \right)^{\frac12} \left\| (\partial_t\Gamma_t)^\perp\right\|_{L^2(\d s)} \\
        &\ge \frac{\theta}{C_{\rm LS}}  \left\| (\partial_t\Gamma_t)^\perp\right\|_{L^2(\d s)},
    \end{split}
\end{equation*}
for any $t \in (\tau,S)$, where we denoted $ \left\| (\partial_t\Gamma_t)^\perp\right\|_{L^2(\d s)}^2 \eqdef \sum_i \int |(\partial_t\gamma^i_t)^\perp|^2 \de s = \sum_i  \int |(\partial_t\widetilde\gamma^i_t)^\perp|^2 \de s$, where we could apply the \L ojasiewicz--Simon inequality in \cref{cor:LojaNetworkMinimali} thanks to \ref{it:r3}. From the above estimate we get
\begin{equation*}
    \begin{split}
        -\frac{\d}{\d t} H 
        &\ge \frac{\theta}{C_{\rm LS}}\left( \sum_i  \int |(\partial_t\widetilde\gamma^i_t)^\perp|^2 \de s  \right)^{\frac12} \\
        &= \frac{\theta}{C_{\rm LS}}\left( \sum_i  \int |\partial_t \NN^i_t \scal{\nu^i_t,\nu^i_*} + \partial_t\TT^i_t \scal{\nu^i_t,\tau^i_*} |^2 \de s  \right)^{\frac12} \\
        &\ge \frac{\theta}{C_{\rm LS}}\left( \frac12\sum_{m \in J_G} \sum_{\ell \in I_m}  \int |\partial_t \NN^\ell_t \scal{\nu^\ell_t,\nu^\ell_*} + \partial_t\TT^\ell_t \scal{\nu^\ell_t,\tau^\ell_*} |^2 \de s  \right)^{\frac12} \\
        &\overset{\ref{it:r1}}{\ge}  \frac{\theta}{C_{\rm LS}\sqrt{C_G}} \left( \sum_i \int |\partial_t \NN^i_t|^2 \de s \right)^{\frac12} \\
        &\overset{\ref{it:r2}}{\ge}  \frac{\theta}{C_{\rm LS}\sqrt{C_G\, c_2}} \left( \sum_i \int |\partial_t \NN^i_t|^2 \de x \right)^{\frac12},
    \end{split}
\end{equation*}
for any $t \in (\tau,S)$. Hence
\begin{align}
        \left\|\overline{\NN}_s - \overline{\NN}_\tau\right\|_{L^2(\d x)} 
        &\nonumber = \left\|\int_\tau^s \partial_t \overline{\NN}_t \de t  \right\|_{L^2(\d x)} 
        \le \int_\tau^s \left\|  \partial_t \overline{\NN}_t  \right\|_{L^2(\d x)} \de t \\
        &\nonumber =  \int_\tau^s \left( \sum_i \int_0^1 |\partial_t \NN^i_t|^2 \de x \right)^{\frac12} \de t\\
        &\nonumber\le \frac{C_{\rm LS}\sqrt{C_G\, c_2}}{\theta} \left( H(\tau)-H(s) \right) \\
        & \le \frac{C_{\rm LS}\sqrt{C_G\, c_2}}{\theta} H(\tau), \label{eq:zwzEstimate}
\end{align}
for any $s \in (\tau,S)$. Recalling \eqref{eq:VicinanzaTau} and \ref{it:etaD}, we conclude that
\begin{equation*}
    \left\|\overline{\NN}_s - \overline{\NN}_\tau\right\|_{L^2(\d x)} \le \frac{r^{\frac1\zeta} }{\left(100\, \lambda C_3\right)^{\frac1\zeta}}.
\end{equation*}
for any $s \in (\tau,S)$. Exploiting the interpolation inequality \eqref{eq:InterpolationSobolevGenerale} with $k=2$ we obtain
\begin{equation*}
\begin{split}
     \left\|\overline{\NN}_s - \overline{\NN}_\tau\right\|_{H^2(\d x)} & \le \frac{r}{100\,  C_3} \left\|\overline{\NN}_s - \overline{\NN}_\tau\right\|_{H^3(\d x)}^{1-\zeta} \\
     &\overset{\eqref{eq:UniformBoundNTempo}}{\le}  \frac{r}{100\,  C_3} (2C_3)^{1-\zeta} \\
     &\le \frac{r}{50},
\end{split}
\end{equation*}
for any $s \in (\tau,S)$. Since $\|\overline{N}_\tau\|_{H^2(\d x)}<\tfrac{r}{2}$, a simple contradiction argument implies that $S=T$ and $\|\overline{N}_t\|_{H^2(\d x)}<\tfrac{r}{2}+\tfrac{r}{50}$ for any $t \in [\tau,T)$. Hence \cref{thm:LongTimeGeneral} implies that $T=+\infty$.

\medskip

We claim that $H(t)\searrow 0$ as $t\to+\infty$. Indeed, since $S=T=+\infty$, we now know that \eqref{eq:UniformBoundNTempo} holds for any time. Hence there exists a sequence of times $t_n\to+\infty$ such that the parametrizations $\widetilde\gamma^i_{t_n}$ converge in $C^2$ to limit parametrizations $\widetilde\gamma^i_\infty\eqdef \gamma^i_* + \NN^i_\infty \nu^i_* + \TT^i_\infty \tau^i_*$ with $\overline{N}_\infty \in B_r(0) \subset V$. Moreover, $\widetilde\gamma^i_\infty$ parametrize a minimal network $\widetilde\Gamma_\infty$. Hence using \ref{it:r3} and \cref{cor:LojaNetworkMinimali} the length of $\widetilde\Gamma_\infty$ has to be equal to the length of $\Gamma_*$. As $H$ is nonincreasing, then $H(t)\searrow 0$ as $t\to+\infty$.

\medskip

Exploiting the fact that $H(t)$ is infinitesimal as $t$ diverges, estimating as in \eqref{eq:zwzEstimate} for large times implies that the curve $\overline{\NN}_t$ is Cauchy in $L^2(\d x)$, and thus there exists its full limit $\overline{\NN}_\infty$ in $L^2(\d x)$ as $t\to+\infty$. Interpolating using \eqref{eq:InterpolationSobolevGenerale} and \eqref{eq:UniformBoundNTempo}, we then conclude that convergence holds in $H^m$ for any $m$. 

\medskip

We are now left to prove the claim \eqref{eq:UniformBoundNTempo}. For the sake of clarity, we consider the case $m=3$ only, the general case following by induction. We differentiate the curvature $\widehat{\boldsymbol{k}}^i_t$ of $\widehat\gamma^i_t$ and we multiply by the normal $\widehat\nu^i_t$ to get the identity
\begin{equation}\label{eq:zxc}
\begin{split}
    \scal{\partial_x \widehat{\boldsymbol{k}}^i_t,\widehat\nu^i_t} 
    &= \left\langle \partial_x \left(|\partial_x \widehat\gamma^i_t|^{-2} \partial_x^2 \widehat\gamma^i_t - |\partial_x \widehat\gamma^i_t|^{-4} \scal{\partial^2_x \widehat\gamma^i_t , \partial_x \widehat\gamma^i_t}\partial_x \widehat\gamma^i_t\right) ,
    \widehat\nu^i_t
    \right\rangle \\
    &= \left\langle |\partial_x \widehat\gamma^i_t|^{-2} \partial_x^3 \widehat\gamma^i_t 
    -2|\partial_x \widehat\gamma^i_t|^{-4} \scal{\partial^2_x \widehat\gamma^i_t , \partial_x \widehat\gamma^i_t} \partial_x^2 \widehat\gamma^i_t 
    - |\partial_x \widehat\gamma^i_t|^{-4} \scal{\partial^2_x \widehat\gamma^i_t , \partial_x \widehat\gamma^i_t}\partial_x^2 \widehat\gamma^i_t,
    \widehat\nu^i_t
    \right\rangle \\
    & = |\partial_x \widehat\gamma^i_t|^{-2} \scal{\partial_x^3\NN^i_t \nu^i_* + \partial_x^3 \TT^i_t \tau^i_* , \widehat\nu^i_t} + \\
    &\qquad- \left\langle 2|\partial_x \widehat\gamma^i_t|^{-4} \scal{\partial^2_x \widehat\gamma^i_t , \partial_x \widehat\gamma^i_t} \partial_x^2 \widehat\gamma^i_t 
    +|\partial_x \widehat\gamma^i_t|^{-4} \scal{\partial^2_x \widehat\gamma^i_t , \partial_x \widehat\gamma^i_t}\partial_x^2 \widehat\gamma^i_t,
    \widehat\nu^i_t
    \right\rangle.
\end{split}
\end{equation}
Taking absolute values and recalling \ref{it:r1}, \ref{it:r2}, we deduce that
\begin{equation*}
    \| \partial^3_x \NN^i_t \|_{L^1(\d x)} \le C(r,\Gamma_*)\left(1 + \int \left|\partial_s  \widehat{\boldsymbol{k}}^i_t\right| \de s \right),
\end{equation*}
where $C(r,\Gamma_*)>0$ here is a constant that may change from line to line.

Recalling \cite[Proposition 5.8]{mannovplusch}, we know that along a motion by curvature the $L^2(\d s)$-norms of derivatives of the curvature are bounded by the $L^2(\d s)$-norms of the curvature and by the inverse of the length of the edges. Hence the assumption $\NN_t \in B_r(0)\subset V$ guarantees that $\int |\partial_s  \widehat{\boldsymbol{k}}^i_t| \de s \le C(r,\Gamma_*)$. In particular $\| \NN^i_t \|_{W^{3,1}(\d x)} \le C(r,\Gamma_*)$, and thus $\|\NN^i_t \|_{W^{2,\infty}(\d x)}\le  C(r,\Gamma_*)$. Therefore we can improve the estimate on $\partial_x^3 \NN^i_t$ by first taking squares and then integrating in \eqref{eq:zxc}, which yields
\begin{equation*}
    \|\partial_x^3 \NN^i_t\|_{L^2(\d x)} \le C(r,\Gamma_*),
\end{equation*}
thus proving the claim \eqref{eq:UniformBoundNTempo}.
\end{proof}

An immediate consequence is the next result, which promotes subconvergence of the motion by curvature to full convergence.

\begin{teo}\label{thm:Convergence}
Let $\Gamma_t: G\to \R^2$ be a smooth motion by curvature defined on $[0,+\infty)$. Let $\Gamma_\infty:G\to \R^2$ be a minimal network such that $\Gamma_{t_n} \to \Gamma_\infty$ in $H^2$ for some sequence $t_n\nearrow +\infty$ as $n\to+\infty$.
Then $\Gamma_{t} \to \Gamma_\infty$ smoothly as $t\to+\infty$, up to reparametrization.
\end{teo}

\begin{proof}
The statement immediately follows from \cref{thm:Stability}.
\end{proof}

We conclude this part by collecting some observations implied by the previous stability results.

\begin{rem}
Theorem~\ref{thm:Convergence} 
can be combined with~\cite[Proposition 13.5]{mannovplusch} in the following way. If $\Gamma_t: G\to \R^2$ is a motion by curvature of a tree-like network, i.e., $G$ has no cycles, defined on $[0,+\infty)$, if the sequential limit $\Gamma_\infty$ along a sequence of times $t_n$, which always exists by \cite[Proposition 13.5]{mannovplusch}, is regular, then $\Gamma_\infty$ is the full limit of $\Gamma_t$ as $t\to+\infty$.\\
However, the example in the next section
shows that in general the limit $\Gamma_\infty$ may be degenerate.
\end{rem}

\begin{rem}
If the network $\Gamma_*$ in Theorem~\ref{thm:Stability} is an isolated critical point of the length, then 
$\Gamma_\infty$ coincides with $\Gamma_*$. 
This is always the case if $\Gamma_*$ is a tree, i.e., $G$ has no cycles, since there exist finitely many minimal trees $\widehat{\Gamma}:G\to\R^2$ having the same endpoints of $\Gamma_*$.
\end{rem}

\begin{rem}
In the notation of \cref{thm:Stability}, in some cases we are able to conclude that $\Gamma_\infty$ coincides with $\Gamma_*$, even if $\Gamma_*$ is not an isolated critical point of the length.

Suppose that $\Gamma_*$
is a minimal network composed of a regular hexagon $H$ with area $A_*$ and six straight segments connecting the vertices of a bigger regular hexagon. Then $\Gamma_*$ is not an isolated critical point of the length, indeed
there exists a one-parameter family of critical points with the same length:
all networks composed of concentric hexagons and straight segments connecting the endpoints, see \cref{Fig:Ragnatela}. It can be proved that there are no other minimal networks
with this topology and with the same endpoints. 
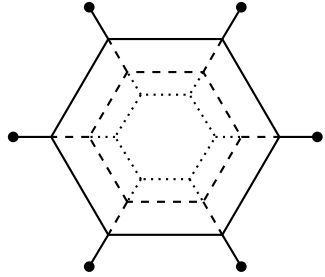
\begin{figure}[H]
	\begin{center}
	\begin{tikzpicture}[scale=2]
\draw[thick]
(-0.5,-0.86)--(-0.375,-0.649)
(-1,0)-- (-0.75,0)
(-0.5,0.86)--(-0.375,0.649)
(0.5,-0.86)--(0.375,-0.649)
(1,0)-- (0.75,0)
(0.375,0.649)--(0.5,0.86) ;
\draw[thick]
(-0.375,-0.649)--(-0.75,0)--(-0.375,0.649)--
(0.375,0.649)--(0.75,0)--(0.375,-0.649)--(-0.375,-0.649);
\draw[thick,dashed]
(-0.25,-0.43)--(-0.375,-0.649)
(-0.5,0)-- (-0.75,0)
(-0.25,0.43)--(-0.375,0.649)
(0.25,-0.43)--(0.375,-0.649)
(0.5,0)-- (0.75,0)
(0.375,0.649)--(0.25,0.43);
\draw[thick,dashed]
(-0.25,-0.43)--(-0.5,0)--(-0.25,0.43)--
(0.25,0.43)--(0.5,0)--(0.25,-0.43)--(-0.25,-0.43);
\draw[thick,dotted]
(-0.16,-0.28)--(-0.33,0)--(-0.16,0.28)--
(0.16,0.28)--(0.33,0)--(0.16,-0.28)--(-0.16,-0.28);
\draw[thick,dotted]
(-0.16,-0.28)--(-0.25,-0.43)
(-0.33,0)-- (-0.5,0)
(-0.16,0.28)--(-0.25,0.43)
(0.16,-0.28)--(0.25,-0.43)
(0.33,0)-- (0.5,0)
(0.25,0.43)--(0.16,0.28) ;
\fill[black](1,0) circle (1pt);  
\fill[black](-1,0) circle (1pt);  
\fill[black](0.5,-0.86) circle (1pt);  
\fill[black](-0.5,-0.86) circle (1pt);  
\fill[black](0.5,0.86) circle (1pt);  
\fill[black](-0.5,0.86) circle (1pt);  
\end{tikzpicture}
	\end{center}
	\caption{Three different minimal networks with the same endpoints and topology. All these networks have the same length.}\label{Fig:Ragnatela}
\end{figure}
In the above notation, suppose now that $\Gamma_0$
is regular network with the same endpoints 
and the same topology of 
$\Gamma_*$, sufficiently close to $\Gamma_*$ in $H^2$,
and such that the area enclosed
by the loop equals $A_*$.
Then $\Gamma_\infty$ coincides with $\Gamma_*$.
Indeed the area enclosed by a loop composed of six curves is preserved during the evolution
(see~\cite[Section 8.2]{mannovplusch}) and $\Gamma_*$ is the unique minimal network
with area $A_*$ among the one-parameter family of possible minimal networks.
\end{rem}

\section{Convergence to a degenerate network in infinite time}
\label{sec:Example}

In this section we construct an example of a motion by curvature existing for all times, with uniformly bounded curvature, smoothly converging to a degenerate network. More precisely, there holds the following result.

\begin{teo}\label{thm:EsempioConvergenceInfiniteTime}
There exists a smooth regular network $\Gamma_0:G\to \R^2$ such that the motion by curvature $\Gamma_t$ starting from $\Gamma_0$ exists for every time, the length of each curve $\gamma^i_t$ is strictly positive for any time, the curvature of each curve $\gamma^i_t$ is uniformly bounded from above, and $\Gamma_t$ smoothly converges to a degenerate network $\Gamma_\infty$ as $t\to+\infty$, up to reparametrization. Specifically, the length of a distinguished curve $\gamma^0_t$ tends to zero as $t\to+\infty$.
\end{teo}

\begin{proof}
The proof of the statement follows by putting together the observations in \ref{Step1}, \ref{Step2}, and \ref{Step3} below.
\end{proof}

From now on and for the rest of this section, let $\Gamma_0:G\to \R^2$ be a smooth regular network as in \cref{Fig:Network}. We assume that $\Gamma_0$ is composed of five curves, it is symmetric with respect to horizontal and vertical axes, the middle curve $\gamma^0$ is a segment, and the remaining four curves are convex, i.e., their oriented curvature has a sign. Moreover, the network has four endpoints located at the vertices of a rectangle of sides of length $2/\sqrt{3}$ and $2$, so that the diagonals of the rectangle meets forming angles of $\tfrac23 \pi$ and $\tfrac\pi3$, see \cref{Fig:Network}.

We want to show that the motion by curvature $\Gamma_t$ starting from such a datum $\Gamma_0$ satisfies the statement of \cref{thm:EsempioConvergenceInfiniteTime}. The candidate limit is given by the degenerate network defined by the diagonals of the rectangle, that is, the dotted lines in \cref{Fig:Network}.

\medskip

\begin{figure}[H]
	\begin{center}
		\begin{tikzpicture}[scale=5]
		\draw
		(-1,-0.5)to[out=5, in=210, looseness=0.8]
		(0,-0.15)to[out=90, in=270, looseness=0]
		(0,0.15);
		\draw[xscale=-1]
		(-1,-0.5)to[out=5, in=210, looseness=0.8]
		(0,-0.15);
		\path[font=\normalsize]
		(-1,-0.5)node[]{$\bullet$};
		\path[font=\normalsize]
		(1,-0.5)node[]{$\bullet$};
		\draw [decorate,decoration={brace,amplitude=5pt,mirror,raise=1ex}]
        (-1,-0.5) -- (0,-0.5)
        node[midway,yshift=-1.5em]{$1$};
        \draw [decorate,decoration={brace,amplitude=5pt,raise=1ex}]
        (-1,-0.5) -- (-1,0)
        node[midway,xshift=-1.7em]{$\frac{1}{\sqrt{3}}$};
        \draw[dotted]
        (-1,0)--(0,0);
        \draw[dotted]
        (0,-0.5)--(0,-0.15);
        %\path[font=\normalsize]
	    %(0,0)node[right]{$0$};
		\draw[cm={cos(-30) ,-sin(-30) ,sin(-30) ,cos(-30) ,(0 cm,-0.15 cm)},-stealth](0,0) -- (0.2,0);
		\draw[xscale=-1, cm={cos(-30) ,-sin(-30) ,sin(-30) ,cos(-30) ,(0 cm,-0.15 cm)},-stealth](0,0) -- (0.2,0);
		\draw[-stealth](0,0)--(0,0.1);
		\draw[cm={cos(-120) ,-sin(-120) ,sin(-120) ,cos(-120) ,(0 cm,-0.15 cm)},-stealth](0,0) -- (0.2,0);
		\draw[cm={cos(-240) ,-sin(-240) ,sin(-240) ,cos(-240) ,(0 cm,-0.15 cm)},-stealth](0,0) -- (0.2,0);
		\draw[cm={cos(-90) ,-sin(-90) ,sin(-90) ,cos(-90) ,(0.8 cm,-0.05 cm)},-stealth](0,0) -- (0.2,0);
		\path[font=\normalsize]
		(0.8,0)node[right]{$\omega$};
		\path[font=\normalsize]
		(-0.7,-0.45)node[above]{$\gamma^1_t$};
		\path[font=\normalsize]
		(0.7,-0.45)node[above]{$\gamma^2_t$};
		\path[font=\normalsize]
		(0,0.14)node[right]{$\gamma^0_t$};
		\path[font=\small]
		(0.25,-0.05)node[]{$\tau^1_t(1)$};
		\path[font=\small]
		(-0.25,0)node[below]{$\tau^2_t(1)$};
		\path[font=\small]
		(-0.15,0.01)node[above]{$\nu^1_t(1)$};
		\path[font=\small]
		(-0.15,-0.3)node[below]{$\nu^2_t(1)$};
		\end{tikzpicture}
	\end{center}
	\caption{Bottom half of the motion by curvature starting from a network as in \cref{Fig:Network}, specifying notation and orientation of the edges.}\label{Fig:Symbols}
\end{figure}
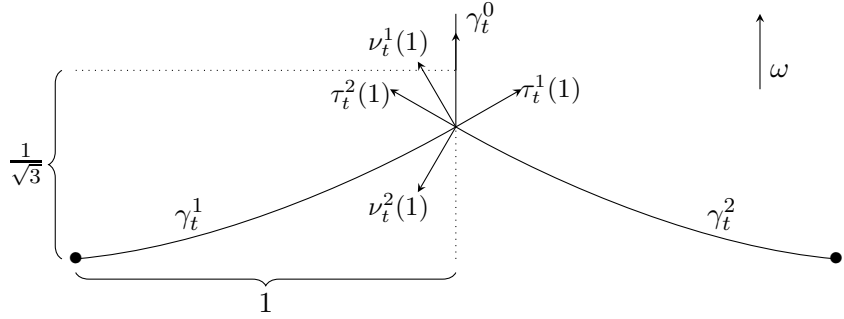

By symmetry, it is sufficient to study the evolution of the middle curve and of the two bottom curves in \cref{Fig:Network}. To fix the notation, we recall such part of the graph in \cref{Fig:Symbols}. Observe that the straight middle curve $\gamma^0$ is parametrized from the bottom to the top, while the convex curves $\gamma^1, \gamma^2$ have the endpoint $1$ at the junction. This is in contrast with the usual choice we adopted of setting endpoints $1$ at the endpoints of the network, however we choose this parametrization here in order to simplify useless presence of minus signs in the computations below. Finally, we denote by
\begin{equation*}
    \omega\eqdef (0,1),
\end{equation*}
the vertical unit vector, coinciding with the tangent vector of the curve $\gamma^0$.

Recalling \cref{rem:SoluzioniConfrontoCompatibilita}, we can assume that the motion by curvature is smooth and evolves by the special flow, i.e., $\partial_t\gamma^i_t= |\partial_x\gamma^i_t|^{-2} \partial^2_x \gamma^i_t$ for any $i$. Decomposing $\partial_t \gamma^i_t$ in tangential and normal components, we denote
\begin{equation*}
\partial_t\gamma^i_t = \widetilde{k}_i \nu^i_t + \lambda_i \tau^i_t,    
\end{equation*}
where we denote by $\widetilde{k}_i$ the oriented curvature of $\gamma^i_t$, i.e., $\widetilde{k}_i\eqdef \scal{\boldsymbol{k}^i_t, \nu^i_t}$. We drop subscript $t$ in $\widetilde{k}_i$ and $\lambda_i$ for ease of notation.

At least for short times, by choice of the initial datum, we can consider the functions $v_i$ defined by
\begin{equation*}
    v_i\eqdef \frac{1}{\scal{\nu^i_t,\omega}},
\end{equation*}
for $i=1,2$. We further assume that
\begin{equation}\label{eq:InclinazioneTauInizio}
    \scal{\tau^1_t(0),\omega} \big|_{t=0} >0.
\end{equation}

We preliminarily observe that, by symmetry and choice of orientations, we have $\widetilde{k}_1=-\widetilde{k}_2$ and $\partial_s \widetilde{k}_1 = - \partial_s \widetilde{k}_2$ at any time and point. Moreover, symmetry and evolution of curvature imply that $\gamma^0_t$ is a vertical segment for any time; then $\partial_t\gamma^0_t(t,0)$ and $\omega$ are parallel, hence $\lambda_1(t,1)= \scal{\partial_t\gamma^0_t(t,0), \tau^1_t} = \scal{\partial_t\gamma^0_t(t,0), \tau^2_t}   = \lambda_2(t,1)$ for any $t\in[0,T)$. On the other hand, the boundary condition obtained by the derivative $\partial_t \scal{\tau^1_t(1),\tau^2_t(1)}=0$, see \cite{mannovplusch}, reads
\begin{equation*}
    \partial_s \widetilde{k}_1(t,1) + \lambda_1(t,1)  \widetilde{k}_1(t,1) =    \partial_s \widetilde{k}_2(t,1) + \lambda_2(t,1)  \widetilde{k}_2(t,1) .
\end{equation*}
Therefore we get that
\begin{equation}\label{eq:ConditionTripuntoSimmetria}
    \partial_s \widetilde{k}_i(t,1) + \lambda_i(t,1)  \widetilde{k}_i(t,1) = 0,
\end{equation}
for $i=1,2$ for any $t \in [0,T)$. Finally, recalling from \cite[Section 3]{mannovplusch} that tangential velocities at a junction can be expressed in terms of normal velocities, which easily follows from identity $\partial_t \gamma^1_t(1)=\partial_t \gamma^2_t(1)$, we have that
\begin{equation}\label{eq:TangVelviaNormVel}
    \lambda_1(t,1) = - \frac{\widetilde{k}_2(t,1)}{\sqrt{3}} = \frac{\widetilde{k}_1(t,1)}{\sqrt{3}} ,
    \qquad\qquad
    \lambda_2(t,1) =  \frac{\widetilde{k}_1(t,1)}{\sqrt{3}} = -\frac{\widetilde{k}_2(t,1)}{\sqrt{3}} .
\end{equation}

\begin{enumerate}[label={\normalfont{\color{blue}Step \arabic*}}, wide, labelwidth=!, labelindent=10pt]

\item\label{Step1} Letting $T>0$ the maximal time of existence of the flow, we want to prove that the functions $v_i$ are defined on $[0,T)$ and
\begin{equation}\label{eq:Step1}
    \widetilde{k}_1 \ge 0 ,\qquad\qquad
    1\le v_1 \le \frac{2}{\sqrt{3}},
\end{equation}
for any $x \in [0,1]$ and $t \in[0,T)$. In particular, the curves $\gamma^1_t, \gamma^2_t$ can be parametrized by convex graphs on a fixed interval for any time.

\medskip

By basic computations on the evolution of geometric quantities, see \cite{mannovplusch, DzKuSc02}, one easily obtains
%\begin{equation*}
%\begin{split}
%    \partial_t v_1 &= (v_1)^2 \left(\partial_s\widetilde{k}_1 + \lambda_1\widetilde{k}_1 \right)\scal{\tau^1_t, \omega},\\
%    \partial_s v_1 &= (v_1)^2 \widetilde{k}_1 \scal{\tau^1_t, \omega}, \\
%    \partial_s^2 v_1 &= 2v_1 \partial_s v_1 \,  \widetilde{k}_1\scal{\tau^1_t,\omega} + (v_1)^2 \partial_s  \widetilde{k}_1 \, \scal{\tau^1_t, \omega} + (v_1)^2 ( \widetilde{k}_1)^2 \scal{\nu^1_t,\omega} \\
%    &= 2 \frac{(\partial_s v_1)^2}{v_1} + (v_1)^2 \partial_s  \widetilde{k}_1 \, \scal{\tau^1_t, \omega} + v_1 ( \widetilde{k}_1)^2 .
%\end{split}
%\end{equation*}
%Hence
\begin{equation}\label{eq:CaloreV1}
    (\partial_t - \partial^2_s) v_1 = -v_1(\widetilde{k}_1)^2 -2 \frac{(\partial_s v_1)^2}{v_1}  + \lambda_1 \partial_s v_1.
\end{equation}
Recalling that
\begin{equation}\label{eq:CaloreK}
    (\partial_t - \partial^2_s) \widetilde{k}_1 = \lambda_1 \partial_s \widetilde{k}_1 + (\widetilde{k}_1)^3,
\end{equation}
we obtain
\begin{equation*}
    \begin{split}
         (\partial_t - \partial^2_s)(v_1 \widetilde{k}_1)
         %&= \lambda_1 \partial_s(v_1 \widetilde{k}_1) -2 \frac{(\partial_s v_1)^2}{v_1} \widetilde{k}_1 - 2 (\partial_s \widetilde{k}_1) (\partial_s v_1) \\
         %&=\lambda_1\partial_s(v_1 \widetilde{k}_1) - 2 \frac{\partial_s v_1}{v_1} \partial_s (v_1 \widetilde{k}_1) \\
         %&
         = \left[\lambda_1  -2v_1 \widetilde{k}_1 \scal{\tau^1_t, \omega} \right]\partial_s (v_1 \widetilde{k}_1) .
    \end{split}
\end{equation*}
Exploiting \eqref{eq:ConditionTripuntoSimmetria} and \eqref{eq:TangVelviaNormVel}, we see that $(v_1 \widetilde{k}_1)$ satisfies a Neumann boundary condition at $x=1$, that is
\begin{equation*}
\begin{split}
    \partial_s(v_1 \widetilde{k}_1) \big|_{x=1}
    &= v_1 \partial_s \widetilde{k}_1 + (\widetilde{k}_1)^2 (v_1)^2 \scal{\tau^1_t,\omega} \, \big|_{x=1} 
    = -\frac{2}{\sqrt{3}}\lambda_1 \widetilde{k}_1 + (\widetilde{k}_1)^2 \left( \frac{2}{\sqrt{3}}\right)^2 \frac12 \, \bigg|_{x=1}  \\
    &= -\frac23 (\widetilde{k}_1)^2  + \frac23 (\widetilde{k}_1)^2 \, \bigg|_{x=1}   = 0.
\end{split}
\end{equation*}
Let $\overline{T}\le T$ be the maximal time such that $v_1$ is well defined. For $\eps,\delta>0$, we consider the function $f\eqdef v_1 \widetilde{k}_1 + \eps t + \delta$. By the above observations and since $\widetilde{k}_1(t,0)=0$, then $f$ satisfies
\begin{equation*}
    \begin{cases}
    (\partial_t - \partial^2_s) f = \left[\lambda_1  -2v_1 \widetilde{k}_1 \scal{\tau^1_t, \omega} \right]\partial_s f + \eps & \text{ on } [0,\overline{T}) \times [0,1], \\
    f(0,x) \ge \delta & \forall\, x \in [0,1],\\
    f(t,0) \ge \delta & \forall\, t \in [0,\overline{T}), \\
    \partial_s f (t,1) = 0 & \forall\, t \in [0,\overline{T}).
    \end{cases}
\end{equation*}
By a standard argument involving the maximum principle, we can prove that $f>0$ at any $(t,x) \in [0,\overline{T})\times [0,1]$. More precisely, if $\overline{t}>0$ is the first time such that there is $\overline{x}$ such that $f(\overline{t}, \overline{x}) =0$, then $\overline{x} \in(0,1]$. The case $\overline{x}=1$ is excluded as Hopf Lemma (see \cite[Theorem 6, p. 174]{ProtterWein}) would imply $\partial_s f(\overline{t},1) <0$. Also the case $\overline{x}\in(0,1)$ leads to contradiction, as in this case $0 \ge (\partial_t - \partial^2_s) f (\overline{t},\overline{x}) \ge \eps >0 $.

Arbitrariness of $\eps,\delta$ implies that $ v_1 \widetilde{k}_1\ge 0$ on $[0,\overline{T})\times [0,1]$. Since by continuity $v_1$ must be strictly positive on $[0,\overline{T})\times [0,1]$, then $\widetilde{k}_1\ge0$ on $[0,\overline{T})\times [0,1]$. Since convexity is preserved up to time $\overline{T}$ and recalling assumption \eqref{eq:InclinazioneTauInizio}, then
\begin{equation*}
    \begin{split}
        \partial_t \scal{\nu^1_t,\omega}|_{x=0} &= -\partial_s \widetilde{k}_1 \scal{\tau^1_t,\omega}|_{x=0}\le 0, \\
        \partial_s\scal{\nu^1_t,\omega} &= \scal{- \widetilde{k}_1 \tau^1_t, \omega} \le 0,
    \end{split}
\end{equation*}
where we used that $\partial_s\widetilde{k}_1|_{x=0}\ge0$ since $\widetilde{k}_1(0)=0$ is a global minimum for $\widetilde{k}_1$. Therefore the minimum of $\scal{\nu^1_t,\omega}$ is achieved at $x=1$, that is $\tfrac{\sqrt{3}}{2}= \scal{\nu^1_t,\omega}|_{x=1} \le \scal{\nu^1_t,\omega}\le 1$. The positive lower bound on $ \scal{\nu^1_t,\omega}$ implies that $\overline{T}=T$ and completes the proof of the first step.

\item\label{Step2} We claim that there exists a constant $C>0$ such that $\widetilde{k}_1 \le C$ for any $t\in[0,T)$. Moreover, for any $k\ge 1$ there is $C_k>0$ such that $\partial_s^k\widetilde{k}_1 \le C_k$ for any $t\in[0,T)$.

\medskip

By the evolution equations for $v_1$ and $\widetilde{k}_1$, we can compute
\begin{equation}\label{eq:zzxxz}
    \begin{split}
        (\partial_t - \partial^2_s) \left((v_1)^2(\widetilde{k}_1)^2\right)
        %&= 2 \Big(\lambda_1 (v_1)^2\widetilde{k}_1 \partial_s \widetilde{k}_1 
        %-(v_1)^2 (\partial_s \widetilde{k}_1)^2 - 3(\partial_s v_1)^2 (\widetilde{k}_1)^2 +\\
        %&\qquad+ \lambda_1  ( \widetilde{k}_1)^2 v_1 \partial_s v_1  - \partial_s( \widetilde{k}_1^2) \,\partial_s (v_1^2)
        %\Big) \\
        &= 2 \Big(
        \tfrac12\lambda_1 \partial_s \left((v_1)^2(\widetilde{k}_1)^2\right) + \\
        &\qquad
        -(v_1)^2 (\partial_s \widetilde{k}_1)^2 - 3(\partial_s v_1)^2 (\widetilde{k}_1)^2  - \partial_s( \widetilde{k}_1^2) \,\partial_s (v_1^2)
        \Big).
    \end{split}
\end{equation}
By Young inequality we estimate
\begin{equation*}
    \begin{split}
        -2\partial_s( \widetilde{k}_1^2) \,\partial_s (v_1^2) 
        &= - \partial_s( \widetilde{k}_1^2) \,\partial_s (v_1^2) - 4 v_1  \widetilde{k}_1 (\partial_s v_1) (\partial_s  \widetilde{k}_1) \\
        &=  - \partial_s \Big(v_1^2  \widetilde{k}_1^2 \Big) \partial_s(v_1^2) \, v_1^{-2} + ( \widetilde{k}_1)^2 v_1^{-2} \big(\partial_s(v_1^2) \big)^2
        - 4 v_1  \widetilde{k}_1 (\partial_s v_1) (\partial_s  \widetilde{k}_1) \\
        &= -2 v_1^{-1} \,\partial_s  v_1 \,  \partial_s \Big(v_1^2  \widetilde{k}_1^2 \Big) 
        + 4 ( \widetilde{k}_1)^2  \big(\partial_s v_1 \big)^2
        - 4 v_1  \widetilde{k}_1 (\partial_s v_1) (\partial_s  \widetilde{k}_1) \\
        &\le -2 v_1^{-1} \,\partial_s  v_1 \,  \partial_s \Big(v_1^2  \widetilde{k}_1^2 \Big) 
        + 4 ( \widetilde{k}_1)^2  \big(\partial_s v_1 \big)^2 
        + 2(v_1)^2 (\partial_s  \widetilde{k}_1)^2 + 2(\widetilde{k}_1)^2 (\partial_s v_1)^2 \\
        &= 2 \Big( 
        - v_1^{-1} \,\partial_s  v_1 \,  \partial_s \Big(v_1^2  \widetilde{k}_1^2 \Big) 
        + 3 ( \widetilde{k}_1)^2  \big(\partial_s v_1 \big)^2 
        + (v_1)^2 (\partial_s  \widetilde{k}_1)^2
        \Big).
    \end{split}
\end{equation*}
Inserting in \eqref{eq:zzxxz} we get
\begin{equation}\label{eq:CaloreVquadroKquadro}
    \begin{split}
        (\partial_t - \partial^2_s) \left((v_1)^2(\widetilde{k}_1)^2\right)
        &\le 
         2 \Big(
        \tfrac12\lambda_1 \partial_s \left((v_1)^2(\widetilde{k}_1)^2\right) -v_1^{-1} \,\partial_s  v_1 \,  \partial_s \left((v_1)^2(\widetilde{k}_1)^2\right) \Big) \\
        &= \left[
        \lambda_1 -2 v_1^{-1} \partial_s v_1
        \right] \partial_s \left((v_1)^2(\widetilde{k}_1)^2\right).
    \end{split}
\end{equation}
Observe that $v_1=v_2$ by symmetry, hence all the above considerations hold for $v_2$ as well. We further consider
\begin{equation*}
    g_i\eqdef (\widetilde{k}_i)^2 (v_i)^2,
\end{equation*}
for $i=1,2$. Again, actually $g_1=g_2$ by symmetry. Observe that
\begin{equation}\label{eq:Boundary1}
    g_1(t,0)=\partial_s g_1(t,0) = 0,
\end{equation}
as $\widetilde{k}_1(t,0)=0$, for any $t \in [0,T)$. Moreover
\begin{equation}\label{eq:Boundary2}
\begin{split}
    \partial_s g_1(t,1) &= 2 \Big( \widetilde{k}_1 (\partial_s \widetilde{k}_1) (v_1)^2 + (\widetilde{k}_1)^2 v_1 (\partial_s v_1)  \Big) \, \Big|_{(t,1)} \\
    &= 2 \Big( \widetilde{k}_1 (\partial_s \widetilde{k}_1)  (v_1)^2  + (\widetilde{k}_1)^2  v_1 (v_1)^2 \widetilde{k}_1 \scal{\tau^1_t,\omega}  \Big) \, \Big|_{(t,1)}   \\
    &= 2 \Big( \widetilde{k}_1 (\partial_s \widetilde{k}_1)  (2/\sqrt{3})^2  + (\widetilde{k}_1)^3  (2/\sqrt{3})^3 (1/2) \Big) \, \Big|_{(t,1)} \\
    &= \frac{8}{3}\widetilde{k}_1 \Big(  (\partial_s \widetilde{k}_1)   + (\widetilde{k}_1)^2 /\sqrt{3}  \Big) \, \Big|_{(t,1)} =0,
\end{split}
\end{equation}
where the last equality follows from \eqref{eq:ConditionTripuntoSimmetria} and \eqref{eq:TangVelviaNormVel}. Obviously, $\partial_s g_2(t,1) =0$ as well.

Now take $t_0 \in (0,T)$ and let $p_0\in\R^2$ be the mid point of the image of the straight edge $\gamma^0$. Without loss of generality we can assume that $p_0=0$ is the origin of $\R^2$. Hence let
\begin{equation*}
    \rho(t,p) \eqdef \frac{1}{\sqrt{4\pi (t_0-t)}} \exp \left(-\frac{|p|^2}{4(t_0-t)} \right).
\end{equation*}
Denoting $\rho\circ \gamma^i_t\eqdef \rho(t,\gamma^i_t)$, we observe that
\begin{equation}\label{eq:Boundary3}
    \begin{split}
        -\partial_s(\rho \circ \gamma^1_t) \,\big|_{(t,1)}
        &= -\left\langle \nabla\rho |_{(t,\gamma^1_t(1))} , \tau^1_t(1)\right\rangle 
        = \frac{\rho\circ \gamma^1_t}{2(t_0-t)} \scal{\gamma^1_t(1), \tau^1_t(1)} \le 0,
    \end{split}
\end{equation}
for any $t \in (0,t_0)$, where the inequality follows by the choice of the origin of $\R^2$.

Now let $A\eqdef \max_{[0,1]} (\widetilde{k}_1)^2 (v_1)^2\, \big|_{t=0}>0$ and define
\begin{equation*}
    f_i(t,x) \eqdef \left(\max \left\{  (\widetilde{k}_i)^2 (v_i)^2 -A, 0 \right\} \right)^2.
\end{equation*}
Since $F(y)\eqdef \left(\max \left\{ y -A, 0 \right\} \right)^2$ is of class $C^{1,1}$, then $f_i(t,\cdot) \in H^2$ for any $t$ and chain rule holds almost everywhere, i.e., $\partial_s f_i =2 \max \left\{  (\widetilde{k}_i)^2 (v_i)^2 -A, 0 \right\} \partial_s((\widetilde{k}_i)^2 (v_i)^2)$ and $\partial_s^2 f_i =2\left[ \partial_s((\widetilde{k}_i)^2 (v_i)^2) \right]^2 + 2 \max \left\{  (\widetilde{k}_i)^2 (v_i)^2 -A, 0 \right\} \partial_s^2((\widetilde{k}_i)^2 (v_i)^2)$ almost everywhere. Analogously, $f_i$ is differentiable with respect to $t$ at any $(t,x)$ and $\partial_t f_i  = 2 \max \left\{  (\widetilde{k}_i)^2 (v_i)^2 -A, 0 \right\} \partial_t((\widetilde{k}_i)^2 (v_i)^2)$ is continuous on $[0,T)\times [0,1]$.

Recalling \eqref{eq:CaloreVquadroKquadro} and using Young inequality we estimate
\begin{equation}\label{eq:CaloreF1}
    \begin{split}
        (\partial_t - \partial_s^2) f_1
        &=2 \max \left\{  (\widetilde{k}_1)^2 (v_1)^2 -A, 0 \right\} (\partial_t - \partial_s^2) \big( (\widetilde{k}_1)^2 (v_1)^2\big) -2\left[ \partial_s((\widetilde{k}_1)^2 (v_1)^2) \right]^2 \\
        &\le 2 \max \left\{  (\widetilde{k}_1)^2 (v_1)^2 -A, 0 \right\} \left[
        \lambda_1 -2 v_1^{-1} \partial_s v_1
        \right] \partial_s \left((v_1)^2(\widetilde{k}_1)^2\right)  -2\left[ \partial_s((\widetilde{k}_1)^2 (v_1)^2) \right]^2\\
        &\le \frac12\left[
        \lambda_1 -2 v_1^{-1} \partial_s v_1
        \right]^2 f_1,
    \end{split}
\end{equation}
for any $t$ and almost every $x$. We apply the monotonicity-type formula from \cref{thm:Monotonicity} with $f=f_1$ to get
\begin{equation*}
\begin{split}
        \frac{\d}{\d t}  \int_0^1 ( \rho\circ \gamma^1_t) \, f_1 \de s
    &\le \int_0^1 ( \rho\circ \gamma^1_t)(\partial_t-\partial_s^2)f_1   +
    \int_0^1
      \left( \partial_s \lambda_1
     -\frac{\lambda_1}{2(t_0-t)}\scal{\gamma^1_t,\tau^1_t}  \right) ( \rho\circ \gamma^1_t) \, f_1
     \de s
    + \\
    &\qquad+\big( ( \rho\circ \gamma^1_t)\partial_sf_1 - f_1 \partial_s ( \rho\circ \gamma^1_t)\big)\bigg|_0^1,
\end{split}
\end{equation*}
for any $t\in(0,t_0)$. Employing \eqref{eq:Boundary1}, \eqref{eq:Boundary2}, \eqref{eq:Boundary3}, and \eqref{eq:CaloreF1}, we obtain
\begin{equation}\label{eq:zzEvolutionComparison}
    \begin{split}
         \frac{\d}{\d t}  \int_0^1 ( \rho\circ \gamma^1_t) \, f_1 \de s
    &\le \int_0^1
      \left(\frac12\left[
        \lambda_1 -2 v_1^{-1} \partial_s v_1
        \right]^2+ \partial_s \lambda_1
     -\frac{\lambda_1}{2(t_0-t)}\scal{\gamma^1_t,\tau^1_t}  \right) ( \rho\circ \gamma^1_t) \, f_1
     \de s
    + \\
    &\qquad- f_1(t,1) \partial_s(\rho\circ \gamma^1_t)\big|_{(t,1)} \\
    &\overset{\eqref{eq:Boundary3}}{\le} \int_0^1
      \left(\frac12\left[
        \lambda_1 -2 v_1^{-1} \partial_s v_1
        \right]^2+ \partial_s \lambda_1
     -\frac{\lambda_1}{2(t_0-t)}\scal{\gamma^1_t,\tau^1_t}  \right) ( \rho\circ \gamma^1_t) \, f_1
     \de s \\
     &\le C(t_0)  \int_0^1 ( \rho\circ \gamma^1_t) \, f_1 \de s,
    \end{split}
\end{equation}
where $C(t_0)>0$ is some constant depending on the flow and on the choice of $t_0$. Since $\int_0^1 ( \rho\circ \gamma^1_t) \, f_1 \de s \big|_{t=0} =0$ by definition of $A$, the differential inequality in \eqref{eq:zzEvolutionComparison} implies that $\int_0^1 ( \rho\circ \gamma^1_t) \, f_1 \de s  =0$ for any $t \in [0,t_0)$. This means that $f_1(t,x)=0$ for any $x$ and $t \in [0,t_0)$. By arbitrariness of $t_0$, we get that
\begin{equation*}
    (\widetilde{k}_i)^2 (v_i)^2(t,x) \le  \max_{[0,1]} (\widetilde{k}_1)^2 (v_1)^2\, \big|_{t=0},
\end{equation*}
for any $x$ and $t \in [0,T)$, $i=1,2$. Taking into account \eqref{eq:Step1}, the claimed uniform upper bound on $\widetilde{k}_i$ follows.
The second part of the claim in \ref{Step2} follows by adapting the above reasoning on derivatives $\partial_s^k\widetilde{k}_i$ in place of $\widetilde{k}_i$ or, more easily, by observing that estimates on derivatives $\partial_s^k\widetilde{k}_i$ are independent of the length of $\gamma^0_t$. Indeed, by locality and uniqueness of the flow, the evolution of $\gamma^1_t, \gamma^2_t$ does not change if $\gamma^1_t, \gamma^2_t$ are considered to be edges of a completely analogous network considered in \cref{Fig:Network} except that the length of $\gamma^0_0$ is taken arbitrarily large (see also the discussion in \cref{rem:FamigliaEsempi}). In such a case the upper bound previously proved on the curvature together with lower bounds away from zero on the length of each edge imply uniform bounds on the derivatives $\partial_s^k\widetilde{k}_i$ (independently of ${\rm L}(\gamma^0_t)$) by classical results like \cite[Proposition 5.8]{mannovplusch}.

\item\label{Step3} We want to show that the length of each curve is strictly positive for any time, $T=+\infty$, the length of $\gamma^0_t$ converges to $0$ as $t\to+\infty$, and the curves $\gamma^1_t, \gamma^2_t$ smoothly converge to (half of) the diagonals of the rectangle having vertices at the endpoints of the network, up to reparametrization.

\medskip
By \ref{Step1}, we can parametrize $\gamma^1_t$ as the graph of a function $u:[0,T)\times [0,1]\to \R$, as in \cref{Fig:Graph}.

\begin{figure}[H]
	\begin{center}
		\begin{tikzpicture}[scale=4]
		\draw
		(-1,-0.5)to[out=5, in=210, looseness=0.4]
		(0,-0.03);
		\path[font=\normalsize]
		(-1,-0.5)node[]{$\bullet$};
		\draw[thick,-{latex}]
		(-1.1,-0.5) -- (0.3,-0.5);
		\draw[thick,-{latex}]
		(-1,-0.6) -- (-1,0.3);
		\draw[dashed]
		(0,-0.5)--(0,0.2);
		\draw[dashed]
		(-1,0)--(0.15,0);
		\draw[dotted]
		(-1,-0.5)--(0,0);
		\path[font=\normalsize]
		(0,-0.5)node[below]{$1$};
		\path[font=\normalsize]
		(-0.5,-0.32)node[below]{$u(t,x)$};
		\path[font=\normalsize]
		(-1,0)node[left]{$\frac{1}{\sqrt{3}}$};
		\end{tikzpicture}
	\end{center}
	\caption{Continuous line: graph parametrization of an edge of $\Gamma_t$. Dotted line: straight limit curve of the flow}\label{Fig:Graph}
\end{figure}
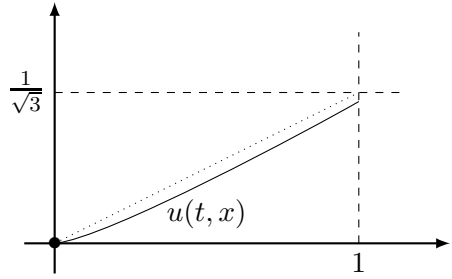
Parametrizing as a graph as in \cref{Fig:Graph} the evolution of an edge $E^i$, whose parametrization evolves according to $\partial\gamma^i_t = \partial^2_x\gamma^i/|\partial_x\gamma^i|^2$, the function $u$ solves the problem
\begin{equation*}
    \begin{cases}
    \partial_t u = \frac{\partial^2_xu}{1+(\partial_x u)^2} & \text{for } (t,x) \in [0,T)\times [0,1], \\
    u(t,0)=0 , \\
    \partial_x u(t,1) = \tan(\pi/6) = 1/\sqrt{3}, \\
    u(0,x)=u_0(x).
    \end{cases}
\end{equation*}
By the above steps, $\partial^2_xu \ge 0$ and $0\le \partial_x u \le \partial_x u(t,1)= 1/ \sqrt{3}$, for any $t \in [0,T)$.

We compare the evolution of $u$ with upper and lower barriers given by solutions of heat-type equations. More precisely, as $\partial^2_xu\ge0$ by convexity, we have that
\begin{equation*}
   \frac{3}{4}\partial^2_xu \le \partial_t u \le \partial^2_xu,
\end{equation*}
at any time and point. Hence we define $v,w$ the solutions to the problems
\begin{equation*}
    \begin{cases}
     \partial_t v = \frac{3}{4}\partial^2_xv & \text{on } [0,+\infty)\times [0,1], \\
    v(t,0)=0 , \\
    \partial_x v(t,1) = 1/\sqrt{3}, \\
    v(0,x)=u_0(x).
    \end{cases}
    \qquad
    \begin{cases}
    \partial_t w= \partial^2_xw &  \text{on } [0,+\infty)\times [0,1], \\
    w(t,0)=0 , \\
    \partial_x w(t,1) = 1/\sqrt{3}, \\
    w(0,x)=u_0(x).
    \end{cases}
\end{equation*}
It is well known that $v,w$ exist for any time and converge to the function $u_\infty(x) \eqdef x/\sqrt{3}$ with an exponential rate in infinite time.

We can consider the function
\begin{equation*}
z(t,x)\eqdef \begin{cases}
u(t,x)-w(t,x) & x \in [0,1], \\
u(t,2-x)-w(t,2-x) & x \in (1,2],
\end{cases}
\end{equation*}
which is the even reflection of the function $u-w$ about the point $x=1$. Hence $z$ is of class $C^2$ and solves
\begin{equation*}
\begin{cases}
\partial_t z \le \partial^2_x u - \partial^2_x w = \partial^2_x z &  \text{on } [0,+\infty)\times [0,2], \\
    z(t,0)=z(t,2)=0 , \\
    z(0,x)=0 & \forall\, x \in [0,2].
\end{cases}
\end{equation*}
By the maximum principle, see \cite[Theorem 2.1.1, Lemma 2.1.3]{MantegazzaBook}, we get that $z\le0$ at any time and point, that is $u(t,x)\le w(t,x)$.

By analogous comparison with $v$, we deduce that $v(t,x)\le u(t,x)\le w(t,x)$. Therefore the length of $\gamma^0_t$ is strictly positive for any $t \in[0,T)$, which, together with \ref{Step2} and \cref{thm:LongTimeGeneral}, implies $T=+\infty$. Moreover, the above comparison analysis completes the proof of \ref{Step3}.
\end{enumerate}

\appendix
\section{Tools needed in some proofs}\label{sec:Appendix}

\subsection*{Quantitative implicit function theorem}

For the convenience of the reader, we sketch the proof of a quantitative implicit function theorem. Specifically, a lower bound on the width of the domain of the implicit function $f$ is given in terms of bounds on the norms of the derivatives of the starting map $F$. The proof is a simplified finite-dimensional version of the general \cite[Theorem 15.8]{Deimling}. An analogous argument can be found in unpublished lecture notes by C. Liverani.

\begin{teo}\label{thm:ImplicitFunction}
Let $n,m\in \N$, $n,m\ge1$, and $(x_0,y_0)\in \R^n\times \R^m$. Denote $Q^n_r\eqdef \{x \in \R^n \st |x-x_0|<r \}$ and $Q^m_r\eqdef \{y \in \R^m \st |y-y_0|<r \}$, for any $r>0$.

Let $F:U\to \R^m$ be a $C^1$ function, where $U\subset \R^n\times\R^m$ is a neighborhood of $(x_0,y_0)$, and assume that $F(x_0,y_0)=0$. Suppose that

\begin{itemize}
    \item $\partial_yF(x_0,y_0)$ is invertible, and let $S\eqdef \|[\partial_yF(x_0,y_0)]^{-1}\|$;
    
    \item $\rho>0$ is such that $\|{\rm id} - [\partial_yF(x_0,y_0)]^{-1}\partial_yF(x,y)\| \le \tfrac12$ for $(x,y) \in Q^n_\rho\times Q^m_\rho$ and $\overline{Q}^n_\rho\times \overline{Q}^m_\rho \Subset U$.
\end{itemize}

Hence, denoting $N\eqdef \sup \{ \|\partial_x F(x,y)\| \st (x,y) \in Q^n_\rho\times Q^m_\rho \}$, there exists $r=r(\rho,S,N) \in (0,\rho]$ such that there exists a unique function $f:Q^n_r\to Q^m_\rho$ such that $f(x_0)=y_0$ and $F(x,y)=0$ for $x \in Q^n_r$ if and only if $y=f(x)$.

%Moreover, $f$ is of class $C^1$ and $\partial_x f(x) = -[\partial_yF(x,f(x))]^{-1}\partial_xF(x,f(x))$. If also $F\in C^k$ for $k \in \N$, then $f \in C^k$.
\end{teo}

\begin{proof}
We just prove that the radius $r$ for the domain $Q^n_r$ of the implicit function $f$ can be chosen depending only on $\rho,S,N$.

Without loss of generality, let $(x_0,y_0)=(0,0)$. Let $r=\min\{\rho/(2SN) ,\rho/2\}$ . For any $x\in \overline{Q}^n_r$ consider the function $\phi_x:\overline{Q}^m_\rho\to \R^m$ given by $\phi_x(y)\eqdef y -[\partial_yF(x_0,y_0)]^{-1} F(x,y) $. We observe that
\[
\|\partial_y \phi_x \| = \|{\rm id} - [\partial_yF(x_0,y_0)]^{-1}\partial_yF(x,y)\| \le \frac12,
\]
for any $y \in Q^m_\rho$. Hence
\[
\begin{split}
|\phi_x(y)| 
&\le |\phi_x(y)-\phi_x(0)|+ |\phi_x(0)| 
\le  \frac12\left|y \right|+ \left|[\partial_yF(x_0,y_0)]^{-1} F(x,0)  \right|\\
&\le \frac{\rho}{2} + S \left|F(x,0) - F(0,0) \right| 
\le \frac{\rho}{2} + S N r \le \rho,
\end{split}
\]
for any $y \in \overline{Q}^m_\rho$. Therefore $\phi_x:\overline{Q}^m_\rho\to\overline{Q}^m_\rho$ is a contraction, and thus there exists a unique $y \in \overline{Q}^m_\rho$ such that $\phi_x(y)=y$, that is, $F(x,y)=0$. Defining $f(x)=y$ the unique $y$ such that $\phi_x(y)=y$, we see that $f$ is defined on $Q^n_r$, and the claim follows.
\end{proof}

\subsection*{A monotonicity-type formula}

We derive here an evolution formula in the spirit of the celebrated Huisken Monotonicity Formula~\cite[Theorem 3.1]{HuiskenMonotonicity}, see also~\cite[Section 1]{EckerHuisken89} 
and~\cite[Theorem 3.1.5, Exercise 3.1.6]{MantegazzaBook}.

Let $\rho:[0,T)\times\R^2\to\R$ be defined by
\begin{equation}\label{eq:DefRho}
    \rho(t,p) \eqdef \frac{1}{\sqrt{4\pi (t_0-t)}} \exp \left(-\frac{|p-p_0|^2}{4(t_0-t)} \right),
\end{equation}
which satisfies
\begin{equation*}
    \nabla\rho = -\frac{p-p_0}{2(t_0-t)}\rho,
    \qquad
    \nabla^2 \rho = -\frac{\rho}{2(t_0-t)}{\rm Id} + \frac{\rho}{4(t_0-t)^2}\, (p-p_0)\otimes (p-p_0).
\end{equation*}
In particular $\partial_t \rho = - \Delta\rho - \frac{\rho}{2(t_0-t)}$.

Moreover consider a smooth evolution of an immersed curve $\gamma:[0,T)\times[0,1]\to\R^2$ by
\begin{equation}\label{eq:zzapp}
    \partial_t\gamma = \boldsymbol{k} + \lambda \tau.
\end{equation}

\begin{lemma}\label{thm:Monotonicity}
Let $\rho$ and $\gamma$ be as in \eqref{eq:DefRho} and \eqref{eq:zzapp}, for some $p_0 \in \R^2$ and $t_0>0$. Let $f:[0,T)\times[0,1]\to \R$ be a function such that $f(t,\cdot) \in H^2(0,1)$ for any $t$ and differentiable with respect to $t$ with $\partial_t f$ continuous, with $T\ge t_0$. Denoting $\rho\circ\gamma\eqdef \rho(t,\gamma(t,x))$, then
\begin{equation*}
\begin{split}
        \frac{\d}{\d t}  \int_0^1 ( \rho\circ \gamma) \, f \de s
    &= \int_0^1 ( \rho\circ \gamma)(\partial_t-\partial_s^2)f -  \left|\boldsymbol{k}+ \frac12 \frac{(\gamma-p_0)^\perp}{t_0-t} \right|^2  ( \rho\circ \gamma) \, f \de s + \\
    &\qquad +
    \int_0^1
      \left( \partial_s \lambda
     -\frac{\lambda}{2(t_0-t)}\scal{\gamma-p_0,\tau}  \right) ( \rho\circ \gamma) \, f
     \de s
    + \\
    &\qquad+\big( ( \rho\circ \gamma)\partial_sf - f \partial_s ( \rho\circ \gamma)\big)\bigg|_0^1,
\end{split}
\end{equation*}
for any $t\in [0,t_0)$.
\end{lemma}

\begin{proof}
If $p =\gamma_t(x)$ then
\begin{equation*}
    \nabla^2\rho \big|_{(t,p)} (\nu_t(x),\nu_t(x)) = - \frac{\rho}{2(t_0-t)} + \frac{\rho}{4(t_0-t)^2}|(p-p_0)^\perp|^2,
\end{equation*}
where $(\cdot)^\perp$ denotes projection along $\nu_t(x)$. Since $\partial_t \rho = - \Delta\rho - \frac{\rho}{2(t_0-t)}$, recalling the relation between Euclidean and intrinsic Laplacian on a submanifold \cite[Lemma 3.1.2]{MantegazzaBook}, we have
\begin{equation*}
\begin{split}
    \frac{\d}{\d t} (\rho \circ \gamma) &=
    -\Delta \rho \big|_{(t,\gamma)} - \frac{\rho \circ \gamma}{2(t_0-t)} + \big\langle\nabla\rho\big|_{(t,\gamma)} , \boldsymbol{k} + \lambda \tau \big\rangle  \\
    &=
    -\partial^2_s(\rho\circ \gamma) - \nabla^2\rho \big|_{(t,\gamma)} (\nu_t,\nu_t) +  \big\langle\nabla\rho\big|_{(t,\gamma)} , \boldsymbol{k} \big\rangle - \frac{\rho \circ \gamma}{2(t_0-t)} + \big\langle\nabla\rho\big|_{(t,\gamma)} , \boldsymbol{k} + \lambda \tau \big\rangle  
    \\&=- \partial^2_s (\rho \circ \gamma) +  \left( -\frac{\scal{\gamma-p_0, \boldsymbol{k}}}{t_0-t}- \frac14 \frac{|(\gamma-p_0)^\perp|^2}{(t_0-t)^2}\right) \rho\circ \gamma - \frac{\lambda}{2(t_0-t)}\scal{\gamma-p_0, \tau} \rho\circ \gamma.
\end{split}
\end{equation*}
Recalling that $\partial_t(\d s) = (\partial_s\lambda - |\boldsymbol{k}|^2) \de s$, the desired formula follows by directly taking the derivative with respect to $t$ and integrating by parts twice in order to transfer the derivatives with respect to $s$ from $\rho$ to $f$.
\end{proof}

\begin{rem}
Lemma~\ref{thm:Monotonicity} is a generalization 
of~\cite[Equation (7), page 455]{EckerHuisken89}
to the case of a curve with boundary and evolving with a tangential velocity $\lambda$ different from zero. 
We refer also to~\cite[Lemma 6.3]{MantegazzaNovagaTortorelli} for the case in which 
$f\equiv 1$.
\end{rem}

\section{Minimal networks on surfaces}
\label{sec:AppendixSurfaces}

In this section we discuss to what extent the theory developed in this work can be adapted to the case of networks in surfaces.

We consider $2$-dimensional complete Riemannian manifolds without boundary, denoted by $(\Sigma, g)$. The obvious adaptation of definitions given in \cref{sec:Networks} allows to speak of networks $\Gamma:G\to \Sigma$, as well as of motion by curvature
(for short time existence see for instance~\cite[Section 8.4]{liramazplusae}). A minimal network $\Gamma_*:G\to \Sigma$
is a collection of geodesic arcs
in $\Sigma$ meeting at triple junctions
forming equal angles.
%shall be a network such that the parametrization of every edge is a geodesic in $\Sigma$ and such that edges meeting at triple junctions form three angles equal to $\tfrac23\pi$.

We expect that minor technical modifications 
of our 
%to the 
arguments 
%above 
lead to the validity of a \L ojasiewicz--Simon inequality as in \cref{thm:LojaIntro} for minimal networks on any analytic surface $(\Sigma, g)$. If such a minimal network is also a local minimizer for the length functional with respect to perturbations sufficiently small in $H^2$ which do not move endpoints, then the stability result as in \cref{thm:StabilityIntro} holds. In particular, relevant examples of Riemannian surfaces where also stability can be proved are given by simply connected analytic surfaces with non-positive sectional curvature, such as the $2$-dimensional hyperbolic space or complete minimal immersions of $\R^2$ in $\R^n$.

\medskip

Let us now describe the main modifications that one should carry out %in order 
to deduce the previous claims.

Concerning preliminary results, suitably adapting the arguments from \cite{GMP}, a short time existence theorem for the flow as \cref{wellposedness} can be proved. Moreover, characterization of singularities as described in \cref{thm:LongTimeGeneral} can be deduced analogously.

As discussed in the work, a graph parametrization like the one established in \cref{sec:GraphParametrization} is \emph{necessary} in order to apply the %most
recent abstract theory 
that implies
%implying 
a \L ojasiewicz--Simon inequality. 

The results of \cref{sec:GraphParametrization} can be directly adapted to networks on $\Sigma$ by employing the exponential map $\exp$ on $(\Sigma, g)$. In fact, let $\Gamma_*, \Gamma:G\to \Sigma$ be a minimal network and a network, respectively, let $m\eqdef \pi(e^i,i)=\pi(e^j,j)=\pi(e^k,k)$ be a junction and denote $m_*\eqdef \Gamma_*(m)$. Assuming that the distance between $\Gamma(m)$ and $m_*$ is less than the injectivity radius ${\rm inj}(m_*)$ of $\Sigma$ at $m_*$, the image $\gamma^\ell(e^\ell)$, for $\ell \in \{i,j,k\}$, can be written as
\[
\gamma^\ell(e^\ell) = \exp_{m_*} \left( \NN^\ell(e^\ell)\nu^\ell_*(e^\ell) + \TT^\ell(e^\ell)\tau^\ell_*(e^\ell) \right),
\]
for some $\NN^\ell(e^\ell), \TT^\ell(e^\ell)$. Hence applying the inverse $\exp_{m_*}^{-1}$ on equalities
\[
\exp_{m_*} \left( \NN^\ell(e^\ell)\nu^\ell_*(e^\ell) + \TT^\ell(e^\ell)\tau^\ell_*(e^\ell) \right)
= 
\exp_{m_*} \left( \NN^s(e^s)\nu^s_*(e^s) + \TT^s(e^s)\tau^s_*(e^s) \right),
\]
for $\ell\neq s$, $\ell,s \in \{i,j,k\}$, readily implies the same linear identities obtained in \cref{lem:BoundaryRelationsTN}. Hence, conversely, inverting such linear relations as in \cref{lem:SufficientConditionsNetwork}, one proves the direct analog of \cref{lem:SufficientConditionsNetwork}. This implies that the linear operators of \cref{def:LinearOperators} can also be used in the setting of networks on $\Sigma$. 

Finally, assuming that parametrizations of $\Gamma$ are sufficiently close in $H^2$ to parametrizations of $\Gamma_*$, meaning that $\sum_i \| \exp_{\gamma^i_*(\cdot)}^{-1}(\gamma^i(\cdot))\|_{H^2}$ is bounded above by a constant also depending on the injectivity radius of $\Gamma_*(G)$, a version of \cref{prop:ParametrizzazioneTN} holds for networks on $\Sigma$. More precisely, curves $\gamma^i$ can be written as
\[
x\mapsto\exp_{\gamma^i_*(x)}( \NN^i(x) \nu^i_*(x) + \TT^i(x)\tau^i_*(x)),
\]
up to reparametrization, where $\TT^i$'s are adapted to the $\NN^i$'s. In order to perform the proof of \cref{prop:ParametrizzazioneTN} on $\Sigma$, it suffices to adapt the argument in neighborhoods of junctions: %and 
this can be carried out by passing into a local chart.

\smallskip

Variations of parametrizations of $\Gamma_*$ analogous to the ones in \cref{prop:FirstVariation} and in \cref{prop:SecondVariationBilinear} take the forms
\begin{align*}
    \gamma^{i,\eps}(x) &= \exp_{\gamma^i_*(x)}( (\NN^i(x) + \eps X^i(x) ) \nu^i_*(x) + \TT^{i,\eps}(x)\tau^i_*(x)), \nonumber\\
\gamma^{i,\eps,\eta}(x) &= \exp_{\gamma^i_*(x)}(  (\eps X^i(x) + \eta Z^i(x))  \nu^i_*(x) + \TT^{i,\eps,\eta}(x)\tau^i_*(x)),
\end{align*}
respectively. Carrying out computations for first and second variations, see \cite[Chapter 9]{DoCarmo} and \cite[Theorem 10.22, Proposition 10.24]{LeeBook}, one obtains the same formulae given in \cref{prop:FirstVariation} and in \cref{prop:SecondVariationBilinear}, except that now the second variation formula \eqref{eq:SecondVariationBilinear} also contains the additive term
\[
\sum_i \int_0^1 -K(\gamma^i_*) \, X^i\, Z^i\, |\partial_x \gamma^i_*| \de x,
\]
where $K(p)$ is the sectional curvature of $\Sigma$ at $p$. However, the linear operator
\[
[H^2(0,1)]^N \ni (X^1,\ldots, X^N) \mapsto \left( -|\partial_x \gamma^1_*| \, K(\gamma^1_*) \, X^1\, , \ldots,  - |\partial_x \gamma^N_*| \,K(\gamma^N_*) \, X^N\, \right) \in [L^2(0,1)]^N,
\]
is compact. Therefore the second variation operator for the length functional differs by a compact operator from the one considered in \cref{sec:Loja}. Since Fredholmness is stable under compact perturbations, i.e., a linear operator $T$ between Banach spaces is Fredholm of index $l$ if and only if $T+T'$ is Fredholm of index $l$, for any compact operator $T'$ (see \cite[Section 19.1]{HormanderIII}), the Fredholmness property required by \cref{prop:LojaAbstract} follows.\\
Assuming that $(\Sigma, g)$ is analytic, all the functional analytic properties required on first and second variations by \cref{prop:LojaAbstract} can be derived as done in \cref{sec:Loja}. Observe that analyticity of the metric $g$ is required as the exponential map shall appear in the expression for the first variation (compare, e.g., with \cite[proposition 3.20]{PozzettaLoja}). Eventually, we deduce that a \L ojasiewicz--Simon inequality as in \cref{thm:LojaIntro} holds for any for minimal network on an analytic surface $(\Sigma, g)$.

\smallskip

Concerning the stability of minimal networks on surfaces, we cannot expect that a version of \cref{thm:StabilityIntro} always holds. Indeed, differently from the case of $\R^2$ (\cref{lemma:MinimalMinimizing}), a minimal network on a surface does not necessarily minimize the length among small perturbations, as this already happens for geodesics.\\ However, assuming that that a minimal network $\Gamma_*:G\to \Sigma$ locally minimizes the length with respect to perturbations having $H^2$-norm sufficiently small and that do not move endpoints, then arguments in the proof of \cref{thm:Stability} can be adapted (see, e.g., \cite[Theorem 4.5]{PozzettaLoja}) to deduce the desired stability. From the technical viewpoint, observe that local minimality with respect to small perturbations is manifestly needed in the argument so that, in the notation of \cref{thm:Stability}, the difference $({\rm L}(\Gamma_t)-{\rm L}(\Gamma_*))$ is non-negative, and thus $H(t)\eqdef ({\rm L}(\Gamma_t)-{\rm L}(\Gamma_*))^\theta$ is well defined.

\smallskip

To conclude, we observe that minimal networks minimize the length among perturbations having $C^0$-norm sufficiently small and which do not move endpoints on simply connected surfaces with non-positive sectional curvature. Indeed, this minimizing property is proved in \cite[Theorem 3.7]{MartelliNovagaPludaRiolo} for surfaces with constant non-positive sectional curvature and it is based on a contradiction argument in combination with the fact that $\delta(t)\eqdef  d(\gamma(t),\sigma(t))$ is convex if $d$ is the geodesic distance on $(\Sigma,g)$ and $\gamma,\sigma$ are minimizing geodesics. The very same argument can be generalized to simply connected surfaces with non-positive sectional curvature taking into account that every geodesic on such a surface is minimizing \cite[Theorem 9.2.2]{BuragoBuragoIvanov} and that convexity for a function $\delta(t)$ as before holds in this generality as well \cite[Lemma 9.2.3]{BuragoBuragoIvanov}.

\bigskip
\noindent\textbf{Conflict of interest.} On behalf of all authors, the corresponding author states that there is no conflict of interest.

\bigskip
\noindent\textbf{Data availability statement.} The manuscript has no associated data.

\bibliographystyle{plain}
\addcontentsline{toc}{section}{References}
\bibliography{references} 

\end{document}